\documentclass{article}
\usepackage[utf8]{inputenc}
\usepackage[margin=1in]{geometry} 
\usepackage{amsmath,amsthm,amssymb,amsfonts}
\usepackage{tikz}
\usetikzlibrary{shapes.geometric, arrows}
\usepackage{verbatim}
\usepackage{dsfont}
\usepackage{cite}
\usepackage{algpseudocode}
\usepackage{algorithm}
\usepackage{subcaption}
\usepackage[colorlinks=true]{hyperref}
\hypersetup{urlcolor=blue, citecolor=blue, linkcolor=blue}
\usepackage[capitalise,noabbrev,nameinlink]{cleveref}
\usepackage[bottom]{footmisc}
\usepackage{multirow}
\usepackage[group-digits=false]{siunitx}
\usepackage{color}

\theoremstyle{plain}
\newtheorem{thm}{Theorem}[section] 
\newtheorem{lem}[thm]{Lemma}
\newtheorem{prop}[thm]{Proposition}

\theoremstyle{definition}
\newtheorem{defn}[thm]{Definition} 
 
\newtheorem{assump}[thm]{Assumption}
\newtheorem{remark}[thm]{Remark}

\numberwithin{equation}{subsection}

\tikzstyle{startstop} = [rectangle, rounded corners, minimum width=3cm, minimum height=1cm,text centered, text width = 4.5cm, draw=black, fill=blue!30]
\tikzstyle{startstop2} = [rectangle, rounded corners, minimum width=3cm, minimum height=1cm,text centered, text width = 4cm, draw=black, fill=red!30]
\tikzstyle{process} = [rectangle, minimum width=3cm, minimum height=1cm, text centered, text width = 8.5cm, draw=black, fill=green!30]
\tikzstyle{decision} = [diamond, minimum width=2cm, minimum height=.5cm, text centered, text width = 1.2cm, draw=black, fill=yellow!30]
\tikzstyle{arrow} = [thick,->,>=stealth]
\tikzstyle{io} = [trapezium, trapezium left angle=70, trapezium right angle=110, minimum width=3cm, minimum height=.5cm, text centered, text width=5cm, draw=black, fill=gray!30]
\tikzstyle{point} = [circle, minimum width=.8cm, minimum height=.8cm, text centered, draw=black, fill=red!30]
\tikzstyle{line} = [draw, -latex']

\begin{document}

\title{Generalized Pseudospectral Shattering and Inverse-Free Matrix Pencil Diagonalization}
\author{James Demmel\footnote{Department of Mathematics, University of California Berkeley} \and Ioana Dumitriu\footnote{Department of Mathematics, University of California San Diego}
\and Ryan Schneider\footnotemark[1] \footnote{Corresponding author (ryan.schneider@berkeley.edu)}}
\date{}

\maketitle

\begin{abstract}
    We present a randomized, inverse-free algorithm for producing an approximate diagonalization of any $n \times n$ matrix pencil $(A,B)$. The bulk of the algorithm rests on a randomized divide-and-conquer eigensolver for the generalized eigenvalue problem originally proposed by Ballard, Demmel, and Dumitriu [Technical Report 2010]. We demonstrate that this divide-and-conquer approach can be formulated to succeed with high probability provided the input pencil is sufficiently well-behaved, which is accomplished by generalizing the recent pseudospectral shattering work of Banks, Garza-Vargas, Kulkarni, and Srivastava [Foundations of Computational Mathematics 2022].  In particular, we show that perturbing and scaling $(A,B)$ regularizes its pseudospectra, allowing divide-and-conquer to run over a simple random grid and in turn producing an accurate diagonalization of $(A,B)$ in the backward error sense. The main result of the paper states the existence of a randomized algorithm that with high probability (and in exact arithmetic) produces invertible $S,T$ and diagonal $D$ such that $||A - SDT^{-1}||_2 \leq \varepsilon$ and $||B - ST^{-1}||_2 \leq \varepsilon$ 
    in at most $O \left(\log^2 \left( \frac{n}{\varepsilon} \right) T_{\text{MM}}(n) \right)$ operations, where $T_{\text{MM}}(n)$ is the asymptotic complexity of matrix multiplication. This not only provides a new set of guarantees for highly parallel generalized eigenvalue solvers but also establishes nearly matrix multiplication time as an upper bound on the complexity of inverse-free, exact-arithmetic matrix pencil diagonalization. 
\end{abstract}

\small
\hspace{2.5mm} \textbf{Keywords:} Matrix pencil, generalized eigenvalue problem, diagonalization, pseudospectra \\
\vspace{2mm}
\indent \hspace{1.5mm} \textbf{AMS MSC2020 Codes:} 15A22, 65F15 \\
\vspace{2mm}
\indent \hspace{1.5mm} Communicated by Peter Buergisser.
\normalsize
\pagebreak

\tableofcontents

\pagebreak

\section{Introduction}
\subsection{The Generalized Eigenvalue Problem}
\label{subsection: the GEP}
Given matrices $A,B \in {\mathbb C}^{n \times n}$ with $\det(A - xB) \in {\mathbb C}[x]$ not identically zero, the generalized eigenvalue problem seeks pairs $(\lambda, v) \in {\mathbb C} \times {\mathbb C}^n$ with $v \neq 0$ such that
\begin{equation}
Av = \lambda Bv .
\label{eqn: generalized eigenvalue problem}
\end{equation}
As the name implies, this is a generalization of the standard, single-matrix eigenvalue problem, which we recover by setting $B = I$. Solutions to \eqref{eqn: generalized eigenvalue problem} yield a generalized eigenvalue $\lambda$ and a generalized eigenvector $v$ of the matrix pencil $(A,B)$. As in the standard eigenvalue problem, generalized eigenvalues and eigenvectors arise in a variety of research disciplines, ranging from signal processing \cite{ESPRIT,SIGNAL_PENCIL} and binary classification \cite{SVM, GEP_Class} to linear differential equations \cite{ODE_book}. \\
\indent Unlike the single-matrix case, however, the generalized eigenvalue problem may admit eigenvalues at infinity, which arise when $B$ is singular. Even more strange behavior occurs when the characteristic polynomial $\det(A-xB)$ \textit{is} identically zero, in which case the pencil $(A,B)$ is \textit{singular} and \eqref{eqn: generalized eigenvalue problem} cannot be used to define eigenvalues and eigenvectors. While the singular case is not the primary focus of this paper, our results apply generally, and we will return to the challenge of computing eigenvalues/eigenvectors for singular pencils later on. \\
\indent When $B$ is invertible, the pencil $(A,B)$ and the matrix $B^{-1}A$ have the same set of eigenvalues and eigenvectors. This prompts a natural approach to the problem:\ explicitly form $B^{-1}A$ and apply any method for solving the single-matrix eigenvalue problem. Though straightforward, this naive approach is unsatisfactory for a variety of reasons. First, many applications feature pencils where one (or both) of the matrices involved is singular. Even when this is not the case, numerical stability concerns still prompt us to avoid taking inverses (see \cref{subsection: Numerical Stability}). With this in mind, our focus will be on methods that work with the pencil $(A,B)$ rather than the matrix $B^{-1}A$. \\
\indent The standard method for solving the generalized eigenvalue problem is the QZ algorithm of Moler and Stewart \cite{QZ}, which obtains generalized eigenvalues and eigenvectors of $(A,B)$ by first producing its generalized Schur form via shifted QR. A significant amount of research has focused on exploring the performance of this approach as well as developing a variety of extensions \cite{QZ1,QZ2,QZ3}. The main drawback to QZ is its relative difficulty to parallelize. Efforts to do so \cite{Parallel_QZ} require a non-trivial reworking of the numerical details and to our knowledge are not widely used.  Methods for extracting only a certain subset of generalized eigenvalues/eigenvectors for sparse problems include the trace minimization algorithm \cite{TRACE_MIN}, projection methods \cite{projection}, and extensions of the Lanczos procedure \cite{Lanczos}. \\
\indent In this paper, we present a naturally parallelizable method that, given any matrix pencil $(A,B)$, constructs a full set of approximate generalized eigenvalues and eigenvectors. For simplicity, we will drop the label ``generalized" in the subsequent sections and simply refer to eigenvalues and eigenvectors of $(A,B)$.
\subsection{Divide-and-Conquer Eigensolvers}
To construct an easily parallelizable method, we employ a divide-and-conquer strategy that recursively splits the original problem $(A,B)$ into smaller ones with disjoint spectra. By ensuring that these subproblems share no eigenvalues (and have no overlapping eigenspaces), we can handle them in parallel. Moreover, as long as the spectrum of each subproblem is contained in its predecessor's, we can reconstruct a full set of eigenvalues/eigenvectors for $(A,B)$ from those of the smallest subproblems, which are either $1 \times 1$ and therefore trivial or small enough to handle with existing techniques. \\
\indent One step of divide-and-conquer consists of the following:
\begin{enumerate}
\item Divide the eigenvalues into two disjoint sets (opposite sides of a line, inside/outside a circle, etc.).
\item Compute unitary projectors onto the corresponding eigenspaces.
\item Apply the projectors to the original problem, splitting it in two.
\end{enumerate}
Note that if we are able to consistently split the eigenvalues into sets of roughly equal size, only $O(\log(n))$ steps of divide-and-conquer are required to find a full set of eigenvalues/eigenvectors. \\
\indent Parallel divide-and-conquer eigensolvers have a rich history in the literature. Early approaches for dense matrices followed work of Beavers and Denman \cite{BEAVERS1974143}, which demonstrated that the matrix sign function of Roberts \cite{SGN} could be used to compute projectors onto certain eigenspaces. While it is possible to adapt the sign function for the generalized eigenvalue problem \cite{Penicl_sgn}, it traditionally requires taking inverses which, as mentioned above, we'd prefer to avoid. An alternative, inverse-free method was outlined by Malyshev \cite{MALYSHEV} and subsequently explored, for both pencils and individual matrices, in work of Bai, Demmel, and Gu \cite{Bai:CSD-94-793}. In this approach, spectral projectors are approximated by repeatedly and implicitly squaring $A^{-1}B$, where eigenvalues are split by the unit circle. \\ 
\indent The requirement to have unitary projectors in each step of divide-and-conquer is non-trivial. In fact, neither the sign function nor the inverse-free, repeated squaring approach of Malyshev produce explicitly unitary projectors; instead, rank-revealing factorizations are applied to guarantee a unitary result. The rank-revealing aspect is important here, as the rank of each projector tells us how many eigenvalues lie in a corresponding region of the complex plane. While there are many ways to compute a rank-revealing factorization \cite{RRQR1,ULV,RRQR2}, we are specifically interested in a randomized, rank-revealing URV factorization introduced by Demmel, Dumitriu, and Holtz \cite{2007}, which is simple to implement and particularly compatible with inverse-free computations \cite{Ballard2010MinimizingCF}. \\
\indent The method we present combines this randomized, rank-revealing factorization with the inverse-free approach of Malyshev. A technical report of Ballard, Demmel, and Dumitriu \cite{Ballard2010MinimizingCF} already explored this combination, though it left one fundamental question unanswered:\ how do we reliably split the spectrum at each step? The importance of this cannot be understated, as divide-and-conquer can only hope to be efficient if the split is significant (i.e., close to 50/50) at each step. And while we might hope that randomly selecting a dividing line/circle may work with high probability, the method will fail if whatever we choose intersects the spectrum or even pseudospectrum. Our most important contribution, we resolve this problem by adapting recent work of Banks et al.\ \cite{banks2020pseudospectral}, demonstrating that by adding a bit of randomness to the problem we can state the inverse-free procedure in a way that succeeds with high probability.

\subsection{Pseudospectra and Random Matrix Theory}
The pseudospectrum of a matrix or matrix pencil consists of all eigenvalues of nearby matrices/pencils. We will define them precisely in \cref{section:pseudospectra}, but for now we can think of the $\epsilon$-pseudospectrum as containing the eigenvalues of matrices or pencils that are, in a spectral norm sense, $\epsilon$-away. In general, each pseudospectrum consists of connected components in ${\mathbb C}$ containing at least one true eigenvalue, to which they collapse as $\epsilon \rightarrow 0$. For more background, see the standard reference \cite{TrefethenEmbree+2020}. \\
\indent Looking toward floating-point implementations where error may accumulate at each step, we expect that divide-and-conquer eigensolvers are most likely to fail when the $\epsilon$-pseudospectra for small $\epsilon$ (say close to machine precision) cover large areas of the complex plane, in which case it becomes difficult to select a dividing line/circle that will significantly split the spectrum without causing the process to break down. This is particularly important for matrix pencils, since infinite eigenvalues are a possibility. \\
\indent To work around this, our algorithm begins by adding small perturbations to the input matrices, swapping the pencil $(A,B)$ for $(\widetilde{A}, \widetilde{B}) = (A + \gamma G_1, B + \gamma G_2)$, where $\gamma$ is a tuning parameter and $G_1$ and $G_2$ are independent, complex Ginibre random matrices -- i.e., matrices whose entries are independent and identically distributed (i.i.d.)\ standard complex Gaussian with mean zero and variance $1/n$. It is well-known that random matrices with i.i.d.\ entries have a number of useful properties \cite{Ed88, circular, edelman_rao_2005}; not only are they almost surely invertible, but in many cases distributions for their smallest singular values are known \cite{RUDELSON2008600, Tao_Vu}. The perturbation applied to $(A,B)$ leverages these properties to regularize the problem, replacing what could be a poorly-conditioned input pencil with one whose eigenvalue problem is well-behaved (in a way to be defined precisely later).  \\
\indent The concept of using a random perturbation to regularize a matrix -- considered a linear algebra extension of smoothed analysis -- has been well-studied \cite{gaussian_reg1,gaussian_reg2,gaussian_reg3}. Of particular relevance here, Banks et al.\ \cite{banks2020gaussian} showed that adding a Ginibre perturbation to a matrix provides a probabilistic bound on its smallest singular value while also controlling its eigenvector condition number. In a subsequent work \cite{banks2020pseudospectral} they demonstrated that these perturbations also regularize $\epsilon$-pseudospectra, breaking them into small connected components around distinct eigenvalues. Banks et al.\ used this phenomena, referred to as \textit{pseudospectral shattering}, as the basis for a divide-and-conquer eigensolver, which proceeds as follows:
\begin{enumerate}
\item Perturb the input matrix, shattering its pseudospectra into distinct components.
\item Construct a random grid that covers the spectrum of the matrix. Shattering implies that a fine enough grid does not intersect the pseudospectra, and moreover each eigenvalue of the perturbed matrix is contained in a unique grid box with high probability.
\item Run divide-and-conquer (using the sign function), each time splitting the eigenvalues by a grid line. 
\end{enumerate}
This rough outline obscures the numerical details but captures the core idea:\ by perturbing, we obtain a matrix whose pseudospectra are sufficiently well-behaved for divide-and-conquer to succeed using a simple random grid.  We generalize these results here, proving our own version of pseudospectral shattering for matrix pencils and using it to complete the divide-and-conquer algorithm posed by Ballard, Demmel, and Dumitriu \cite{Ballard2010MinimizingCF}. \\
\indent This is not the first approach that solves the generalized eigenvalue problem by perturbing the input matrices. A method of Hochstenbach, Mehl, and Plestenjak \cite{rank_completing_perturb}, for example, finds eigenvalues via QZ by first applying rank-completing perturbations, if applicable. Note that because our perturbations are random, and because we use a random rank-revealing factorization, convergence is only possible with high probability. 

\subsection{Matrix Pencil Diagonalization}
The decision to perturb the input matrices also impacts accuracy, where only backward-error guarantees will be available. In particular, the eigensolver we present approximates the eigenvalues/eigenvectors of the nearby pencil $(\widetilde{A}, \widetilde{B})$ rather than those of $(A,B)$. While perturbation results for the generalized eigenvalue problem are well-developed in the literature, see for example \cite[Chapter VI]{stewart1990matrix}, little can be said about the accuracy of an individual eigenvalue or eigenvector found this way in the absence of additional information (e.g., a condition number). Ultimately, the eigenvalues of $(\widetilde{A}, \widetilde{B})$ may vary significantly from those of $(A,B)$ even when $||A - \widetilde{A}||_2$ and $||B - \widetilde{B}||_2$ are small. Take the following $2 \times 2$ matrices as an example:
\begin{equation}
    A = B = \begin{pmatrix} 1 & 0 \\ 0 & 10^{-10} \end{pmatrix}, \; \; \widetilde{A} = \begin{pmatrix} 1 & 0 \\ 0 & 1.5 \times 10^{-10} \end{pmatrix}, \; \; \widetilde{B} = \begin{pmatrix} 1 & 0 \\ 0 & 5 \times 10^{-11} \end{pmatrix} .
    \label{eqn: two by two example}
\end{equation}
By construction, $||A - \widetilde{A}||_2 = ||B - \widetilde{B}||_2 = 5 \times 10^{-11}$ while $\Lambda(A,B) = \left\{ 1 \right\}$ and $\Lambda(\widetilde{A}, \widetilde{B}) = \left\{ 1, 3 \right\} $ -- i.e., a tiny perturbation in $A$ and $B$ results in macroscopic changes to the eigenvalues. This is a consequence of the instability of the second eigenvalue of $(A,B)$, which itself is due to the near singularity of $B$. Guarantees are even worse when the matrix $B$ is actually singular. In particular, our algorithm always produces a set of finite eigenvalues, even though the input pencil may have one or many eigenvalues at infinity. \\
\indent With this in mind, we consider the related problem of joint diagonalization. A pair of nonsingular matrices $S$ and $T$ jointly diagonalize $A$ and $B$, or equivalently diagonalize the pencil $(A,B)$, if $S^{-1}AT$ and $S^{-1}BT$ are both diagonal. The connection to the generalized eigenvalue problem is fairly straightforward, as the pencils $(A,B)$ and $(S^{-1}AT, S^{-1}BT)$ have the same set of eigenvalues. Moreover any matrix pencil with a full set of distinct eigenvalues can be diagonalized by matrices $S$ and $T$ containing its left and right eigenvectors, respectively \cite[Corollary VI.1.12]{stewart1990matrix}. \\
\indent It is in this context that our method produces provably accurate results, in a backward-error sense, regardless of the input pencil $(A,B)$. The main result of the paper (presented technically as \cref{thm: RPD guarantees} and \cref{prop: op count} below) may be summarized as follows:

\begin{thm}\label{thm: main result}
    Let $(A,B) \in {\mathbb C}^{n \times n} \times {\mathbb C}^{n \times n}$ be any matrix pencil and let $\varepsilon < 1$ be any desired (backward) diagonalization accuracy. There exists an exact-arithmetic, inverse-free, and randomized divide-and-conquer algorithm that takes as inputs $(A,B)$ and $\varepsilon$ and produces matrices $D,S,T \in {\mathbb C}^{n \times n}$ such that $D$ is diagonal, $S$ and $T$ are nonsingular, and 
    $$ ||A - SDT^{-1}||_2 \leq \varepsilon \; \; \text{and} \; \; ||B - ST^{-1}||_2 \leq \varepsilon $$
    with probability at least $1 - O(\frac{1}{n})$. This algorithm requires at most $O(\log^2(\frac{n}{\varepsilon}) T_{\text{\normalfont MM}}(n))$ operations, where $T_{\text{\normalfont MM}}(n)$ is the complexity of $n \times n$ matrix multiplication. 
\end{thm}

Note that \cref{thm: main result} is fully general, and in particular applies to edge cases like $A = B = 0$. As we will see, this result provides forward error guarantees on approximate eigenvalues and eigenvectors when $B$ is invertible and $(A,B)$ has a full set of distinct eigenvalues (\cref{thm: forward error bound}).

\subsection{Complexity and Numerical Stability}
\label{subsection: Numerical Stability}
The complexity $O(\log^2(\frac{n}{\varepsilon}) T_{\text{MM}}(n))$  in \cref{thm: main result} implies that our algorithm runs in \textit{nearly matrix multiplication time}, meaning its asymptotic complexity is the same as matrix multiplication up to logarithmic factors. This is accomplished by showing that (1) one step of divide-and-conquer runs in nearly matrix multiplication time, and (2) that only $O(\log(n))$ steps are needed to obtain the approximate diagonalization. As alluded to earlier, (2) is guaranteed as long as the eigenvalue split at each step is significant. Claim (1), meanwhile, stems from using only QR and matrix multiplication, both of which can be implemented stably in $O(T_{\text{MM}}(n))$ operations (see work of Demmel, Dumitriu, and Holtz \cite{2007}). In fact, matrix multiplication and QR have $O(T_{\text{MM}}(n))$ implementations that are respectively forward and backwards stable -- i.e., the relative forward/backward error is bounded by machine precision ${\bf u}$ times an appropriate condition number. \\
\indent This combination of efficiency and stability is our motivation for stating an inverse-free approach. While it is possible to invert a matrix in $O(T_{\text{MM}}(n))$ operations, such implementations are only logarithmically stable \cite[\S 3]{2007}, meaning a polylogarithmic increase in precision is required to obtain the same error as a truly backward stable version. Accordingly, the pseudospectral method of Banks et al., which requires inversion to compute the sign function, is built on a logarithmically stable black-box inversion algorithm \cite[Definition 2.7]{banks2020pseudospectral}. This lack of stability carries over into their floating-point analysis, which yields a slightly weaker-than-logarithmic stability guarantee for the divide-and-conquer routine as a whole. \\
\indent We do not present a full finite-precision analysis of our algorithm in this paper (though we do outline progress in that direction in \cref{appendix: finite arithmetic}). Nevertheless, we anticipate that the inverse-free approach presented here can produce a logarithmically stable method for computing an approximate diagonalization of any matrix or pencil, including the singular case. In addition, we present examples in \cref{section: shattering,section experiments} in which the use of inversion (in double precision) introduces nontrivial -- and potentially catastrophic -- errors. In light of these examples, we make the case that an inverse-free algorithm is likely to perform significantly better on arbitrary pencils. As a result, we position our work as a more stable generalization of the algorithm presented by Banks et al. \cite{banks2020pseudospectral}. \\
\indent Up to this point, efficiency has been characterized by the number of arithmetic operations performed. We could alternatively focus on communication costs, both between multiple processors and/or between levels of memory hierarchy. This is especially relevant for the algorithm presented here, which is highly parallel and primed for memory-constrained applications; in particular, divide-and-conquer can be used to reduce a problem to fit in available fast memory. With this in mind, we could consider a communication-optimal implementation in the vein of Ballard et al. \cite{minimizing_comm_NLA}. In fact, the variant of our approach presented by Ballard, Demmel, and Dumitriu \cite{Ballard2010MinimizingCF} was shown to achieve asymptotic communication lower bounds for $O(n^3)$ eigensolvers (again, a consequence of avoiding inversion). While we won't analyze in detail the communication costs of our algorithm, we mention this here to emphasize the flexibility of the high level strategy:\ the method may be easily reformulated to minimize communication costs if desired.

\subsection{Outline}
To improve readability, we summarize the main contributions of each of the remaining sections here.
\begin{enumerate}
\item In \cref{section: background}, we define the pseudospectrum of a matrix pencil and explore some of its properties, culminating in our version of Bauer-Fike \cite{BF} for pencil pseudospectra (\cref{thm: pencil BF}). 
\item In \cref{section: shattering}, we prove pseudospectral shattering for matrix pencils (\cref{thm: shattering}). The proof presented is a fairly straightforward generalization of the single-matrix case done by Banks et al. \cite{banks2020pseudospectral}.
\item In \cref{section: divide and conquer}, we state our divide-and-conquer eigensolver  and prove that it succeeds with high probability in exact arithmetic (\cref{thm: EIG succeeds}).
\item In \cref{section diag}, we use the eigensolver to construct an algorithm that produces an accurate diagonalization of any matrix pencil with high probability, thereby proving \cref{thm: main result}.
\item In \cref{section experiments} we present a handful of numerical experiments to demonstrate that the algorithm works in finite precision. 
\item Finally in the appendices, we consider alternative versions of Bauer-Fike for matrix pencils (\cref{appendix: BF}) and discuss progress towards a floating-point analysis of the algorithm (\cref{appendix: finite arithmetic}).
\end{enumerate}

\section{Background}
\label{section: background}
\subsection{Matrix Pencils, Eigenvalues, and Eigenvectors}
In this section, we review necessary background information from linear algebra as presented in the standard references \cite{stewart1990matrix, horn_johnson_2012}. Here and in the subsequent sections, $||\cdot ||_2$ is the spectral norm, $A^H$ and $A^{-H}$ denote the Hermitian transpose and inverse Hermitian transpose of $A$, respectively, and $\kappa(A)$ is the spectral norm condition number of $A$. Throughout, $(A,B)$ is a matrix pencil with $A,B \in {\mathbb C}^{n \times n}$,  $\Lambda(A)$ and $\Lambda(A,B)$ denote the spectrum of $A$ and $(A,B)$, and $\rho(A)$ is the spectral radius. Note that while non-square pencils are of interest in certain applications, we do not consider them here. Finally, we assume familiarity with the basic definitions of the single-matrix eigenvalue problem.
\begin{defn}
    The matrix pencil $(A,B)$ is regular if its characteristic polynomial $\det(A-x B) \in {\mathbb C}[x]$ is not identically zero. If $(A,B)$ is not regular, it is singular.
\end{defn}
In the following definitions we assume that $(A,B)$ is regular. For a discussion of singular pencils, see \cref{rem: singular_pencils}.
\begin{defn} 
    A nonzero vector $v \in {\mathbb C}^n$ is a right eigenvector of $(A,B)$ if $Av = \lambda Bv$ for some $\lambda \in {\mathbb C}$. In this case the eigenvalue $\lambda$ is a root of the characteristic polynomial $\det(A-x B)$ and $(\lambda, v)$ is an eigenpair of $(A,B)$. Similarly, a nonzero vector $w$ is a left eigenvector of $(A,B)$ if $w^HA = \lambda w^HB$. 
    \label{defn: generalized eigenvals/eigenvect}
\end{defn}
Technically, \cref{defn: generalized eigenvals/eigenvect} only defines finite eigenpairs; missing are the infinite eigenvalues that correspond to nontrivial vectors in the null space of $B$, which in some sense formally satisfy $\frac{1}{\lambda} Av = Bv$ (though of course infinity cannot be a root of the characteristic polynomial). Ultimately, eigenvectors with infinite eigenvalues are included for completeness:\ since the characteristic polynomial of $(A,B)$ has degree less that $n$ if $B$ is singular, we only obtain a full set of $n$ eigenvalues for $(A,B)$ if we include those at infinity. As is standard, we will use ``eigenvector" to refer to right eigenvectors unless stated otherwise. \\
\indent Throughout, we are interested in factorizations of $A$ and $B$ that reveal the eigenvalues of $(A,B)$. First, we consider a generalization of Schur form to matrix pencils.
\begin{defn}
    $(A,B) = (UT_AV^H, UT_BV^H)$ is a generalized Schur form of $(A,B)$ if $U,V \in {\mathbb C}^{n \times n}$ are unitary and $T_A$ and $T_B$ are upper triangular.
\end{defn}
\indent $(A,B)$ and its Schur form are equivalent, meaning they have the same set of eigenvalues -- given by $T_{A}(i,i)/T_B(i,i)$ -- where zeros on the diagonal of $T_B$ correspond to eigenvalues at infinity. It is not hard to show that $(A,B)$ is guaranteed to have a Schur form (see for example \cite[Theorem VI.1.9]{stewart1990matrix}). As in the single-matrix case, we can build a generalized Schur decomposition by constructing $V$ so that any leading set of columns is an orthonormal basis for a space spanned by a collection of right eigenvectors of $(A,B)$. In fact, the columns of $V$ span right deflating subspaces for the pencil.
\begin{defn}
    Subspaces ${\mathcal X}, {\mathcal Y} \subseteq {\mathbb C}^n$ are right and left deflating subspaces of $(A,B)$, respectively, if $\dim({\mathcal X}) = \dim({\mathcal Y})$ and 
    $\text{span} \left\{Ax, Bx : x \in {\mathcal X} \right\} = {\mathcal Y} .$
    \label{defn: left and right deflating subspaces}
\end{defn}
The right and left deflating subspaces of $(A,B)$ generalize the invariant eigenspaces of an individual matrix. That is, a collection of eigenvectors of $(A,B)$ span a right deflating subspace that is ``mapped" to a corresponding left deflating subspace by $A$ and $B$. Despite the similarity in naming, the left deflating subspace of $(A,B)$ corresponding to a set of right eigenvectors is not usually spanned by left eigenvectors. \\
\indent Note that if $(A,B) = (UT_AV^H, UT_BV^H)$ then the first $k$ columns of $U$ span the left deflating subspace corresponding to the first $k$ columns of $V$ (this follows from $AV = UT_A$ and $BV = UT_B$ since $T_A$ and $T_B$ cannot both have a zero in the same diagonal entry). Moreover, if we take the first $k$ columns of $U$ as a right deflating subspace of $(A^H,B^H)$ then the first $k$ columns of $V$ span the corresponding left deflating subspace since $A^HU = VT_A^H$ and $B^HU = VT_B^H$. In this way, taking a Hermitian transpose swaps left and right deflating subspaces. \\
\indent Moving beyond Schur form, we next generalize diagonalization to pencils.
\begin{defn}
    Nonsingular matrices $S$ and $T$ diagonalize $(A,B)$ if $S^{-1}AT = D_1$ and $S^{-1}BT = D_2$ are both diagonal.
    \label{defn: pencil diagonalization}
\end{defn}
Once again, $(A,B)$ and its diagonalization $(D_1, D_2)$ have the same spectrum, with infinite eigenvalues corresponding to zeros on the diagonal of $D_2$. Moreover, the columns of $T$ and rows of $S^{-1}$ contain right and left eigenvectors of $(A,B)$. As in the single-matrix problem, not every pencil is diagonalizable. Instead, the diagonalization given in \cref{defn: pencil diagonalization} is a special case of the Kronecker canonical form \cite{kronecker}, which generalizes Jordan form to matrix pencils. Of particular use here, note that if $B$ is invertible and $T$ contains right eigenvectors of $(A,B)$ then $S = BT$ produces a diagonalization with $D_2 = I$. \\
\indent To prove pseudospectral shattering in \cref{section: shattering}, we will need the following definitions related to the eigenvalues/eigenvectors of an individual matrix. While each of these generalizes to pencils in a fairly straightforward way, we state them first for the single-matrix case to match how they are used in the analysis to come. 
\begin{defn}
    The eigenvector condition number of a diagonalizable matrix $M$ is  $\kappa_V(M) = \inf_{V} \kappa(V) ,$ where the infimum is taken over nonsingular matrices $V$ that diagonalize $M$.
    \label{defn: eigenvector condition number}
\end{defn}
\begin{defn}
    Let $\lambda_1(M), \ldots \lambda_n(M)$ be the eigenvalues of $M$ repeated according to multiplicity. The eigenvalue gap of $M$ is $ \text{gap}(A) = \min_{i \neq j} |\lambda_i(M) - \lambda_j(M)|. $
    \label{defn: eigenvalue gap}
\end{defn}
Under the right conditions, we can also define a condition number for each eigenvalue of a diagonalizable matrix.
\begin{defn}
    If $\lambda_i$ is an eigenvalue of a matrix $M$ with distinct eigenvalues and $v_i$ and $w_i$ are right and left eigenvectors corresponding to $\lambda_i$ normalized so that $w_i^Hv_i= 1$, then the condition number of $\lambda_i$ is 
    $$ \kappa(\lambda_i) = ||v_iw_i^H||_2 = ||v_i||_2 ||w_i||_2 .$$
\end{defn}
\indent Together, the eigenvector/eigenvalue condition numbers and the eigenvalue gap give an indication of how well-behaved the corresponding eigenvalue problem is. To generalize this to the eigenvalues/eigenvectors of a matrix pencil, we can adapt $\kappa_V(M)$ and $\text{gap}(M)$ as follows.
\begin{defn}
   \label{defn: pencil eigenvector condition number}
    The (right) eigenvector condition number of a diagonalizable pencil $(A,B)$ is $\kappa_V(A,B) = \inf_T \kappa(T)$, where the infimum is taken over all invertible $T$ containing a full set of right eigenvectors of $(A,B)$.
\end{defn}
\begin{defn}
    Let $\lambda_1(A,B), \ldots, \lambda_n(A,B)$ be the eigenvalues of $(A,B)$ repeated according to multiplicity. The eigenvalue gap of $(A,B)$ is $ \text{gap}(A,B) = \min_{i \neq j} | \lambda_i(A,B) - \lambda_j(A,B)|. $
    \label{defn: pencil eigenvalue gap}
\end{defn}
\indent Of course, these definitions are not the only possible condition numbers for generalized eigenvalues and eigenvectors \cite{Higham_condition_number, ANGUAS2019170}. We justify our choice in \cref{section RPD}, where we show that a backward error diagonalization of a suitably well-behaved pencil provides forward error approximations of its eigenvalues and eigenvectors in terms of $\kappa_V(A,B)$ and $\text{gap}(A,B)$ (\cref{thm: forward error bound}). Note that when $B$ is invertible, $\kappa_V(A,B) = \kappa_V(B^{-1}A)$ and $\text{gap}(A,B) = \text{gap}(B^{-1}A)$. \\
\indent Finally, we consider an important building block of divide-and-conquer eigensolvers:\ spectral projectors. 
\begin{defn}
    $P$ is a spectral projector of the matrix $M$ if $P^2 = P$ and $MP= PM$.
    \label{defn: spectral projector}
\end{defn}
The simplest spectral projector is $vw^H$ for $v$ and $w$ right and left eigenvectors of $M$ corresponding to the same eigenvalue $\lambda$ and scaled so that $w^Hv = 1$. In this case, we have by Cauchy's integral formula
\begin{equation}
    vw^H = \frac{1}{2 \pi i} \oint_{\Gamma}(z - M)^{-1} dz
    \label{eqn: Cauchy + resolvent}
\end{equation}
for $\Gamma$ any closed, rectifiable curve that separates $\lambda$ from the rest of the spectrum of $M$. More generally, any spectral projector can be expressed as
\begin{equation}
    P =  V\begin{pmatrix} I_k & 0 \\ 0 & 0 \end{pmatrix} V^{-1}
    \label{eqn: spectral projector}
\end{equation}
for $V$ a matrix containing the right eigenvectors of $M$, again appropriately scaled. \\
\indent As mentioned above, an important step of any divide-and-conquer eigensolver is the approximation of spectral projectors onto certain sets of eigenvectors. For our purposes, since $(A,B)$ and $B^{-1}A$ have the same eigenvectors if $B$ is nonsingular, it will be sufficient to approximate the spectral projectors of $B^{-1}A$ without explicitly forming it. 

\begin{remark}\label{rem: singular_pencils}
    When $(A,B)$ is singular, $\lambda \in {\mathbb C}$ is a finite eigenvalue if $\text{rank}(A-\lambda B) < \text{rank}(A - x B)$, where the latter is computed over the field of fractions of ${\mathbb C}[x]$. From there, all of the definitions given above can be suitably adapted to singular pencils (see for example \cite{Forney,reducing_subspaces,DOPICO202037}). Because our algorithm begins with a perturbation -- which guarantees almost surely that the perturbed pencil is regular -- a deeper discussion of singular pencils is beyond the scope of this paper. Nevertheless, our main error bound still applies to input pencils that are singular, and in \cref{section experiments} we will explore the practical challenges of recovering eigenvalues in this case. 
\end{remark}

\indent To wrap up this section, we state a few general linear algebra results/facts that will be useful later on. First is the stability of singular values, which follows from Weyl's inequality \cite[Corollary IV.4.9]{stewart1990matrix} 
\begin{lem}[Stability of Singular Values] \label{lem: stability of singular values}
    If $\sigma_i(M)$ is the $i$-th singular value of a matrix $M$, then for any $M_1,M_2 \in {\mathbb C}^{n \times n}$
    $$| \sigma_i(M_1) - \sigma_i(M_2)| \leq ||M_1 - M_2||_2, \; \; \;  1 \leq i \leq n. $$
\end{lem}
Second, we note that the spectral norm is often generalized to pencils by setting $||(A,B)||_2$ equal to the norm of the $n \times 2n$ matrix obtained by concatenating $A$ and $B$. Throughout we use the loose upper bound $||(A,B)||_2 \leq ||A||_2 + ||B||_2$.

\subsection{Pseudospectra}
\label{section:pseudospectra}
In this section, we define pseudospectra for both individual matrices and matrix pencils. We start with the single-matrix case, where the definition is standard.
\begin{defn} \label{defn: matrix pseudospectrum}
For any $\epsilon > 0$, the $\epsilon$-pseudospectrum of $M$ is
$$ \Lambda_{\epsilon}(M) = \left\{ z : \text{there exists} \; u \neq 0 \; \text{with}\;  (M + E)u = zu \; \text{for some} \; ||E||_2 \leq \epsilon \right\}. $$
\end{defn}
Each pseudospectrum $\Lambda_{\epsilon}(M)$ consists of connected components in ${\mathbb C}$ containing at least one eigenvalue of $M$. As $\epsilon \rightarrow 0$ these connected components collapse to the true eigenvalues of the matrix. If $M$ is diagonalizable, we can quantify this explicitly via the Bauer-Fike Theorem \cite{BF}.

\begin{thm}[Bauer-Fike]\label{thm: bauer fike}
If $\lambda_1, \ldots, \lambda_n$ are the eigenvalues of $M \in {\mathbb C}^{n \times n}$ and $V$ is any invertible matrix that diagonalizes $M$, then 
$$ \bigcup_{i=1}^n B_{\epsilon}(\lambda_i) \subseteq \Lambda_{\epsilon}(M) \subseteq \bigcup_{i=1}^n B_{\epsilon \kappa(V)}(\lambda_i) $$
where $B_r(z)$ is the ball of radius $r$ centered at $z$.
\end{thm}
To extend this to matrix pencils, we must first decide how to define $\Lambda_{\epsilon}(A,B)$ -- the $\epsilon$-pseudospectrum of the pencil $(A,B)$. Unlike the single matrix case there is no standard choice here (see Trefethen and Embree \cite[Chapter 45]{TrefethenEmbree+2020} for a discussion of various options). In this paper, we follow a definition originally due to Frayss\'{e} et al.\ \cite{Fraysse96spectralportraits}.

\begin{defn}\label{defn:pseudospectrum}
For any $\epsilon > 0$, the $\epsilon$-pseudospectrum of $(A,B)$ is
$$ \Lambda_{\epsilon}(A,B) = \left\{ z :  \text{there exists}  \; u \neq 0 \; \text{with}  \; (A + \Delta A)u = z(B + \Delta B)u \; \text{for some} \; ||\Delta A||_2, ||\Delta B||_2 \leq \epsilon \right\}  .$$
\end{defn}

\indent There are several equivalent ways to state \cref{defn:pseudospectrum}, which we use throughout when convenient.

\begin{figure}[t]
    \centering
    \begin{subfigure}{.5\linewidth}
        \centering
        \includegraphics[width=.95\linewidth]{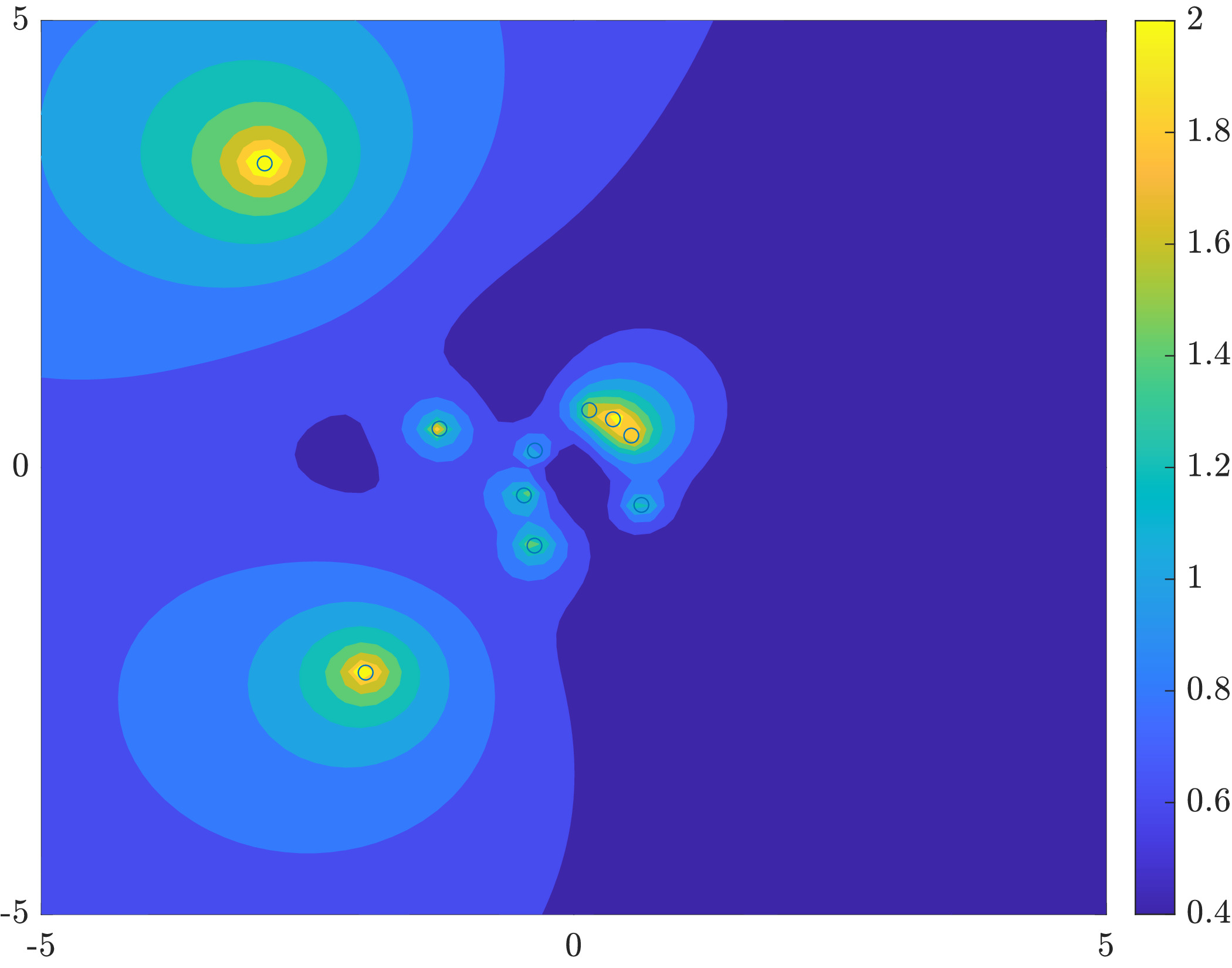}
        \caption{Pseudospectra of $(A,B)$}
    \end{subfigure}%
    \begin{subfigure}{.5\linewidth}
        \centering
        \includegraphics[width=.95\linewidth]{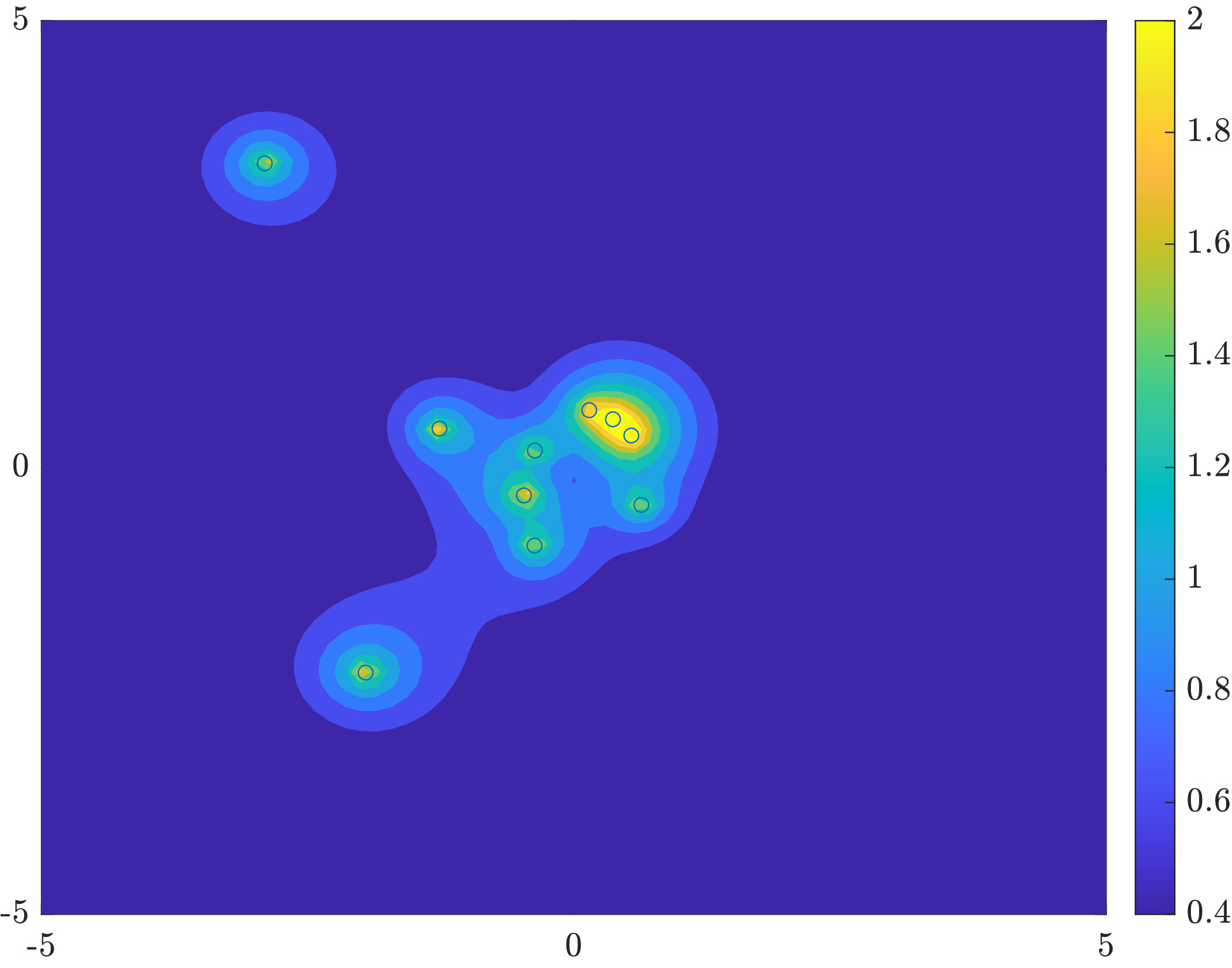}
        \caption{Pseudospectra of $B^{-1}A$}
    \end{subfigure}
\caption{Pseudospectra of $(A,B)$ and $B^{-1}A$ for Gaussian $A,B \in {\mathbb C}^{10 \times 10}$ following \cref{defn:pseudospectrum} and \cref{defn: matrix pseudospectrum}, respectively. Eigenvalues are plotted with open circles. Pseudospectra were obtained by graphing the level curves of$\log_{10}\left[(1 + |z|)||(A - zB)^{-1}||_2\right]$ and $\log_{10}\left[||(B^{-1}A - zI)^{-1}||_2\right]$ in Matlab R2023a.}
    \label{fig: pseudospectra comparison}
\end{figure}

\begin{thm}[Frayss\'{e} et al.\ 1996]\label{thm: equiv pseudo}
The following are equivalent:
\begin{enumerate}
    \item $z \in \Lambda_{\epsilon}(A,B)$.
    \item There exists a unit vector $u$ such that $||(A-zB)u||_2 \leq \epsilon(1 + |z|)$.
    \item $||(A - zB)^{-1}||_2  \geq \frac{1}{\epsilon(1 + |z|)}$.
    \item $\sigma_n(A - zB) \leq \epsilon (1 + |z|)$.
\end{enumerate}
\end{thm}

\indent \cref{fig: pseudospectra comparison} plots the pseudospectra of $(A,B)$ and $B^{-1}A$ for one (randomly chosen) pair $A,B \in {\mathbb C}^{10 \times 10}$. Clearly, the pseudospectra of $(A,B)$ and $B^{-1}A$ differ significantly. Informally, we might say that the pseudospectra of $(A,B)$ are less well-behaved than those of $B^{-1}A$, despite the fact that both coalesce around the same set of eigenvalues.  In the remainder of this section, we pursue bounds on $\Lambda_{\epsilon}(A,B)$.

\begin{lem}\label{lem: bounded pseudo}
 $\Lambda_{\epsilon}(A,B)$ is bounded if and only if $\epsilon  < \sigma_n(B)$.
 \end{lem}
\begin{proof}
Let $\epsilon < \sigma_n(B)$ and suppose $\Lambda_{\epsilon}(A,B)$ is unbounded. By \cref{thm: equiv pseudo}, any nonzero $z \in \Lambda_{\epsilon}(A,B)$ satisfies
\begin{equation}
    \frac{1}{\epsilon(1 + |z|)} \leq ||(A - zB)^{-1}||_2 = \frac{1}{|z|} \left| \left| \left(\frac{1}{z}A - B \right)^{-1} \right| \right|_2 \;\;  \Longrightarrow \; \; \frac{|z|}{\epsilon ( 1+ |z|)} \leq \left| \left| \left(\frac{1}{z}A - B \right)^{-1} \right| \right|_2.
    \label{eqn: apply equiv pseudo def}
\end{equation}
Since $\Lambda_{\epsilon}(A,B)$ is unbounded, we can take a limit $z \rightarrow \infty$ in this inequality to obtain
\begin{equation}
    \frac{1}{\epsilon} \leq ||(-B)^{-1}||_2 = \frac{1}{\sigma_n(B)},
    \label{eqn: take limit to infinity}
\end{equation}
which implies $\epsilon \geq \sigma_n(B)$, a contradiction. \\
\indent To prove the converse, suppose now $\epsilon \geq \sigma_n(B)$ and let $U \Sigma V^H$ be the singular value decomposition of $B$. If $D \in {\mathbb C}^{n \times n}$ is a diagonal matrix with $D(i,i) = 0$ for $1 \leq i \leq n-1$ and $D(n,n) = - \sigma_n(B)$, then $||UDV^H||_2 = \sigma_n(B) \leq \epsilon$ and therefore eigenvalues of $(A, B + UDV^H)$ belong to $\Lambda_{\epsilon}(A,B)$. But by construction $B + UDV^H$ is singular, which means $(A, B + UDV^H)$ is singular and/or  has an eigenvalue at infinity. In both cases $\Lambda_{\epsilon}(A,B)$ must be unbounded.
\end{proof}

We next derive an upper bound on $\Lambda_{\epsilon}(A,B)$ when it exists. Keeping in mind that we plan to perturb our matrices, and will therefore have access to a lower bound on $\sigma_n(B)$ with high probability, we prove the following lemma.

\begin{lem}\label{lem: pseudo upper bound}
If $\epsilon < \sigma_n(B)$ then any $ z \in \Lambda_{\epsilon}(A,B)$ satisfies
$$ |z| \leq \frac{\epsilon ||B^{-1}||_2 + ||B^{-1}A||_2}{1 - \epsilon ||B^{-1}||_2}. $$
\end{lem}
\begin{proof}
Suppose $|z| > ||B^{-1}A||_2$ and consider $\frac{B^{-1}A}{z} - I$. For any vector $x$,
\begin{equation}
    \left| \left| \frac{B^{-1}A}{z} x \right| \right|_2 \leq \frac{||B^{-1}A||_2}{|z|} ||x||_2
    \label{eqn: vector matrix norm compatible}
\end{equation}
by matrix/vector norm compatibility, so
\begin{equation}
    ||x||_2 - \left| \left| \frac{B^{-1}A}{z} x \right| \right|_2 \geq ||x||_2 - \frac{||B^{-1}A||_2}{|z|} ||x||_2 \geq 0,
    \label{eqn: apply prev inequality}
\end{equation}
where we know $||x||_2 - (||B^{-1}A||_2/|z|)||x||_2$ is positive since $|z| > ||B^{-1}A||_2$. Thus, applying this result and the reverse triangle inequality,
\begin{equation}
    \sigma_n \left( \frac{B^{-1}A}{z} - I \right) = \min_{||x||_2 = 1} \left| \left| \left( \frac{B^{-1}A}{z} - I \right)x \right| \right|_2 \geq \min_{||x||_2 = 1} \left[ ||x||_2 - \frac{||B^{-1}A||_2}{|z|} ||x||_2 \right] = 1 - \frac{||B^{-1}A||_2}{|z|} 
    \label{eqn: reverse triangle}
\end{equation}
and therefore
\begin{equation}
    \left| \left| \left( \frac{B^{-1}A}{z} - I \right)^{-1} \right| \right|_2 = \frac{1}{\sigma_n \left( \frac{B^{-1}A}{z} - I \right)} \leq \frac{1}{1 - \frac{||B^{-1}A||_2}{|z|}} .
    \label{eqn: norm inequality}
\end{equation}
If $z \in \Lambda_{\epsilon}(A,B)$, we then have
\begin{equation} 
    \frac{1}{\epsilon (1 + |z|)} \leq \frac{||B^{-1}||_2}{|z|} \left| \left| \left( \frac{B^{-1}A}{z} - I \right)^{-1} \right| \right|_2  \leq \frac{||B^{-1}||_2}{|z|} \frac{1}{1 - \frac{||B^{-1}A||_2}{|z|}} = \frac{||B^{-1}||_2}{|z| - ||B^{-1}A||}, 
\end{equation}
which, rearranging to solve for $|z|$, is equivalent to 
\begin{equation}
    |z| \leq \frac{\epsilon ||B^{-1}||_2 + ||B^{-1}A||_2}{1 - \epsilon ||B^{-1}||_2}.
    \label{eqn: final upper bound}
\end{equation}
Since we assumed $|z| > ||B^{-1}A||_2$, we have proved that $z \in \Lambda_{\epsilon}(A,B)$ for $\epsilon < \sigma_n(B)$ implies
\begin{equation}
    |z| \leq \max \left\{ ||B^{-1}A||_2, \; \frac{\epsilon ||B^{-1}||_2 + ||B^{-1}A||_2}{1 - \epsilon ||B^{-1}||_2} \right\} = \frac{\epsilon ||B^{-1}||_2 + ||B^{-1}A||_2}{1 - \epsilon ||B^{-1}||_2}
    \label{eqn: take maximum of bounds}
\end{equation}
which finishes the proof.
\end{proof}

While this gives a nice upper bound for the pseudospectrum as a whole, the upper bound in Bauer-Fike is in terms of balls around the eigenvalues. If we add an assumption that $B^{-1}A$ is diagonalizable, we can convert the bound just obtained to a similar result, obtaining our version of Bauer-Fike for matrix pencils.

\begin{thm}[Bauer-Fike for Matrix Pencils]\label{thm: pencil BF}
Suppose $B^{-1}A$ is diagonalizable with eigenvalues $\lambda_1, \ldots, \lambda_n$ and invertible right eigenvector matrix $V$. For $\epsilon < \sigma_n(B)$ let 
$$ r_{\epsilon} = \epsilon \kappa(V) ||B^{-1}||_2 \left( 1 + \frac{\epsilon ||B^{-1}||_2 + ||B^{-1}A||_2}{1 - \epsilon ||B^{-1}||_2} \right) $$
and further set 
$$r_i = \begin{cases} 
      \frac{1}{||B||_2} & \text{if} \; \; A = 0 \\
      \max \left\{ \frac{1}{||B||_2}, \frac{|\lambda_i|}{||A||_2} \right\} & \text{otherwise}
   \end{cases} $$
for $1 \leq i \leq n$. Then,
$$ \bigcup_{i=1}^n B_{\epsilon r_i}(\lambda_i) \subseteq \Lambda_{\epsilon}(A,B) \subseteq \bigcup_{i=1}^n B_{r_{\epsilon}}(\lambda_i). $$ 
\end{thm}
   \begin{proof} To obtain the first containment we note that for any $|\Delta \lambda| \leq \epsilon/||B||_2$ the pencil $(A + \Delta \lambda B,B)$ has eigenvalues $\lambda_i + \Delta \lambda$ while $\Lambda(A+\Delta \lambda B, B) \subseteq \Lambda_{\epsilon}(A,B)$. Similarly, if $A \neq 0$ then $(A + \Delta \lambda A, B)$ has eigenvalues $\lambda_i(1 + \Delta \lambda)$ and $\Lambda(A + \Delta \lambda A,B) \subseteq \Lambda_{\epsilon}(A,B)$ as long as $|\Delta \lambda| \leq \epsilon / ||A||_2$. For the remaining containment we appeal to \cref{thm: equiv pseudo}:\ any $z \in \Lambda_{\epsilon}(A,B)$ satisfies
\begin{equation}
    \frac{1}{\epsilon(1 + |z|)} \leq ||(A - zB)^{-1}||_2 \leq ||B^{-1}||_2 ||(B^{-1}A - zI)^{-1} ||_2 .
    \label{eqn: equiv pseudo plus norm submult}
\end{equation}
Applying the fact that $V$ diagonalizes $B^{-1}A$, meaning $B^{-1}A = V\Lambda V^{-1}$ for a diagonal matrix $\Lambda$, this expression becomes
\begin{equation}
    \frac{1}{\epsilon(1 + |z|)} \leq ||B^{-1}||_2 ||(V \Lambda V^{-1} - zI)^{-1} ||_2 \leq \kappa(V) ||B^{-1}||_2 ||(\Lambda - zI)^{-1}||_2 .
    \label{eqn: apply diagonalizability}
\end{equation}
Inverting and rearranging, we then have
\begin{equation}
    \sigma_n(\Lambda - zI) \leq \epsilon \kappa(V) ||B^{-1}||_2 (1 + |z|) .
    \label{eqn: invert prev inequality}
\end{equation}
We complete the proof by noting
\begin{equation}
    \sigma_n(\Lambda - zI) = \min_{\lambda_i \in \Lambda(A,B)} |\lambda_i - z| 
    \label{eqn: eigenvalue distance}
\end{equation}
and replacing $|z|$ in \eqref{eqn: invert prev inequality} with the upper bound provided by \cref{lem: pseudo upper bound}.
\end{proof}

\indent \cref{thm: pencil BF} is not the only version of Bauer-Fike derived for the generalized eigenvalue problem (see \cref{appendix: BF}). Nevertheless, it is useful here since it leverages the assumptions we'll have available to us in the following sections. In particular, the perturbations applied to $(A,B)$ to obtain $(\widetilde{A}, \widetilde{B})$ guarantee with high probability that $\widetilde{B}$ is nonsingular (with a lower bound on $\sigma_n(\widetilde{B})$ known) and that $\widetilde{B}^{-1}\widetilde{A}$ is diagonalizable.  

\section{Pseudospectral Shattering for Matrix Pencils}
\label{section: shattering}
Following the convention used by Banks et al. \cite{banks2020pseudospectral}, define a grid on the complex plane as follows.
\begin{defn} 
    The grid $g = \text{grid}(z_0, \omega, s_1, s_2)$ is the boundary of a $s_1 \times s_2$ lattice in the complex plane consisting of $\omega \times \omega$-sized squares with lower left corner $z_0 \in {\mathbb C}$. The grid lines of $g$ are parallel to either the real or the complex axis.
    \label{defn: grid}
\end{defn}
In this section, we show that with high probability the $\epsilon$-pseudospectrum of a scaled and perturbed matrix pencil is shattered with respect to a random grid $g$ (constructed according to \cref{defn: grid}) for sufficiently small $\epsilon$, where shattering is defined as follows.
\begin{defn}
    $\Lambda_{\epsilon}(A,B)$ is shattered with respect to a grid $g$ if $ \Lambda_{\epsilon}(A,B) \cap g = \emptyset$ and each eigenvalue of $(A,B)$ is contained in a unique box of $g$.
    \label{defn: shattering}
\end{defn}
Throughout, we consider perturbations made by Ginibre random matrices and make use of results from random matrix theory for this ensemble.
\begin{defn}
    A random matrix $G \in {\mathbb C}^{n \times n}$ is Ginibre if its entries are i.i.d.\ with distribution ${\mathcal N}_{\mathbb C}(0, \frac{1}{n})$.
    \label{defn: Ginibre}
\end{defn}
\indent Note that ${\mathcal N}_{\mathbb C}(0, \frac{1}{n})$ in \cref{defn: Ginibre} is a complex normal distribution. In practice, we can construct a Ginibre matrix by drawing a complex Gaussian matrix and normalizing by $\frac{1}{\sqrt{n}}$. \\
\indent For the remainder of this section, let $(\widetilde{A}, \widetilde{B}) = (A + \gamma G_1, B + \gamma G_2)$ for $G_1$, $G_2$ two independent, complex Ginibre matrices, $0 < \gamma < \frac{1}{2}$, and $A,B \in {\mathbb C}^{n \times n}$ with $||A||_2, ||B||_2 \leq 1$. The norm requirement here is applied to simplify the analysis and can be assumed without loss of generality since our input matrices can be normalized if necessary. We are interested in showing shattering for the pseudospectra of the scaled pencil $(\widetilde{A}, n^{\alpha} \widetilde{B})$, where $\alpha$ is a positive constant. As we will see, the choice to scale $\widetilde{B}$ is required to obtain usable bounds on the locations of the eigenvalues of $(\widetilde{A}, n^{\alpha} \widetilde{B})$. Later on, we will comment on how this choice of scaling affects the results. 

\subsection{Preliminaries}
We begin by discussing some singular value bounds that will be useful later. The first of these is a bound on the norm (i.e., largest singular value) of a complex Ginibre matrix due to Banks et al. \cite[Lemma 2.2]{banks2020gaussian}.
\begin{lem}[Banks et al. 2021] For $G$ an $n \times n$ complex Ginibre matrix, $ {\mathbb P} \left[ ||G||_2 \geq 2\sqrt{2} + t \right] \leq 2e^{-nt^2} .$
\label{lem: ginibre norm bound}
\end{lem}
The same paper also provides a lower tail bound on the smallest singular value of a matrix after it is perturbed by a Ginibre matrix \cite[Lemma 3.3]{banks2020gaussian}.

\begin{lem}[Banks et al. 2021] Let $G$ be an $n \times n$ complex Ginibre matrix. Then for any $M \in {\mathbb C}^{n \times n}$ and any $\gamma, t > 0$, $ {\mathbb P} \left[ \sigma_n(M + \gamma G) < t \right] \leq n^2 \frac{t^2}{\gamma^2} .$
\label{lem: small ball estimate} 
\end{lem}

\indent Finally, we use the following result from their subsequent work \cite[Corollary 3.3]{banks2020pseudospectral}. 
\begin{lem}[Banks et al.\ 2022]\label{lem: banks_3.3}
Let $G$ be an $n \times n$ complex Ginibre matrix. Then for any $M \in {\mathbb C}^{n \times n}$ and any  $\gamma, t > 0$, ${\mathbb P} \left[ \sigma_{n-1}(M + \gamma G) < t \right] \leq 4 (\frac{tn}{\gamma})^8$.
\end{lem}
\indent These are the main results from random matrix theory that we will use. Note that \cref{lem: small ball estimate} implies that $\widetilde{B}$ is almost surely invertible, meaning with probability one we may assume that $\widetilde{B}^{-1}$ exists. Moreover, \cref{lem: ginibre norm bound} and \cref{lem: small ball estimate} allow us to control the location of the eigenvalues of $(\widetilde{A}, n^{\alpha} \widetilde{B})$ as follows.

\begin{lem} The eigenvalues of $(\widetilde{A}, n^{\alpha} \widetilde{B})$ are contained in $B_3(0)$ -- the ball of radius three centered at the origin -- with probability at least $ 1 - \frac{n^{2 - 2\alpha}}{\gamma^2} - 2e^{-n} .$
\label{lem: spectrum in disk}
\end{lem}
\begin{proof}
Consider $X = n^{-\alpha} \widetilde{B}^{-1} \widetilde{A}$. Since $X$ and $(\widetilde{A}, n^{\alpha} \widetilde{B})$ have the same eigenvalues, it is sufficient to bound the probability that $\Lambda(X)$ is not contained in $B_3(0)$. To do this, we note
\begin{equation}
    \rho(X) \leq ||X||_2 = n^{-\alpha}|| \widetilde{B}^{-1} \widetilde{A}||_2 \leq \frac{1}{n^{\alpha} \sigma_n(\widetilde{B})} ||\widetilde{A}||_2 \leq \frac{1}{n^{\alpha} \sigma_n(\widetilde{B})} \left( ||A||_2 + \gamma ||G_1||_2 \right) .
    \label{eqn: norm bound on X}
\end{equation}
Now $||A||_2 \leq 1$ so, conditioning on the event $||G_1||_2 \leq 4$, we have $||A + \gamma G_1||_2 \leq 3$. Thus, if $\Lambda(X) \nsubseteq B_3(0)$ and therefore $\rho(X) > 3$, we obtain $n^{\alpha} \sigma_n(\widetilde{B}) < 1$. Consequently,
\begin{equation}
    {\mathbb P}\left[ \Lambda(X) \nsubseteq B_3(0) \; | \; ||G_1||_2 \leq 4 \right] \leq {\mathbb P} \left[ n^{\alpha} \sigma_n(\widetilde{B}) < 1 \right] \leq \frac{n^{2-2\alpha}}{\gamma^2}, 
    \label{eqn: probability bound for X}
\end{equation}
where the last inequality comes from \cref{lem: small ball estimate}.
By Bayes' Theorem, we therefore have
\begin{equation}
    {\mathbb P} \left[ \Lambda(X) \nsubseteq B_3(0), \; ||G_1||_2 \leq 4 \right] \leq \frac{n^{2 - 2\alpha}}{\gamma^2} .
    \label{eqn: disk bound part 1}
\end{equation}
At the same time, applying \cref{lem: ginibre norm bound}, we have
\begin{equation}
    {\mathbb P} \left[ \Lambda(X) \nsubseteq B_3(0), \; ||G_1||_2 > 4 \right] \leq {\mathbb P} \left[ ||G_1||_2 > 4 \right] \leq 2e^{-n(4 - 2 \sqrt{2})^2} \leq 2e^{-n}.
    \label{eqn: disk bound part 2}
\end{equation}
Putting these two results together, we conclude
\begin{equation}
    {\mathbb P} \left[ \Lambda(X) \nsubseteq B_3(0) \right] \leq \frac{n^{2 - 2\alpha}}{\gamma^2} + 2e^{-n}
    \label{eqn: disk bound total}
\end{equation}
which completes the proof.
\end{proof}
With this in place, we can now consider building our random grid $g$ over $B_3(0)$ which (with an appropriate choice of $\alpha$ and $\gamma$) will contain every eigenvalue of $(\widetilde{A}, n^{\alpha} \widetilde{B})$ with high probability. Given a grid box size $\omega$, we define the random grid $g =  \text{grid}(z, \omega, \lceil 8 / \omega \rceil, \lceil 8 / \omega \rceil)$ for $z$ a point chosen uniformly at random from the square with bottom left corner $-4-4i$ and side length $\omega$. The construction here is somewhat arbitrary; for convenience we choose $g$ so that it roughly covers the smallest box with integer side length that contains $B_3(0)$ (with some buffer space).

\subsection{Multi-Parameter Tail Bound}
As in  \cref{lem: spectrum in disk}, let $X = n^{-\alpha} \widetilde{B}^{-1} \widetilde{A}$. The bulk of the proof of shattering in the following section rests on a multi-parameter tail bound for $\kappa_V(X)$ and $\text{gap}(X)$, which we derive here. We start with two intermediate results. The first concerns the second to last singular value of $yI - X$ for $y \in {\mathbb C}$.
\begin{lem} \label{lem: tail bound result one}
For any $t > 0$ and any $y \in  {\mathbb C}$, 
$${\mathbb P} \left[ \sigma_{n-1}(yI - X) < t \; | \; ||G_2||_2 \leq 4 \right] \leq 4 \left( \frac{t(1+4 \gamma)n^{\alpha+1}}{\gamma} \right)^8 $$
\end{lem}
\begin{proof} We begin by rewriting $yI - X$ as 
\begin{equation}
    yI - X = n^{-\alpha} \widetilde{B}^{-1} \left( n^{\alpha} y \widetilde{B} - A - \gamma G_1 \right) .
    \label{eqn: rewrite yI-X}
\end{equation}
Using standard singular value inequalities, we then have
\begin{equation}
    \sigma_{n-1}(yI - X) \geq \sigma_n \left(n^{-\alpha} \widetilde{B}^{-1} \right) \sigma_{n-1}\left(n^{\alpha} y \widetilde{B} - A - \gamma G_1 \right).
    \label{eqn: n-1 singular value inequal}
\end{equation}
Now conditioning on $||G_2||_2 \leq 4$, we have $|| \widetilde{B}||_2 \leq 1 + 4 \gamma $. Thus, $\sigma_n(n^{-\alpha} \widetilde{B}^{-1}) \geq \frac{1}{n^{\alpha}(1 + 4 \gamma)} $ and therefore
\begin{equation}
    \sigma_{n-1}(yI - X) \geq \frac{1}{n^{\alpha}(1 + 4 \gamma)} \sigma_{n-1}\left(n^{\alpha} y \widetilde{B} - A - \gamma G_1 \right).
    \label{eqn: n-1 singular value inequal 2}
\end{equation}
Consequently, we observe
\begin{equation}
    {\mathbb P} \left[ \sigma_{n-1}(yI - X) < t  \; | \; ||G_2||_2 \leq 4 \right] \leq {\mathbb P} \left[ \sigma_{n-1} \left( n^{\alpha} y \widetilde{B} - A - \gamma G_1 \right) < t n^{\alpha} (1 + 4 \gamma) \right].
    \label{eqn: n-1 singular value probability bound}
\end{equation}
Setting $M = n^{\alpha} y \widetilde{B} - A$ and applying  \cref{lem: banks_3.3} (noting that $M$ is independent of $G_1$ and that $G_1$ and $-G_1$ have the same distribution) we have
\begin{equation}
    {\mathbb P} \left[ \sigma_{n-1}(yI - X) < t \; | \; ||G_2||_2 \leq 4 \right] \leq 4 \left( \frac{t (1 + 4 \gamma) n^{\alpha + 1}}{\gamma} \right)^8 
    \label{eqn: final n-1 prob bound}
\end{equation}
for any $t > 0$.
\end{proof}
Next, we show that given any measurable set $S \subset {\mathbb C}$, we can bound the expected value of $\sum_{\lambda_i \in S} \kappa(\lambda_i)^2$, where the sum is taken over eigenvalues of $X$. This gives us some control over $\kappa_V(X)$; in particular, taking $V$ to be the eigenvector matrix of $X$ scaled so that each column is a unit vector, we observe 
\begin{equation}
    \kappa_V(X) \leq \kappa(V) \leq ||V||_F ||V^{-1}||_F \leq \sqrt{n \sum_{i=1}^n \kappa(\lambda_i)^2 }.
    \label{eqn: eigenvector condition numbers}
\end{equation}
Note that with probability one $X$ has a full set of distinct eigenvalues and is therefore diagonalizable. The proof of this result follows closely a proof of Banks et al.\ for the single matrix case \cite[Theorem 1.5]{banks2020gaussian}.
\begin{lem}
    Let $\lambda_1, \ldots \lambda_n$ be the eigenvalues of $X$. For any measurable set $S \subset {\mathbb C}$, 
    \begin{equation*}
   {\mathbb E} \left[ \sum_{\lambda_i \in S} \kappa(\lambda_i)^2 \; \vline \; ||G_2||_2 \leq 4 \right] \leq \left( \frac{(1 + 4 \gamma) n^{\alpha+1}}{\gamma} \right)^2 \frac{\text{\emph{ vol}}(S)}{\pi}.
    \end{equation*}
    \label{lem: tail bound result two}
\end{lem}
\begin{proof}
By \cref{defn: matrix pseudospectrum}, we know for any $z \in {\mathbb C}$
\begin{equation}
    {\mathbb P} \left[ z \in \Lambda_{\epsilon}(X)  \right] = {\mathbb P} \left[ \sigma_n(zI - X) < \epsilon \right] .
    \label{eqn: equiv pseudo prob}
\end{equation}
Following the same argument made in \cref{lem: tail bound result one}, swapping \cref{lem: banks_3.3} for \cref{lem: small ball estimate}, we therefore have
\begin{equation}
    {\mathbb P} \left[ z \in \Lambda_{\epsilon}(X) \; | \; ||G_2||_2 \leq 4 \right] \leq \left( \frac{ \epsilon (1 + 4 \gamma) n^{\alpha + 1}}{\gamma} \right)^2 .
    \label{eqn: smallest singular value control}
\end{equation}
Consider now the measurable set $S \subset {\mathbb C}$. Using $\eqref{eqn: smallest singular value control}$, we observe
\begin{equation}
    \aligned
    {\mathbb E} \left[ \text{vol} \left( \Lambda_{\epsilon}(X) \cap S \right) \; | \; ||G_2||_2 \leq 4 \right] & = {\mathbb E} \left[ \int_{S} \mathds{1}_{z \in \Lambda_{\epsilon}(X) \; | \; ||G_2||_2 \leq 4} dz \right] \\
    & = \int_S {\mathbb P} \left[ z \in \Lambda_{\epsilon}(X) \; | \; ||G_2||_2 \leq 4 \right] dz \\
    & \leq  \left( \frac{\epsilon (1 + 4 \gamma) n^{\alpha + 1}}{\gamma} \right)^2 \text{vol}(S). 
    \endaligned
    \label{eqn: psuedo intersection}
\end{equation}
Applying this result alongside \cite[Lemma 3.2]{banks2020gaussian} and Fatou's Lemma for conditional expectation, we conclude
\begin{equation}
    \aligned 
    {\mathbb E} \left[ \sum_{\lambda_i \in S} \kappa(\lambda_i)^2 \;  \vline \; ||G_2||_2 \leq 4 \right] &=  {\mathbb E} \left[ \liminf_{\epsilon \rightarrow 0} \frac{ \text{vol}(\Lambda_{\epsilon}(X) \cap S)}{\pi \epsilon^2} \; \vline \; ||G_2||_2 \leq 4 \right] \\
    & \leq \liminf_{\epsilon \rightarrow 0} \frac{ {\mathbb E} \left[ \text{vol}(\Lambda_{\epsilon}(X) \cap S \; | \;  ||G_2||_2 \leq 4 \right]}{ \pi \epsilon^2} \\
    & \leq \left( \frac{(1 + 4 \gamma) n^{\alpha + 1}}{\gamma} \right)^2 \frac{\text{vol}(S)}{\pi}, 
    \endaligned
    \label{eqn: expectation liminf}
\end{equation}
which completes the proof.
\end{proof}

With these in place, we are now ready to prove the main multi-parameter tail bound. In this result, as in the previous two, $X = n^{-\alpha} \widetilde{B}^{-1} \widetilde{A}$ for our perturbed pencil $(\widetilde{A}, \widetilde{B})$. Again, this follows closely the proof of an equivalent tail bound of Banks et al.\ \cite[Theorem 3.6]{banks2020pseudospectral}.
\begin{thm}
    For any $t, \delta > 0$, 
    \begin{equation*}
        {\mathbb P} \left[ \kappa_V(X) < t, \; \text{\emph{gap}}(X) > \delta \right] \geq \left[ 1 - \frac{9(1+4\gamma)^2 n^{2 \alpha+3}}{t^2 \gamma^2} - 1296 \delta^6 \left( \frac{t(1+4\gamma)n^{\alpha+1}}{\gamma} \right)^8 \right] \left[ 1 - \frac{n^{2-2\alpha}}{\gamma^2} - 4e^{-n} \right]. 
    \end{equation*}
\label{thm: tail bound}
\end{thm}
\begin{proof}
We condition on the events $||G_2||_2 \leq 4$ and $\Lambda_{\epsilon}(X) \subseteq B_3(0)$. Combining \cref{lem: ginibre norm bound} and  \cref{lem: spectrum in disk}, we know these occur with probability at least $1 - \frac{n^{2 - 2\alpha}}{\gamma^2} - 4e^{-n}$. Define now the events:
\begin{itemize}
    \item $E_{\text{gap}} = \left\{ \text{gap}(X) < \delta \right\}$
    \item $E_{\kappa} = \left\{ \kappa_V(X) > t \right\}$.
\end{itemize}
We are interested in bounding the probability of $E_{\text{gap}} \cup E_{\kappa}$. To do this, we construct a minimal $\frac{\delta}{2}$-net ${\mathcal N}$ covering $B_3(0)$, which exists with $|{\mathcal N}| \leq \frac{324}{\delta^2}$ (see for example \cite[Corollary 4.2.11]{vershynin_HDP}). By construction
\begin{equation}
    E_{\text{gap}} \subset \left\{ |D(y, \delta) \cap \Lambda(X) | \geq 2 \; \text{for some} \; y \in {\mathcal N} \right\}
    \label{eqn: Egap containment}
\end{equation}
where $D(y, \delta)$ is the ball of radius $\delta$ centered at $y$. Now if $D(y, \delta)$ contains two eigenvalues of $X$, \cite[Lemma 3.5]{banks2020pseudospectral} implies $\sigma_{n-1}(yI - X) < \delta \kappa_V(X)$, where we note again that $X$ is almost surely diagonalizable. Thus, if $\text{gap}(X) < \delta$ either $\sigma_{n-1}(yI -X) < \delta t$ for at least one $y \in {\mathcal N}$ or $\kappa_V(X) > t$. In other words, defining the event $E_y = \left\{ \sigma_{n-1}(yI - X) < \delta t \right\}$ for $y \in {\mathcal N}$,
\begin{equation}
    E_{\text{gap}} \subset E_{\kappa} \cup \bigcup_{y \in {\mathcal N}} E_y \; \Longrightarrow \; E_{\text{gap}} \cup E_{\kappa} \subset E_{\kappa} \cup \bigcup_{y \in {\mathcal N}} E_y.
    \label{eqn: event containment}
\end{equation}
By a union bound, we then have
\begin{equation}
    {\mathbb P}[E_{\text{gap}} \cup E_{\kappa}] \leq {\mathbb P}[E_{\kappa}] + |{\mathcal N}|\max_{y \in {\mathcal N}} {\mathbb P}[E_y]. 
    \label{eqn: union bound}
\end{equation}
To bound ${\mathbb P}[E_{\kappa}]$, we first use \eqref{eqn: eigenvector condition numbers} to obtain
\begin{equation}
    {\mathbb P}[E_{\kappa}] \leq {\mathbb P} \left[ t < \sqrt{n \sum_{i=1}^n \kappa(\lambda_i)^2} \; \right] = {\mathbb P} \left[ \sum_{i=1}^n \kappa(\lambda_i)^2 > \frac{t^2}{n} \right] .
    \label{E_k bound 1}
\end{equation}
Applying Markov's inequality and \cref{lem: tail bound result two} (recalling that we've conditioned on $\Lambda(X) \subseteq B_3(0)$) we then have
\begin{equation}
    {\mathbb P}[E_{\kappa}] \leq \frac{n}{t^2} {\mathbb E} \left[ \sum_{\lambda_i \in B_3(0)} \kappa(\lambda_i)^2 \right] \leq \frac{9 (1 + 4 \gamma)^2 n^{2 \alpha + 3}}{t^2 \gamma^2}.
    \label{eqn: E_k bound 2}
\end{equation}
Similarly, \cref{lem: tail bound result one} implies
\begin{equation}
    {\mathbb P} [E_y]  \leq 4 \left( \frac{\delta t (1 + 4 \gamma) n^{\alpha+1}}{\gamma} \right)^8 
    \label{eqn: E_y bound}
\end{equation}
for all $y \in {\mathcal N}$. Putting everything together, we conclude
\begin{equation}
    {\mathbb P} [ E_{\text{gap}} \cup E_{\kappa}] \leq \frac{9 (1 + 4\gamma)^2 n^{2 \alpha + 3}}{t^2 \gamma^2} + 4 \left(\frac{\delta t (1 + 4 \gamma)n^{\alpha+1}}{\gamma} \right)^8 \frac{324}{\delta^2} 
    \label{eqn: take union bound}
\end{equation}
so we have shown
\begin{equation}
    \resizebox{.9\linewidth}{!}{
    ${\mathbb P} \left[ \kappa_V(X) < t, \; \text{gap}(X) > \delta \; | \; \Lambda(X) \subseteq B_3(0), \; ||G_2||_2 \leq 4 \right] \geq 1 - \frac{9(1+4 \gamma)^2 n^{2\alpha+3}}{t^2 \gamma^2} - 1296 \delta^6 \left( \frac{t (1 + 4\gamma) n^{\alpha+1}}{\gamma} \right)^8 $}.
    \label{eqn: conditional tail bound}
\end{equation}
Noting 
\begin{equation}
    {\mathbb P}[ \kappa_V(X) < t, \; \text{gap}(X) > \delta] \geq  {\mathbb P} \left[ \kappa_V(X) < t, \; \text{gap}(X) > \delta, \; \Lambda(X) \subseteq B_3(0), \; ||G_2||_2 \leq 4 \right],
    \label{eqn: tail bound plus Bayes}
\end{equation}
we obtain the final result by Bayes' Theorem.
\end{proof}

\begin{remark}
Note that we could replace $\kappa_V(X)$ in \cref{thm: tail bound} with $\kappa(V)$ for any eigenvector matrix $V$ satisfying \eqref{eqn: eigenvector condition numbers}. In particular, as was used to derive the inequality, this applies to the scaling where each column of $V$ has unit length.
\label{rem: condition number}
\end{remark}

Informally, \cref{thm: tail bound} says that for the right choice of $t$ and $\delta$ (i.e., polynomial in $n$ and inverse polynomial in $n$, respectively) the eigenvector condition number of $X$ is not too large and its eigenvalue gap is not too small. Since the eigenvalues and eigenvectors of $X$ and $(\widetilde{A}, n^{\alpha} \widetilde{B})$ are the same, this eigenvector and eigenvalue stability is inherited by our perturbed and scaled pencil. 

\subsection{Shattering}
We are now ready to prove shattering for the pseudospectra of $(\widetilde{A}, n^{\alpha} \widetilde{B})$.

\begin{thm}\label{thm: shattering}
Let $A,B \in {\mathbb C}^{n \times n}$ with $||A||_2 \leq 1$ and $||B||_2 \leq 1$. Let $(\widetilde{A}, \widetilde{B}) = (A + \gamma G_1, B + \gamma G_2)$ for $G_1, G_2$ two independent complex Ginibre matrices and $0 < \gamma < \frac{1}{2}$. Let $g$ be the grid
$$ g = \text{\normalfont grid}(z, \omega, \lceil 8/ \omega \rceil, \lceil 8/\omega \rceil) \; \; \; \text{with} \; \; \; \omega = \frac{\gamma^4}{4}n^{-\frac{8 \alpha + 13}{3}}, $$
where $z$ is chosen uniformly at random from the square with bottom left corner $-4 - 4i$ and side length $\omega$. Then $||\widetilde{A}||_2 \leq 3$, $||\widetilde{B}||_2 \leq 3$, and $\Lambda_{\epsilon}(\widetilde{A}, n^{\alpha} \widetilde{B})$ is shattered with respect to $g$ for 
$$ \epsilon =  \frac{\gamma^5}{64 n^{\frac{11 \alpha + 25}{3}} + \gamma^5} $$ 
with probability at least $ \left[ 1 - \frac{82}{n} - \frac{531441}{16n^2} \right] \left[ 1 - \frac{n^{2 - 2\alpha}}{\gamma^2} - 4e^{-n} \right] .$
\end{thm}
\begin{proof} As in the previous section, let $X = n^{-\alpha} \widetilde{B}^{-1} \widetilde{A}$. We condition on the following events:\ $n^{\alpha}\sigma_n(\widetilde{B}) \geq 1$, $||G_1||_2 \leq 4$, and $||G_2||_2 \leq 4$. As noted above, these events occur with probability at least $ 1 - \frac{n^{2 - 2\alpha}}{\gamma^2} - 4e^{-n} $ and, following the argument made in \cref{lem: spectrum in disk}, guarantee that $\Lambda(X) \subseteq B_3(0)$. Consequently, we know with probability one each eigenvalue of $X$, and therefore of $(\widetilde{A}, n^{\alpha} \widetilde{B})$, is contained in a box of $g$. At the same time,
\begin{equation}
    || \widetilde{A} ||_2 = ||A + \gamma G_1||_2 \leq ||A||_2 + \gamma ||G_1||_2 \leq 3
    \label{eqn: perturb norm change}
\end{equation}
and similarly $||\widetilde{B}||_2 \leq 3$. \\
\indent Suppose now $\kappa_V(X) < \frac{n^{\alpha+2}}{ \gamma}$, $\text{gap}(X) > \gamma^4n^{-\frac{8 \alpha + 13}{3}}$, and 
\begin{equation}
    \min_{\lambda_i \in \Lambda(X)} \text{dist}_g(\lambda_i) > \frac{\omega}{4n^2},
    \label{eqn: min dist requirement}
\end{equation}
where $\text{dist}_g(\lambda_i) = \min_{y \in g} |\lambda_i - y|$. By \eqref{eqn: conditional tail bound}, we know the first two of these occur under our assumptions with probability at least
\begin{equation}
    1 - \frac{9(1+4 \gamma)^2 n^{2 \alpha+3}}{\left( \frac{n^{\alpha+2}}{\gamma} \right)^2 \gamma^2} - 1296 \left(\frac{\gamma^4}{n^{\frac{8\alpha+13}{3}}} \right)^6 \left( \frac{ \left( \frac{n^{\alpha + 2}}{\gamma} \right) (1 + 4 \gamma) n^{\alpha+1}}{\gamma} \right)^8,
    \label{eqn: plug in conditional tail bound}
\end{equation}
which, applying $\gamma < \frac{1}{2}$, simplifies to $1 - \frac{81}{n} - \frac{531441}{16n^2}.$ Similarly, by a simple geometric argument (using the fact that the eigenvalues of $X$ are uniformly distributed in their grid boxes),
\begin{equation}
    {\mathbb P} \left[ \min_{\lambda_i \in \Lambda(X)} \text{dist}_g(\lambda_i) > \frac{\omega}{4n^2} \right] \geq 1 - \frac{1}{n}.
    \label{eqn: eig grid distance prob}
\end{equation}
Thus, these events occur under our assumptions with probability at least $1 - \frac{82}{n} - \frac{531441}{16n^2}$, which means by Bayes' Theorem that we have $n^{\alpha} \sigma_n(\widetilde{B}) \geq 1$, $||G_1||_2 \leq 4$, $||G_2||_2 \leq 4$, $\kappa_V(X) < \frac{n^{\alpha + 2}}{\gamma}$, $\text{gap}(X) > \gamma^4n^{-\frac{8 \alpha + 13}{3}}$, and \eqref{eqn: min dist requirement} with probability at least $\left[ 1 - \frac{82}{n} - \frac{531441}{16n^2} \right] \left[ 1 - \frac{n^{2 - 2\alpha}}{\gamma^2} - 4e^{-n} \right] $. \\
\indent To complete the proof, we now show that these events guarantee shattering. To do this, we first observe that if $\text{gap}(X) > \gamma^4n^{-\frac{8 \alpha + 13}{3}}$ then no two eigenvalues can share a grid box of $g$. At the same time, \eqref{eqn: min dist requirement} implies that the ball of radius $\frac{\omega}{ 4n^2}$ around each eigenvalue does not intersect $g$. Therefore, it is sufficient to show that $\Lambda_{\epsilon}(\widetilde{A}, n^{\alpha} \widetilde{B})$ is contained in these balls, which we can do by appealing to \cref{thm: pencil BF}, our version of Bauer-Fike for matrix pencils. In particular, noting that we can replace $\kappa(V)$ with $\kappa_V(X)$ by taking an infimum, \cref{thm: pencil BF} implies that $(\widetilde{A}, n^{\alpha} \widetilde{B})$ is contained in balls of radius
\begin{equation}
    r_{\epsilon} = \epsilon \kappa_V(X) ||n^{-\alpha}\widetilde{B}^{-1}||_2 \left( 1 + \frac{\epsilon ||n^{-\alpha}\widetilde{B}^{-1}||_2 + ||X||_2}{1 - \epsilon ||n^{- \alpha} \widetilde{B}^{-1}||_2} \right).
    \label{eqn: BF for shattering}
\end{equation}
Applying our bounds $||n^{-\alpha} \widetilde{B}^{-1}||_2 = \frac{1}{n^{\alpha} \sigma_n(\widetilde{B})} \leq 1$, $||X||_2 \leq 3$, and $\kappa_V(X) \leq \frac{n^{\alpha + 2}}{ \gamma}$, we obtain shattering as long as 
\begin{equation}
    r_{\epsilon} \leq  \epsilon \left( \frac{n^{\alpha + 2}}{\gamma} \right) \left( \frac{4}{1 - \epsilon} \right) \leq \frac{ \omega}{4n^2} = \frac{\gamma^4}{16n^{\frac{8 \alpha + 19}{3}}}  
    \label{eqn: shattering critera}
\end{equation}
or, equivalently, $\epsilon \leq \frac{\gamma^5}{64 n^{\frac{11 \alpha + 25}{3}} + \gamma^5}$.
\end{proof}
This result clarifies the impact of the $n^{\alpha}$ scaling. On one hand, increasing $\alpha$ drives the term $\frac{n^{2 - 2 \alpha}}{\gamma^2}$ in our probability bound to zero, assuming $\gamma$ is fixed. Said another way, a larger choice of $\alpha$ allows us to take $\gamma$ smaller without losing our guarantee of shattering. In particular, we need
\begin{equation}
    \gamma > n^{1-\alpha}
    \label{eqn: gamma and alpha}
\end{equation}
to ensure that the probability in \cref{thm: shattering} is not vacuous. This is important, as we'd like to perturb our matrices as little as possible. Nevertheless, we pay a penalty for increasing $\alpha$ in $\omega$ and $\epsilon$, both of which shrink as $\alpha$ increases. This trade-off reflects a fundamental geometric reality:\ the more we scale by, the closer the eigenvalues of $(\widetilde{A}, n^{\alpha}\widetilde{B})$ are driven  to zero (and therefore to each other), meaning we'll need a finer grid and a smaller pseudospectrum to guarantee shattering. 

\begin{figure}[t]
    \centering
    \begin{subfigure}{.5\linewidth}
        \centering
        \includegraphics[width= .95\linewidth]{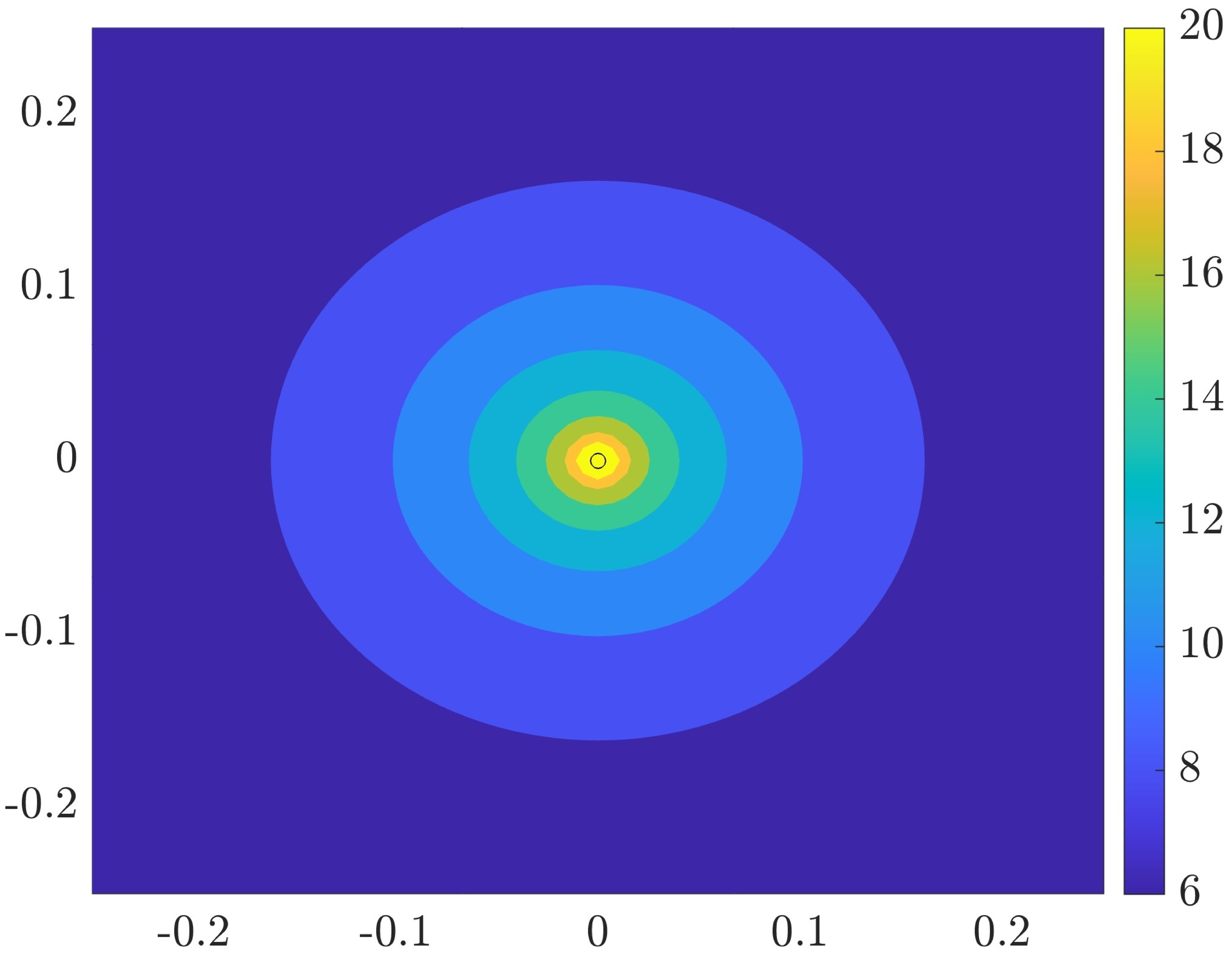}
        \caption{Pseudospectra of $(A,B)$}
    \end{subfigure}%
    \begin{subfigure}{.5\linewidth}
        \centering
        \includegraphics[width= .95\linewidth]{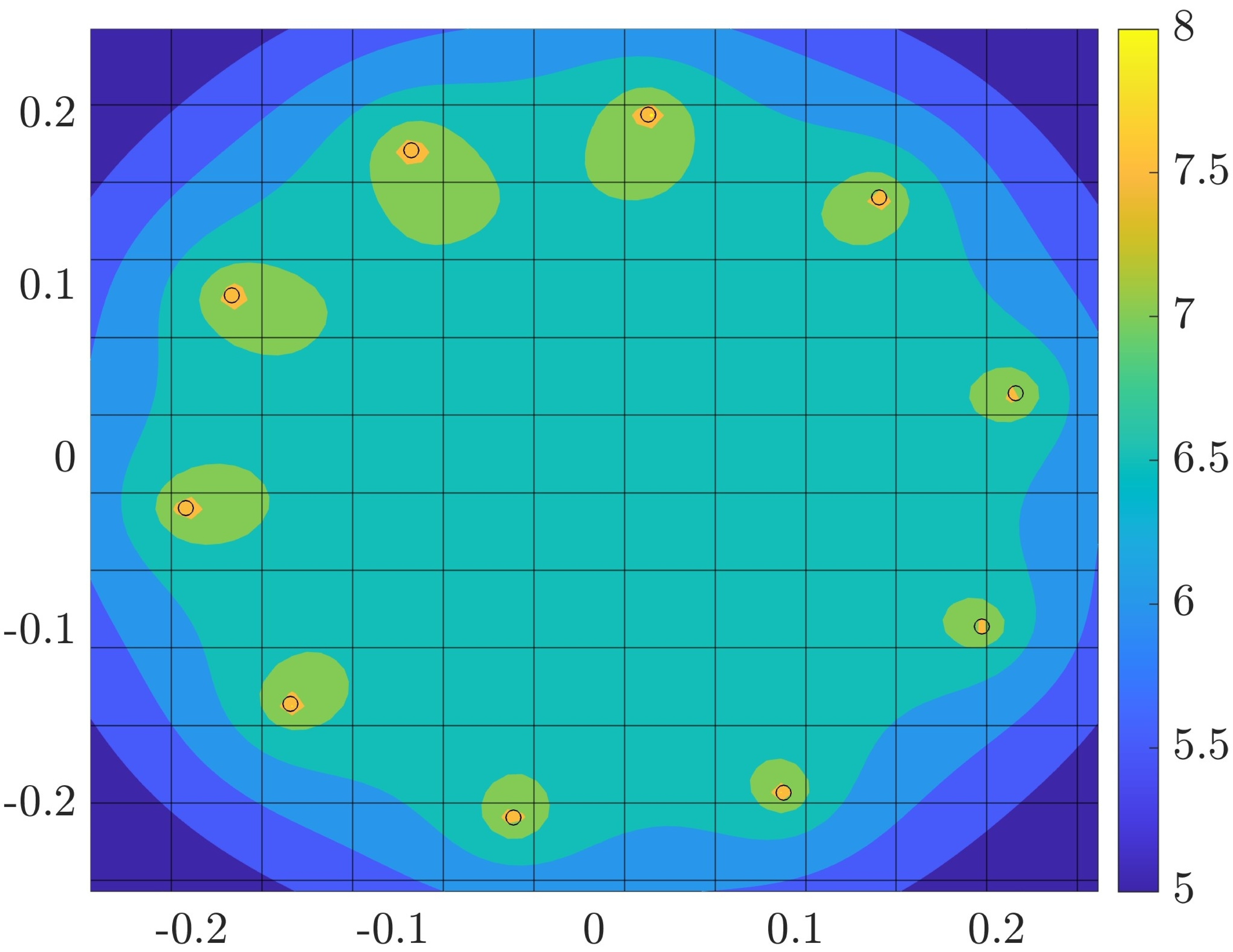}
        \caption{Pseudospectra of $(\widetilde{A}, \widetilde{B})$}
    \end{subfigure}
\caption{Pseudospectra of a pencil $(A,B)$ before and after perturbation. Here, $A$ is a $10 \times 10$ Jordan block, $B$ is the identity matrix, and the perturbed matrices are $\widetilde{A} = A + 10^{-7}G_1$ and $\widetilde{B} = B + 10^{-7}G_2$ for $G_1$ and $G_2$ two independent, complex Gaussian matrices. Once again eigenvalues are plotted with open circles. A grid that shatters the $10^{-8}$-pseudospectrum of $(\widetilde{A}, \widetilde{B})$ is superimposed over figure (b). }
    \label{fig: pseudospectral shattering}
\end{figure}

\cref{fig: pseudospectral shattering} demonstrates pseudospectral shattering on a $10 \times 10$ example. We observe here all of the same phenomena as in single-matrix shattering:\ small perturbations turn a pencil with a poorly conditioned eigenvalue problem -- here only one eigenvalue repeated ten times -- to one with a full set of distinct eigenvalues. Note in particular how much faster the pseudospectra of $(\widetilde{A}, \widetilde{B})$ coalesce around the eigenvalues as $\epsilon \rightarrow 0$. \\
\indent Because \cref{thm: tail bound} applies to the product matrix $X = n^{-\alpha} \widetilde{B}^{-1} \widetilde{A}$, our proof of shattering only works with $(\widetilde{A}, n^{\alpha} \widetilde{B})$ as a pencil via the version of Bauer-Fike used. Consequently, swapping \cref{thm: pencil BF} for \cref{thm: bauer fike}, the same proof implies the following pseudospectral shattering result for $X$.
\begin{prop} \label{prop: shattering_for_prod}
Let $A,B \in {\mathbb C}^{n \times n}$ and let $0 < \gamma < \frac{1}{2}$. Suppose $(\widetilde{A}, \widetilde{B}) = (A + \gamma G_1, B + \gamma G_2)$ for $G_1, G_2$ two independent complex Ginibre matrices and let $X = n^{- \alpha} \widetilde{B}^{-1} \widetilde{A}$ for $\alpha > 0$. Then $\Lambda_{\epsilon}(X)$ is shattered with respect to the grid $g$ (as defined in \cref{thm: shattering}) for $\epsilon =  \frac{\gamma^5}{16}n^{-\frac{11 \alpha + 25}{3}}$ with probability at least $ \left[ 1 - \frac{82}{n} - \frac{531441}{16n^2} \right] \left[ 1 - \frac{n^{2 - 2\alpha}}{\gamma^2} - 4e^{-n} \right] .$
\end{prop}
\cref{prop: shattering_for_prod} implies an alternative strategy for diagonalizing $(A,B)$:\ form the product matrix $X$ and apply single-matrix divide-and-conquer as defined by Banks et al. \cite[Algorithm EIG]{banks2020pseudospectral} using the grid $g$ and the corresponding pseudospectral guarantee. While this is fairly straightforward -- and even implies the same asymptotic complexity for the diagonalization as we derive in the subsequent sections -- it is not viable in general due to potential numerical instability. If $B$ is poorly conditioned and $\gamma$ is small, inverting $\widetilde{B}$ to form $X$ will incur significant error in finite precision.

\begin{figure}[t]
    \centering
    \begin{subfigure}{.5\linewidth}
        \centering
        \includegraphics[width=.82\linewidth]{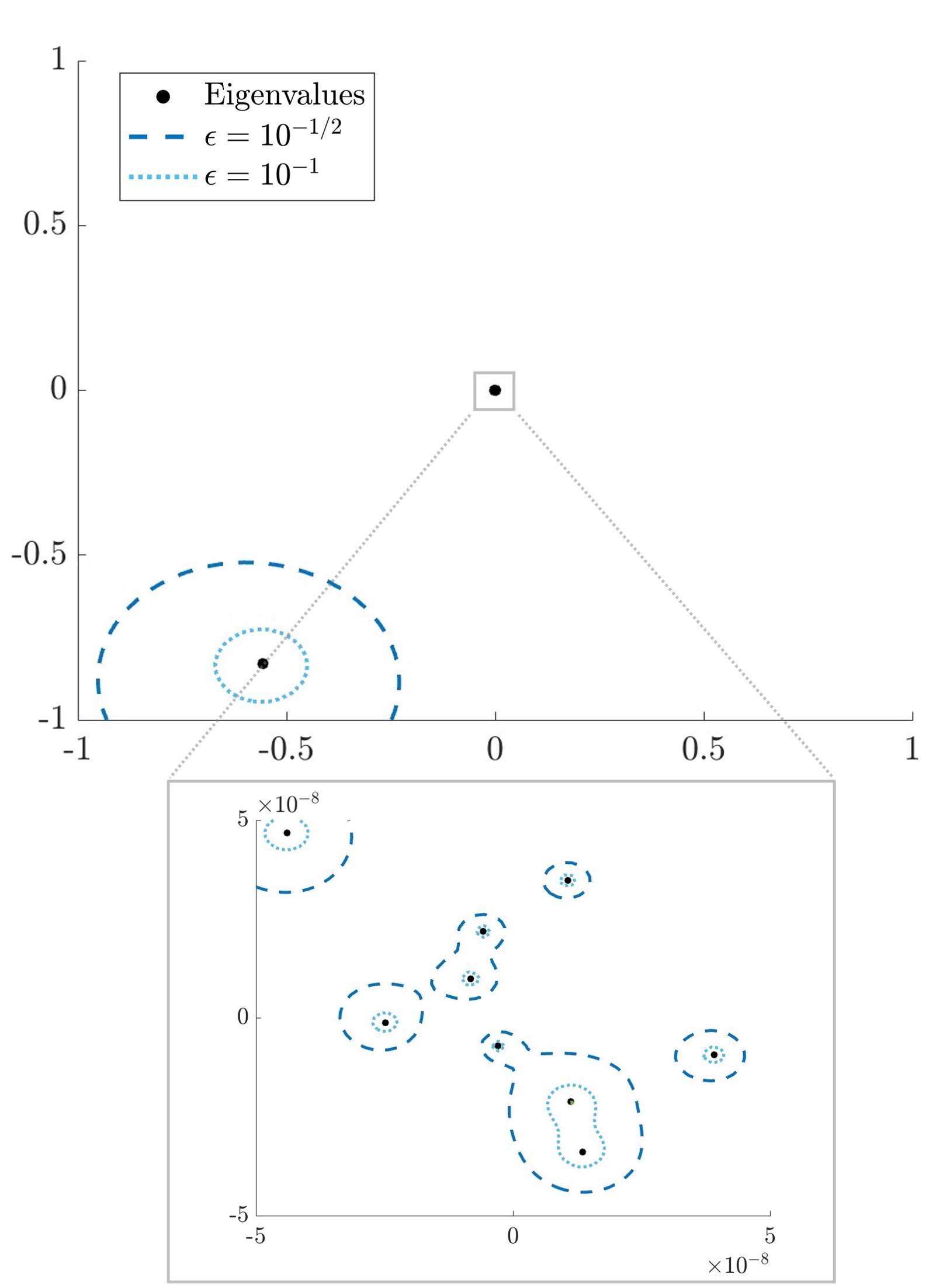}
        \caption{$\Lambda_{\epsilon}(\widetilde{A}, \widetilde{B})$}
    \end{subfigure}%
    \begin{subfigure}{.5\linewidth}
        \centering
        \includegraphics[width=.82\linewidth]{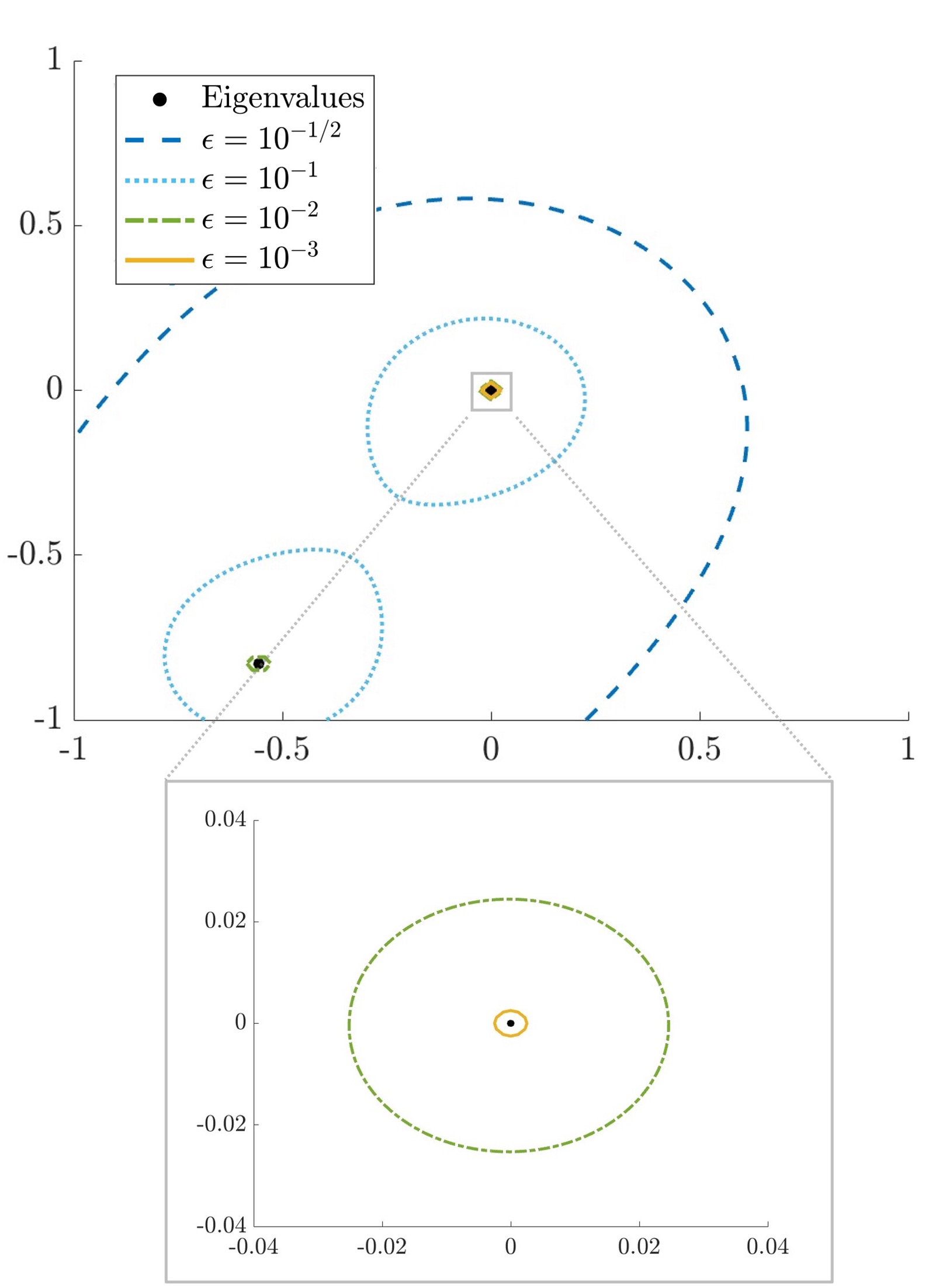}
        \caption{$\Lambda_{\epsilon}(\widetilde{B}^{-1} \widetilde{A}$)}
    \end{subfigure}
\caption{Pseudospectra of a $10 \times 10$ pencil $(\widetilde{A}, \widetilde{B})$ and its corresponding product matrix $\widetilde{B}^{-1}\widetilde{A}$. Here, $\widetilde{A} = A + 10^{-7}G_1$ and $\widetilde{B} = B + 10^{-7}G_2$ for $A$ and $B$ drawn randomly and $B$ modified to be singular (without changing its remaining singular values). In both plots, $\widetilde{B}$ is scaled so that the maximum eigenvalue has modulus one. In addition, we provide a subplot focused around the origin to examine the pseudospectra around the small eigenvalues of $(\widetilde{A}, \widetilde{B})$, which correspond to finite eigenvalues of $(A,B)$. In plot (a), we omit the pseudospectra for $\epsilon = 10^{-2}$ and $\epsilon = 10^{-3}$ since they are too close to the eigenvalues to be visible.}
    \label{fig: zoomed_pseudospectra_comp2}
\end{figure}

\indent This phenomenon is illustrated in \cref{fig: zoomed_pseudospectra_comp2}, which plots the pseudospectra of both $(\widetilde{A}, \widetilde{B})$ and $\widetilde{B}^{-1} \widetilde{A}$ for a $10 \times 10$ example in which $B$ is initially singular. The difference between the two is striking; while $\Lambda_{\epsilon}(\widetilde{A}, \widetilde{B})$ easily separates the eigenvalues, $\Lambda_{\epsilon}(\widetilde{B}^{-1} \widetilde{A})$ covers large regions of the complex plane, even for small $\epsilon$. This indicates that the guarantee provided by \cref{prop: shattering_for_prod} is more vulnerable in finite precision than \cref{thm: shattering}. Indeed, we prove in \cref{appendix: finite arithmetic} that higher precision is needed to guarantee shattering for the product matrix $X$. As mentioned in the introduction, this increase in precision is our motivation for bypassing this approach to the problem and prioritizing inverse-free computations. \\
\indent In the remainder of this section, we develop two perturbation results:\ one for pseudospectral shattering (\cref{lem: shattering preserved}) and one for individual eigenvalues and eigenvectors (\cref{lem: shattered eigenvectors}). These follow closely \cite[Lemmas 5.8 and 5.9]{banks2020pseudospectral} and will be useful in the next section.

\begin{lem}\label{lem: shattering preserved}
Suppose $(A,B)$ is regular and $\Lambda_{\epsilon}(A,B)$ is shattered with respect to a finite grid $g$. If $||A - A'||_2, ||B-B'||_2 \leq \eta < \epsilon$ then each eigenvalue of $(A',B')$ shares a grid box with exactly one eigenvalue of $(A,B)$ and $\Lambda_{\epsilon - \eta}(A',B')$ is also shattered with respect to $g$.
\end{lem}
\begin{proof}
If $z \in \Lambda_{\epsilon - \eta}(A',B')$ then $z$ is an eigenvalue of a pencil $(C,D)$ with $||A' - C||_2, ||B' - D||_2 \leq \epsilon - \eta$. In this case,
\begin{equation}
    ||A - C||_2 \leq || A - A' + A' - C||_2 \leq ||A - A'||_2 + ||A' - C||_2 \leq \eta + \epsilon - \eta = \epsilon
    \label{eqn: pseudo triangle}
\end{equation}
and similarly $||B - D||_2 \leq \epsilon$, which implies $z \in \Lambda_{\epsilon}(A,B)$. Thus, $\Lambda_{\epsilon - \eta}(A',B') \subseteq \Lambda_{\epsilon}(A,B)$, which guarantees that $\Lambda_{\epsilon - \eta}(A',B')$ is also shattered with respect to $g$. To show that each eigenvalue of $(A',B')$ shares a grid box with exactly one eigenvalue of $(A,B)$, consider $A_t = A + t(A'-A)$ and $B_t = B + t(B'-B)$ for $t \in [0,1]$. Since $(A,B)$ is regular (and moreover $\epsilon < \sigma_n(B)$ since $\Lambda_{\epsilon}(A,B)$ is bounded), $(A_t,B_t)$ continuously deforms the eigenvalues of $(A,B)$ to eigenvalues of $(A',B')$ while staying within $\Lambda_{\eta}(A,B) \subseteq \Lambda_{\epsilon}(A,B)$. Since $\Lambda_{\epsilon}(A,B)$ is shattered with respect to $g$ and therefore no two eigenvalues of $(A,B)$ share a grid box, this ensures that each eigenvalue of $(A',B')$ belongs to a grid box with a unique eigenvalue of $(A,B)$
\end{proof}

\begin{lem}\label{lem: shattered eigenvectors}
Suppose $(A,B)$ is regular and $\Lambda_{\epsilon}(A,B)$ is shattered with respect to a finite grid $g$ with boxes of side length $\omega$ . If $||A - A'||_2, ||B - B'||_2 \leq \eta < \epsilon$ then for any right unit eigenvector $v'$ of $(A',B')$ there exists a right unit eigenvector  $v$ of $(A,B)$ such that 
\begin{enumerate}
    \item The eigenvalue of $(A',B')$ corresponding to $v'$ shares a grid box of $g$ with the eigenvalue of $(A,B)$ that corresponds to $v$.
    \item $ ||v' - v|| \leq \frac{ \sqrt{8} \omega}{ \pi} \frac{ \eta }{\epsilon (\epsilon - \eta)} (1 + ||B^{-1}A||_2) ||B||_2  .$
\end{enumerate}
\end{lem}
\begin{proof}
Let $\lambda'$ be the eigenvalue of $(A',B')$ corresponding to $v'$. By \cref{lem: shattering preserved}, $\lambda'$ shares a grid box of $g$ with a unique eigenvalue of $\lambda$ of $(A,B)$. Let $v$ be the right unit eigenvector of $(A,B)$ corresponding to $\lambda$. In addition, let $w'$ and $w$ be the left eigenvectors corresponding to $v'$ and $v$ respectively, normalized so that $ w^Hv = (w')^Hv' = 1 .$ If $\Gamma$ is the contour of the grid box containing both $\lambda$ and $\lambda'$, then by \eqref{eqn: Cauchy + resolvent}
\begin{equation}
    \aligned 
    v'(w')^H - vw^H &=  \frac{1}{2 \pi i} \oint_{\Gamma} (z - (B')^{-1}A')^{-1}dz - \frac{1}{2 \pi i} \oint_{\Gamma} (z - B^{-1}A)^{-1}dz \\
    &= \frac{1}{2 \pi i} \oint_{\Gamma} (z - (B')^{-1}A')^{-1}  - (z - B^{-1}A)^{-1} dz .
    \endaligned
    \label{eqn: eigenvector contour integral}
\end{equation}
By the resolvent identity
\begin{equation}
    \aligned 
    (z - (B')^{-1}A')^{-1} - (z - B^{-1}A)^{-1} & = (z - (B')^{-1}A')^{-1}((B')^{-1}A' - B^{-1}A)(z - B^{-1}A)^{-1} \\
    & = (zB' - A')^{-1}B' ((B')^{-1}A' - B^{-1}A)(zB - A)^{-1} B \\
    & = (zB' - A)^{-1} (A' - B'B^{-1}A)(zB - A)^{-1}B,
    \label{eqn: resolvent identity}
\endaligned 
\end{equation}
so, applying this along with the ML inequality above, we have
\begin{equation}
     || v'(w')^H - vw^H||_2 \leq \frac{2 \omega}{\pi} \sup_{z \in \Gamma} || (zB' - A')^{-1}(A' - B'B^{-1}A)(zB - A)^{-1}B ||_2 .
     \label{eqn: ML plus resolvent}
\end{equation}
Now $\Lambda_{\epsilon}(A,B)$ is shattered with respect to $g$ and therefore does not intersect $\Gamma$, so we know $||(zB - A)^{-1}||_2 \leq \frac{1}{\epsilon (1 + |z|)} \leq \frac{1}{\epsilon}$ for all $z \in \Gamma$. Moreover, again by \cref{lem: shattering preserved}, $\Lambda_{\epsilon - \eta}(A',B')$ is shattered with respect to $g$, so the same argument implies $||(zB' - A')^{-1}||_2 \leq \frac{1}{\epsilon - \eta} $ for all $z \in \Gamma$. Applying this to the previous inequality, we conclude
\begin{equation}
    ||v'(w')^H - vw^H||_2 \leq \frac{2 \omega}{\pi} \frac{1}{\epsilon (\epsilon - \eta)} ||A' - B'B^{-1}A||_2 ||B||_2. 
    \label{eqn: ML plus resolvent simplified}
\end{equation}
Writing $B' = B + E$ for some $||E||_2 \leq \eta$ since $||B - B'||_2 \leq \eta$, we have
\begin{equation}
    ||A'- B'B^{-1}A||_2 = || A' - (B + E)B^{-1}A||_2 \leq ||A' - A||_2 + ||EB^{-1}A||_2 \leq \eta(1  + ||B^{-1}A||_2 )
    \label{eqn: triangle simplification}
\end{equation}
which yields a final upper bound
\begin{equation} 
    ||v'(w')^H - vw^H||_2 \leq \frac{2 \omega}{\pi} \frac{\eta}{\epsilon( \epsilon - \eta)} ( 1 + ||B^{-1}A||_2)) ||B||_2 .
    \label{eqn: eigenvector upper bound}
\end{equation}
Since without loss of generality we can assume $v^Hv' \geq 0$ (if this is not true we can just rotate $v$), and using the fact that $||v|| = ||v'|| = 1$, we complete the proof by observing 
\begin{equation} 
    ||v' - v||_2 = \sqrt{2 - 2 v^Hv'} \leq \sqrt{2} ||v'(w')^H - vw^H||_2. 
    \label{eqn: eigenvector upper bound last step}
\end{equation}
This inequality is nontrivial (see the proof of \cite[Lemma 5.8]{banks2020pseudospectral} for the details).  Combining it with \eqref{eqn: eigenvector upper bound}, we conclude that $v$ is the desired eigenvector of $(A,B)$ with $||v' - v|| \leq \frac{ \sqrt{8} \omega}{\pi} \frac{\eta}{\epsilon (\epsilon - \eta)} ( 1 + ||B^{-1}A||_2) ||B||_2 $.
\end{proof}

\section{Inverse-Free Divide-and-Conquer}
\label{section: divide and conquer}

In this section, we develop a divide-and-conquer eigensolver that, based on the pseudospectral shattering result proved in the previous section, produces approximations of the eigenvalues and eigenvectors of the perturbed and scaled pencil $(\widetilde{A}, n^{\alpha} \widetilde{B})$. We first provide motivation and background before stating the divide-and-conquer algorithm and proving convergence in exact arithmetic. In the next section, we use this divide-and-conquer routine to produce approximate diagonalizations for arbitrary pencils.

\subsection{Motivation}
\label{subsec: divide and conquer motivation}
Suppose we have a regular matrix pencil $(A,B)$ to be split by divide-and-conquer. Let $P_R$ be a rank-$k$ spectral projector onto the span of a set of right eigenvectors of $(A,B)$ and let $U_R \in {\mathbb C}^{n \times k}$ be a matrix whose orthonormal columns span $\text{range}(P_R)$. Similarly, let $U_L \in {\mathbb C}^{n \times k}$ be a matrix whose orthonormal columns span the left deflating subspace of $(A,B)$ corresponding to $\text{range}(P_R)$. Since $(\lambda, v)$ is an eigenpair of $(A,B)$ if and only if $(\lambda, w)$ is an eigenpair of $(U_L^HAU_R, U_L^HAU_R)$ for $v = U_Rw$, we can find eigenpairs of $(A,B)$ by solving the $k \times k$ problem $(U_L^HAU_R, U_L^HAU_R)$, whose spectra is contained in $\Lambda(A,B)$ and whose eigenvectors are in simple correspondence with those of $(A,B)$.  
This observation is the heart of divide-and-conquer, which proceeds as follows:
\begin{enumerate}
    \item Locate a set of eigenvectors of $(A,B)$ and compute the spectral projectors $P_R$ and $P_L$ onto the corresponding right and left deflating subspaces.
    \item Apply a rank-revealing factorization to find matrices $U_R$ and $U_L$ consisting of orthonormal columns that span $\text{range}(P_R)$ and $\text{range}(P_L)$, respectively.
    \item Reduce the problem $(A,B)$ to $(U_L^HAU_R, U_L^HBU_R)$ and repeat for the other set of eigenvectors, then recur.
\end{enumerate}

\indent One question remains: how do we locate a set of eigenvectors to project onto? Pseudospectral shattering provides an answer. If we have a grid that shatters the pseudospectrum of our perturbed and scaled problem $(\widetilde{A}, n^{\alpha} \widetilde{B})$, we can search for a grid line that significantly\footnote{Significant here means that at least a fifth of the eigenvalues lie on either side of the line.} splits the spectrum; applying a M\"{o}bius transformation that sends that grid line to the unit circle and then repeatedly squaring the pencil to drive one set of eigenvalues to zero, we obtain a projector onto the eigenvectors belonging to the surviving eigenvalues. In this approach, finding a grid line that produces a significant split ensures that divide-and-conquer will succeed in $O(\log(n))$ steps. \\
\indent For this to work without drawing a new grid at each step, we need shattering to be preserved as the problem shrinks in size. This is guaranteed by the following lemma. 
\begin{lem}\label{lem:exact projectors}
Let $P_R$ be a spectral projector for a regular matrix pencil $(A,B)$ and let $P_L$ be the projector onto the corresponding left deflating subspace. Let $U_L,U_R \in {\mathbb C}^{n \times k}$ be matrices whose orthonormal columns span the ranges of $P_L$ and $P_R$ respectively. Then
$$ \Lambda_{\epsilon} (U_L^HAU_R, U_L^HBU_R) \subseteq \Lambda_{\epsilon}(A,B) .$$
\end{lem}
\begin{proof}
If $z \in \Lambda_{\epsilon}(U_L^HAU_R, U_L^HBU_R)$ then there exists a unit vector $u \in {\mathbb C}^k$ such that 
\begin{equation}
    ||U_L^H(A-zB)U_Ru||_2 = ||(U_L^HAU_R - zU_L^HBU_R)u||_2 \leq \epsilon (1 + |z|) .
    \label{eqn: inequal one}
\end{equation}
Let $y = U_Ru \in {\mathbb C}^n$. Since $U_R$ has orthonormal columns $||y||_2 = ||U_Ru||_2 = ||u||_2 = 1$. Moreover, $y$ is in the right deflating subspace $\text{range}(P_R)$, which means $(A - zB)y$ belongs to the corresponding left deflating subspace $\text{range}(P_L)$. Since the columns of $U_L$ are an orthonormal basis for this subspace, we conclude
\begin{equation}
    ||(A-zB)y||_2 = ||U_L^H(A - zB)y||_2 \leq \epsilon( 1 + |z|)
    \label{eqn: inequal two}
\end{equation}
and therefore $z \in \Lambda_{\epsilon}(A,B)$.
\end{proof}
\cref{lem:exact projectors} emphasizes the importance of the matrices $U_L$ and $U_R$. While in principle many of the building blocks of divide-and-conquer hold for arbitrary matrices whose columns span deflating subspaces, its proof requires properties of matrices with orthonormal columns, meaning it would not be sufficient to simply compute the projectors $P_R$ and $P_L$. On top of this, $U_L$ and $U_R$ have the added benefit of ensuring that the norms of $U_L^HAU_R$ and $U_LBU_R$ don't grow, allowing us to carry norm assumptions through the recursion. 

\subsection{Numerical Building Blocks}
In this section, we review the numerical building blocks that allow us to compute the projectors $P_R$ and $P_L$ and subsequently the orthonormal bases contained in $U_L$ and $U_R$. Note that, as mentioned in the introduction, each of these uses only QR and matrix multiplication.
\subsubsection{Implicit Repeated Squaring (IRS)}
We begin with Implicit Repeated Squaring (\textbf{IRS}), a routine for repeatedly squaring a product $A^{-1}B$ without forming it. \textbf{IRS} first appeared in the context of eigenvalue problems in initial explorations of inverse-free and parallel methods \cite{MALYSHEV, Bai:CSD-94-793}. It was subsequently stated as it appears in \cref{alg:IRS} under the name \textbf{IRS} in a technical report of Ballard, Demmel, and Dumitriu \cite{Ballard2010MinimizingCF}. 

\begin{algorithm}
\caption{Implicit Repeated Squaring (\textbf{IRS})\\
\textbf{Input:} $A, B \in {\mathbb C}^{n \times n}$ and $p$ a positive integer \\
\textbf{Output:} $A_p, B_p \in {\mathbb C}^{n \times n}$ satisfying $A_p^{-1}B_p = (A^{-1}B)^{2^p}$ }\label{alg:IRS}
\begin{algorithmic}[1]
\State $A_0 = A$
\State $B_0 = B$
\For{$j = 0:p-1$}
    \State $$ \aligned \begin{pmatrix} B_j \\ - A_j 
    \end{pmatrix}  &= \begin{pmatrix} Q_{11} & Q_{12} \\
    Q_{21} & Q_{22} \end{pmatrix} \begin{pmatrix} R_{j} \\ 0  \end{pmatrix} \\
    A_{j+1} & = Q_{12}^H A_j \\
    B_{j+1} & = Q_{22}^H B_j  \\
    \endaligned $$
\EndFor
\State \Return $A_p, B_p$
\end{algorithmic}
\end{algorithm}

\indent Assuming invertibility where needed, we can verify the guarantee of \cref{alg:IRS} by noting that line 4 ensures $Q_{12}^HB_j - Q_{22}^HA_j = 0$ and consequently $B_jA_j^{-1} = Q_{12}^{-H}Q_{22}^H$, which implies
\begin{equation}
    A_{j+1}^{-1} B_{j+1} = A_j^{-1}Q_{12}^{-H}Q_{22}^HB_j = A_j^{-1}B_jA_j^{-1}B_j = (A_j^{-1}B_j)^2 .
    \label{eqn: IRS check}
\end{equation}
\indent The connection to spectral projectors is fairly straightforward:\ if $A$ is nonsingular, $A^{-1}B$ and $(A,B)$ have the same set of eigenvectors, with each finite eigenvalue $\lambda$ of $A^{-1}B$ corresponding to an eigenvalue $\lambda^{-1}$ of $(A,B)$. Repeatedly squaring $A^{-1}B$ drives its eigenvalues to either zero or infinity, assuming $(A,B)$ has no eigenvalues on the unit circle, without changing eigenvectors. Moreover, after $p$ steps of repeated squaring we have
\begin{equation}
    (A_p + B_p)^{-1} A_p = (I + A_p^{-1}B_p)^{-1} = (I + (A^{-1}B)^{2^p})^{-1},
    \label{eqn: IRS check 2}
\end{equation}
so the eigenvalues of $(A_p + B_p)^{-1}A_p$ are close to either zero or one, corresponding to eigenvalues that were originally inside, respectively outside, the unit disk. For clarity, \cref{tab:irs eig} shows the correspondence between eigenvalues of $(A,B)$
and eigenvalues of the outputs of \textbf{IRS}. \\
\indent Since $(A_p + B_p)^{-1}A_p$ has the same set of eigenvectors as $A^{-1}B$, it is not hard to show that $(A_p + B_p)^{-1}A_p$ approaches a projector $P_{R, |z|> 1}$ onto the right deflating subspace of $(A,B)$ spanned by eigenvectors with eigenvalues outside the unit disk \cite[\S 4]{Bai:CSD-94-793}. Nevertheless, for divide-and-conquer to proceed we also need a projector onto the corresponding left deflating subspace. Luckily, this can be done by running $\textbf{IRS}$ on $(A^H,B^H)$ since taking a Hermitian transpose swaps right and left deflating spaces. In fact, if ${\mathcal A}_p$ and ${\mathcal B}_p$ are the outputs of $p$ steps of repeated squaring on $(A^H, B^H)$, then
\begin{equation}
    {\mathcal A}_p^H({\mathcal A}_p + {\mathcal B}_p)^{-H} = (I + (A^{-H}B^H)^{2^p})^{-H}
    \label{eqn: IRS check 3}
\end{equation}
approaches a projector $P_{L, |z|>1}$ onto the left deflating subspace corresponding to $P_{R, |z|>1}$. And while taking a Hermitian transpose does change the eigenvalues, we can see from \cref{tab:irs eig} that repeated squaring on $(A^H, B^H)$ drives corresponding eigenvalues to the same destination (either zero or one). \\
\indent As mentioned, \textbf{IRS} breaks down if the input pencil has an eigenvalue on the unit circle. Considering such a problem to be ill-posed, Bai, Demmel, and Gu introduced the following distance $d_{(A,B)}$, whose reciprocal can be thought of as a condition number for the procedure.
\begin{defn}\label{ill-posed}
The distance from $(A,B)$ to the nearest ill-posed problem is 
$$d_{(A,B)} = \inf \left\{ ||E|| + ||F|| : (A + E) - z(B + F) \; \; \text{is singular for some} \; \; |z| =1 \right\}.  $$ 
\end{defn}
\indent Using $d_{(A,B)}$, they proved quadratic convergence for exact arithmetic \textbf{IRS}. Informally, this result states that to reach accuracy $\delta$ we need to take roughly $ O(\log_2 \left( 1/\delta \right))$ steps of repeated squaring.
\begin{thm}[Bai-Demmel-Gu 1994]\label{thm: IRS error} Let $A_p,B_p$ be the result of applying \text{\normalfont \bf IRS} to $A,B$. If
$$ p \geq \log_2 \left[ \frac{ ||(A,B)||_2 - d_{(A,B)}}{d_{(A,B)}} \right]  $$
then
$$ ||(A_p + B_p)^{-1}A_p - P_{R, |z|>1}||_2 \leq ||P_{R,|z|>1}||_2 \frac{2^{p+3} \left( 1 - \frac{d_{(A,B)}}{||(A,B)||_2} \right)^{2^p}}{\max \left\{ 0, 1 - 2^{p+2} \left( 1 -\frac{d_{(A,B)}}{||(A,B)||_2} \right)^{2^p} \right\}}.$$
\end{thm}
\indent Note that because $d_{(A^H, B^H)} = d_{(A,B)}$  and $||(A,B)||_2 = ||(A^H, B^H)||_2$, \cref{thm: IRS error} can be applied to the left projector by swapping $A$, $B$, $(A_p+B_p)^{-1}A_p$, and $P_{R, |z|>1}$ for $A^H$, $B^H$, $A_p^H(A_p+B_p)^{-H}$, and $P_{L,|z|>1}$, respectively.

\renewcommand{\arraystretch}{1.4}
\begin{table}[t]
    \centering
    \begin{tabular}{cccccccc}
        \hline
        \hline
         \multicolumn{8}{c}{Eigenvalue Correspondence in \textbf{IRS} (\cref{alg:IRS})} \\
         $(A,B)$ & $A^{-1}B$ & $A_p^{-1}B_p$ & $(A_p + B_p)^{-1}A_p$ & $(A^H, B^H)$ & $A^{-H}B^H$ & ${\mathcal A}_p^{-1}{\mathcal B}_p$ & ${\mathcal A}_p^H({\mathcal A}_p + {\mathcal B}_p)^{-H}$ \\
         \hline 
         $\lambda$ & $ \lambda^{-1} $ & $\lambda^{-2^p}$ & $(1 + \lambda^{-2 p})^{-1} $ & $\overline{\lambda}$ & $ ( \overline{\lambda})^{-1} $ & $(\overline{\lambda})^{-2p}$ & $(1 + \lambda^{-2p})^{-1} $  \\
          \hline 
          \hline
    \end{tabular}
    \caption{Correspondence between each eigenvalue $\lambda$ of $(A,B)$ and an eigenvalue of the matrices that appear when \textbf{IRS} is used to compute spectral projectors. Here, $[A_p, B_p] = \textbf{IRS}(A,B,p)$ while $[{\mathcal A}_p, {\mathcal B}_p] = \textbf{IRS}(A^H, B^H, p)$.}
    \label{tab:irs eig}
\end{table}

\subsubsection{RURV and GRURV}
Once we have our spectral projectors $P_R$ and $P_L$, we next need to apply a rank-revealing factorization to locate an orthonormal basis for the corresponding right and left deflating subspaces. To do this, we produce a (random) URV factorization of the projectors.
\begin{defn}
    $A = URV$ is a URV factorization of the matrix $A \in {\mathbb C}^{n \times n}$ if $R$ is upper triangular and  $U, V \in {\mathbb C}^{n \times n}$ are unitary.
    \label{defn: urv}
\end{defn}
The URV factorization was first introduced by  Stewart \cite{ULV} as an alternative to the SVD with many of the same properties. Most important is its rank-revealing capability:\ when $A$ has effective rank $k$ (i.e., there is a significant gap between $\sigma_k(A)$ and $\sigma_{k+1}(A)$), $A = URV$ is a rank-revealing factorization with 
\begin{equation}
    R = \begin{pmatrix} R_{11} & R_{21} \\ 0 & R_{22} \end{pmatrix}
    \label{eqn: rank-revealing urv}
\end{equation}
for $R_{11}$ a $k \times k$ block if $\sigma_k(R_{11})$ is a ``good" approximation to $\sigma_k(A)$ and  $\sigma_1(R_{22})$ is a ``good" approximation of $\sigma_{k+1}(A)$. We will make the meaning of ``good" precise in the guarantees to come.  \\
\indent While there are many ways to compute a rank-revealing URV factorization (see for example \cite[\S 3]{Fierro1999UTVTM}) we are interested in the \textbf{RURV} algorithm of Demmel, Dumitriu, and Holtz \cite{2007}, which is stated below as \cref{alg:RURV}. This randomized algorithm is simple to implement, backwards stable \cite[Theorem 4.5]{grurv}, and capable of producing strongly rank-revealing factorizations (in the sense of Gu and Eisenstat \cite{RRQR2}). In fact, Ballard et al.\ \cite{grurv} proved the following result, which demonstrates that \textbf{RURV} matches the best-known guarantees for deterministic factorizations.
\pagebreak

\begin{algorithm}
\caption{Randomized Rank-Revealing Factorization (\textbf{RURV})\\
\textbf{Input:} $A \in {\mathbb C}^{n \times n}$ \\
\textbf{Output:} $U$ unitary matrix, $R$ upper triangular matrix, and $V$ Haar such that $A = URV$ is a rank-revealing factorization of $A$. }\label{alg:RURV}
\begin{algorithmic}[1]
\State Draw a random matrix $B$ with i.i.d. ${\mathcal N}_{\mathbb C}(0,1)$ entries
\State $[V, \widehat{R}] = \textbf{QR}(B)$
\State $\widehat{A} = A \cdot V^H$
\State $[U,R] = \textbf{QR}(\widehat{A}) $
\State \Return $U,R,V$
\end{algorithmic}
\end{algorithm}

\begin{thm}[Ballard et al.\ 2019]
    Let $A \in {\mathbb C}^{n \times n}$ with singular values $\sigma_1, \sigma_2, \ldots, \sigma_n$. Let $R$ be the matrix produced by applying exact arithmetic \text{\normalfont \bf RURV} to $A$ as in \eqref{eqn: rank-revealing urv}. Assume that $k, n-k > 30$. Then with probability $1 - \delta$ the following occur:
    $$ \aligned
        \frac{\delta}{2.02} \frac{\sigma_k}{\sqrt{k(n-k)}} & \leq \sigma_{k}(R_{11}) \leq \sigma_k \\
        \sigma_{k+1} & \leq \sigma_{1}(R_{22}) \leq 2.02 \frac{\sqrt{k(n-k)}}{\delta} \sigma_{k+1} \\
        ||R_{11}^{-1}R_{12}||_2 & \leq \frac{6.1 \sqrt{k(n-k)}}{\delta} + \frac{\sigma_{k+1}}{\sigma_k} \frac{50 \sqrt{k^3(n-k)^3}}{\delta^3}. 
    \endaligned $$
    \label{thm: rurv guarantees}
\end{thm}

  The proof of \cref{thm: rurv guarantees} boils down to bounding the smallest singular value of a $k \times k$ block of a Haar unitary matrix via \cite[Theorem 5.2]{2007}. The requirement $k, n-k > 30$ comes from the bound used by Ballard et al.\ \cite[Corollary 3.4]{grurv}. Banks et al.\ subsequently demonstrated that this can be relaxed \cite[Proposition C.3]{banks2020pseudospectral}, although the fundamental guarantees are the same:\ with high probability \textbf{RURV} produces a factorization such that $\sigma_k(R_{11})$ and $\sigma_1(R_{22})$ are at worst a multiplicative factor of $O(\sqrt{k(n-k)})$ away from $\sigma_k(A)$ and $\sigma_{k+1}(A)$, respectively.\\
\indent Of course, the projectors we would like to apply \textbf{RURV} to are not simple matrices but rather products of the form $A^{-1}B$ or $AB^{-1}$ (i.e., the approximate projectors $(A_p + B_p)^{-1}A_p$ and ${\mathcal A}_p^H({\mathcal A}_p + {\mathcal B}_p)^{-H}$ from \textbf{IRS}). Since we avoid taking inverses, explicitly forming either of these products is not an option. Instead, a generalized version of \textbf{RURV} -- referred to as \textbf{GRURV} and presented here as \cref{alg:GRURV} -- allows us to apply \textbf{RURV} to an arbitrary product of matrices and their inverses. Note that in this routine \textbf{RULV} is a version of \textbf{RURV} that replaces the QR factorization in line 4 of \cref{alg:RURV} with QL.

\begin{algorithm}[t]
\caption{Generalized Randomized Rank-Revealing Factorization (\textbf{GRURV})\\
\textbf{Input:} $k$ a positive integer, $A_1, A_2, \ldots, A_k \in {\mathbb C}^{n \times n}$, and $m_1, m_2, \ldots, m_k \in \left\{ 1, -1 \right\} $ \\
\textbf{Output:} $U$ unitary, $R_1, R_2, \ldots, R_k$ upper triangular, and $V$ Haar such that  $UR_1^{m_1} R_2^{m_2} \cdots R_k^{m_k}V$ is a rank-revealing factorization of $A_1^{m_1}A_2^{m_2} \cdots A_k^{m_k}$ }\label{alg:GRURV}
\begin{algorithmic}[1]
\If{$m_k = 1$}
    \State $[U,R_k,V] = \textbf{RURV}(A_k)$
\Else 
    \State $[U, L_k, V] = \textbf{RULV}(A_k^H)$
    \State $R_k = L_k^H$
\EndIf
\State $U_{\text{current}} = U$
\For{$i = k-1: 1$} 
    \If{$m_i = 1$} 
        \State $[U, R_i] = \textbf{QR}(A_i \cdot U_{\text{current}}) $
        \State $U_{\text{current}} = U$
    \Else
        \State $[U, R_i] = \textbf{RQ}(U_{\text{current}}^H \cdot A_i)$
        \State $U_{\text{current}} = U^H$
    \EndIf
\EndFor
\State \Return $U_{\text{current}}$, optionally $R_1, R_2, \ldots, R_k$, $V$
\end{algorithmic}
\end{algorithm}

\indent \textbf{GRURV} was first introduced in a technical report of Ballard, Demmel, and Dumitriu \cite{Ballard2010MinimizingCF} specifically to apply \textbf{RURV} to spectral projectors found by \textbf{IRS}. Importantly, exact arithmetic \textbf{GRURV} is essentially equivalent to applying exact arithmetic \textbf{RURV} to the corresponding product \cite[Theorem 5.2]{grurv}, which allows us to access guarantees like \cref{thm: rurv guarantees} for \textbf{GRURV} in exact arithmetic without any additional effort.

\subsubsection{DEFLATE}
We are now ready to state a routine that combines \textbf{IRS} and \textbf{GRURV} to accomplish items (1) and (2) from our outline of divide-and-conquer in \cref{subsec: divide and conquer motivation}. Such a routine was first stated as \textbf{RGNEP} \cite[Algorithm 4]{Ballard2010MinimizingCF}, albeit in a different form than we present below. In particular, \textbf{RGNEP} assumed no knowledge of the number of eigenvalues of $(A,B)$ inside/outside the unit circle (equivalently, the rank of the corresponding spectral projectors), instead multiplying by the full $n \times n$ unitary matrices produced by \textbf{GRURV} and deciding where to split the problem to minimize certain matrix norms. Since we will have access to information about the rank of the projectors being computed (this is how dividing grid lines are selected), we state an alternative \textbf{DEFLATE} (\cref{alg:DEFLATE}), which simply takes the first $k$ columns of the matrices computed by \textbf{GRURV}.

\begin{algorithm}
\caption{Deflating Subspace Finder (\textbf{DEFLATE})\\
\textbf{Input:} $A,B \in {\mathbb C}^{n \times n}$, positive integers $p$ and $k$  \\
\textbf{Requires:} $k \leq n$; $(A,B)$ has no eigenvalues on the unit circle and exactly $k$ eigenvalues outside it. \\
\textbf{Output:} $U_R^{(k)}, U_L^{(k)} \in {\mathbb C}^{k \times n}$ with orthonormal columns that approximately span right and left deflating subspaces of $(A,B)$. }\label{alg:DEFLATE}
\begin{algorithmic}[1]
\State $[A_p, B_p] = \textbf{IRS}(A,B,p)$
\State $U_R = \textbf{GRURV}(2, A_p + B_p, A_p, -1, 1) $
\State $[A_p, B_p] = \textbf{IRS}(A^H, B^H, p)$
\State $U_L = \textbf{GRURV}(2, A_p^H, (A_p + B_p)^H, 1, -1)$
\State $U_R^{(k)} = U_R(: \; , 1:k) $
\State $U_L^{(k)} = U_L(: \; , 1:k) $
\State \Return $U_R^{(k)}$, $U_L^{(k)}$
\end{algorithmic}
\end{algorithm}

For \textbf{DEFLATE} to succeed, we need to know that the first $k$ columns of the U-factor produced by \textbf{GRURV} span the range of the rank-$k$ product it is applied to with high probability. We would also like a guarantee that the result for an approximate projector is close to the corresponding result for a true spectral projector (since \textbf{IRS} can only produce approximations of $P_R$ or $P_L$). \\
\indent If we apply \textbf{RURV} to a rank-$k$ matrix $A$, the first $k$ columns of the U-factor almost surely span $\text{range}(A)$. The intuition behind this is fairly simple:\ multiplying by the Haar matrix in line 3 of \textbf{RURV} ``mixes" the columns of $A$, distributing information so that the first $k$ columns of $A$ -- and therefore the first $k$ columns of $U$ -- are likely to span $\text{range}(A)$. On top of this, we have the following perturbation result in exact arithmetic, originally proved by Banks et al.\ \cite[Proposition C.12]{banks2020pseudospectral}.

\begin{thm}[Banks et al.\ 2022] \label{thm: Banks et al deflate guarnatee}
Let $A,A' \in {\mathbb C}^{n \times n}$ with $||A - A'||_2 \leq \delta$ and $\text{rank}(A) = \text{rank}(A^2) = k$. Let $T$ and $S$ contain the first $k$ columns of the U-factors produced by applying exact arithmetic \text{\normalfont \bf RURV} to $A$ and $A'$, respectively. Then for any $\theta \in (0,1)$ with probability $1 - \theta^2$ there exists a unitary $W \in {\mathbb C}^{k \times k}$ such that 
$$ ||S - TW^H||_2 \leq \sqrt{ \frac{8 \sqrt{k(n-k)}}{\sigma_k(T^HAT)}} \cdot \sqrt{ \frac{\delta}{\theta}} $$
\end{thm}

Letting $A$ be a rank-$k$ spectral projector (in which case $\sigma_k(T^HAT) = 1$), \cref{thm: Banks et al deflate guarnatee} says that the first $k$ columns of the U-factor of an approximation $A'$ of $A$ are close to a rotation/reflection of the first $k$ columns of the U-factor of $A$ as long as $||A-A'||_2$ is sufficiently small. \\ 
\indent Recalling that exact arithmetic \textbf{GRURV} is equivalent to exact arithmetic \textbf{RURV} on the corresponding product, these results generalize directly. Not only can we say that almost surely the first $k$ columns of the U-factor produced by \textbf{GRURV} span the range of the product, but a perturbation result similar to \cref{thm: Banks et al deflate guarnatee} holds. Combining these with our analysis of \textbf{IRS} yields the following exact arithmetic guarantee for \textbf{DEFLATE}.

\begin{thm}\label{thm: deflate}
Suppose $(A,B)$, $p$, and  $k$ satisfy the requirements of \text{\normalfont \bf DEFLATE}, where $p$ is large enough to ensure error in repeated squaring is at most $\delta$ in lines 1 and 3. Let $U_R^{(k)}$ and $U_L^{(k)}$ be the outputs of running this algorithm in exact arithmetic. Then for any $\nu \in (0,1)$, there exist $U_R,U_L \in {\mathbb C}^{n \times k}$ with orthonormal columns spanning right and left deflating subspaces of $(A,B)$, respectively, such that 
\begin{enumerate}
\item  $||U_R^{(k)} - U_R||_2 \leq  \sqrt{8 \sqrt{k(n-k)}} \sqrt{\frac{\delta}{\nu}}$
\item $ ||U_L^{(k)} - U_L||_2 \leq \sqrt{8 \sqrt{k(n-k)}} \sqrt{\frac{\delta}{\nu}} $
\end{enumerate}
with probability at least $1 - 2 \nu^2$.
\end{thm}
\begin{proof}
Consider first $U_R^{(k)}$ and let $A_p$ and $B_p$ be the outputs of applying $p$ steps of exact arithmetic repeated squaring to $A$ and $B$. We know $(A_p+B_p)^{-1}A_p$ approaches a projector $P_{R,|z|>1}$ onto the right deflating subspace spanned by eigenvectors of $(A,B)$ with eigenvalues outside the unit circle. Let $URV = P_{R,|z|>1}$ be a rank-revealing factorization of $P_{R,|z|>1}$ obtained via exact arithmetic \textbf{RURV} and let $U^{(k)}$ contain the first $k$ columns of $U$. Since $k$ is the number of eigenvalues with modulus greater than one, we know $k= \text{rank}(P_{R,|z|>1})$ and moreover $\text{range}(U^{(k)}) = \text{range}(P_{R, |z|>1})$ almost surely. Thus, since $p$ is large enough to ensure $||(A_p + B_p)^{-1}A_p - P_{R, |z|>1}||_2 \leq \delta$ and exact arithmetic \textbf{GRURV} satisfies the same guarantees as exact \textbf{RURV}, we have by \cref{thm: Banks et al deflate guarnatee} that with probability at least $1 - \nu^2$ there exists a unitary $W \in {\mathbb C}^{k \times k}$ such that 
\begin{equation} 
    ||U_R^{(k)} - U^{(k)} W^H||_2 \leq \sqrt{8 \sqrt{k(n-k)}} \sqrt{ \frac{\delta}{\nu}} .
    \label{eqn: exact deflate guarantee}
\end{equation}
Setting $U_R = U^{(k)}W^H$, we have inequality (1) with probability at least $1 - \nu^2$. Repeating the same argument for $U_L^{(k)}$, using this time the fact that $[A_p,B_p] = \textbf{IRS}(A^H, B^H, p)$ implies $||A_p^H(A_p + B_P)^{-H} - P_{L, |z|>1}||_2 \leq \delta$ for $P_{L,|z|>1}$ a projector onto the left deflating subspace corresponding to $P_{R,|z|>1}$, we obtain result (2) also with probability at least $1 - \nu^2$. Taking a union bound, we have both 1 and 2 with probability at least $1 - 2 \nu^2$.
\end{proof}

\subsection{Divide-and-Conquer Routine}
With \textbf{IRS}, \textbf{GRURV}, and \textbf{DEFLATE} (and their exact arithmetic guarantees discussed above) we can now state our divide-and-conquer eigensolver, presented below as \textbf{EIG} (\cref{alg:EIG}). The main inputs for this routine are a pencil $(A,B)$ and a grid $g$ that shatters $\Lambda_{\epsilon}(A,B)$ for some $\epsilon > 0$. In practice, this will be our perturbed and scaled pencil $(\widetilde{A}, n^{\alpha} \widetilde{B})$ (hence the norm assumption on $B$ has a factor of $n^{\alpha}$ attached) and the shattering grid guaranteed by \cref{thm: shattering}. Recall that shattering implies that $(A,B)$ is both regular and diagonalizable. 
\begin{algorithm}
\caption{Divide-and-Conquer Eigensolver (\textbf{EIG}) \\
\textbf{Input:} $n \in {\mathbb N}_{+}$, $A,B \in {\mathbb C}^{m \times m}$, $\epsilon  > 0$, $\alpha > 1$, $g$ an $s_1 \times s_2$ grid with box size $\omega$, $\beta > 0$ a desired eigenvector accuracy, and $\theta \in (0,1)$ a failure probability.\\
\textbf{Requires:} $m \leq n$, $||A||_2 \leq 3$, $||B||_2 \leq 3n^{\alpha}$, $g \subset \left\{ z : |\text{Re}(z)|,  |\text{Im}(z)| < 5 \right\}$, and $\Lambda_{\epsilon}(A,B)$ shattered w.r.t.\ $g$. \\
\textbf{Output:} $T$ an invertible matrix and $(D_1,D_2)$ a diagonal pencil. The eigenvalues of $(D_1,D_2)$ each share a grid box of $g$ with a unique eigenvalue of $(A,B)$ and each column of $T$ is an approximate right unit eigenvector of $(A,B)$.}\label{alg:EIG} 
\begin{algorithmic}[1]
\If{$m=1$}
    \State $T = 1$; $D_1 = A$; $D_2 = B$
\Else
    \State $\zeta = 2 \left( \lfloor \log_2(\max \left\{s_1, s_2 \right\})  + 1 \rfloor\right)  $
    \State $\eta = \min \left\{ \frac{  4\pi}{315 \sqrt{8}} \frac{ \beta \epsilon^2}{\omega n^{\alpha}}, \; \frac{1}{2 \log_{5/4}(n)} \right\}$
    \State $\delta = \min \left\{ \sqrt{ \frac{\theta}{10}} \frac{\epsilon^2}{7200n^{2 \alpha +3 }}, \; \frac{\theta}{2(\theta + 10n^6 \zeta)}, \; \sqrt{\frac{\theta}{10}} \frac{ \eta^2}{288n^{2 \alpha + 3}} \right\}   $
    \State $p = \left\lceil \max \left\{ 7, \; \log_2 \left( \frac{105n^{\alpha}}{\epsilon} - 1 \right), \; -2 \log_2 \left( -\frac{1}{2} \log_2 \left( 1 - \frac{\epsilon}{105n^{\alpha}} \right) \right), \; 1 + \log_2 \left[ \frac{ \log_2 \left( \frac{\delta \pi \epsilon}{12n^{\alpha} m \omega + \delta \pi \epsilon} \right)}{\log_2 \left( 1 - \frac{\epsilon}{105n^{\alpha}}\right)} \right] \right\} \right\rceil$
    \State Choose a grid line $\text{Re}(z) = h$ of $g$
    \State $({\mathcal A}, {\mathcal B}) = (A - (h-1)B, A - (h+1)B)$
    \State $[A_p,B_p] = \textbf{IRS}({\mathcal A}, {\mathcal B}, p)$
    \State $[U,R_1,R_2, V] = \textbf{GRURV}(2, A_p + B_p, A_p, -1, 1)$
    \State $k = \# \left\{ i : \left| \frac{R_2(i,i)}{R_1(i,i)} \right| \geq \sqrt{\frac{\theta}{10 \zeta}}\frac{1 - \delta}{n^3} \right\}  $
    \If{$k < \frac{1}{5}m$ or $k > \frac{4}{5}m$} 
        \State Return to step 8 and choose a new grid line, executing a binary search if necessary. If this fails, \indent  \; \; search over horizontal grid lines $\text{Im}(z) = h$, this time setting $({\mathcal A}, {\mathcal B}) = (A - i(h-1)B, A - i(h+1)B)$.
    \Else 
        \State $[U_R^{(k)}, U_L^{(k)}] = \textbf{DEFLATE}({\mathcal A}, {\mathcal B}, p, k)$
        \State $  ({\mathcal A}, {\mathcal B}) = (A - (h+1)B, A - (h-1)B) $ \Comment{or  $({\mathcal A}, {\mathcal B}) = (A - i(h+1)B, A - i(h-1)B)$}
        \State $[U_R^{(m-k)}, U_L^{(m-k)}] = \textbf{DEFLATE}({\mathcal A}, {\mathcal B}, p, m-k) $
            \State $$ \aligned 
            (A_{11}, B_{11}) &= \left( (U_L^{(k)})^H A U_R^{(k)}, \; (U_L^{(k)})^H B U_R^{(k)} \right)  \\
         (A_{22}, B_{22}) &= \left( (U_L^{(m-k)})^HAU_R^{(m-k)}, \; (U_L^{(m-k)})^HBU_R^{(m-k)} \right) \\
        \endaligned $$
        \State $g_R = \left\{ z \in g : \text{Re}(z) > h \right\}$ \Comment{or $g_R = \left\{ z \in g : \text{Im}(z) > h \right\}$}
        \State $[\widehat{T}, \widehat{D}_1, \widehat{D}_2] = \textbf{EIG}(n, A_{11}, B_{11}, \frac{4}{5} \epsilon, \alpha, g_R, \frac{1}{3} \beta, \theta)$
        \State $g_L = \left\{ z \in g : \text{Re}(z) < h \right\}$ \Comment{or $g_L = \left\{ z \in g : \text{Im}(z) < h \right\}$}
        \State $[\widetilde{T}, \widetilde{D}_1, \widetilde{D}_2] = \textbf{EIG}(n, A_{22}, B_{22}, \frac{4}{5} \epsilon, \alpha, g_L, \frac{1}{3} \beta, \theta)$
    \State $$ \aligned 
    &T = \begin{pmatrix} U_R^{(k)} & U_R^{(m-k)} \end{pmatrix} \begin{pmatrix} \widehat{T} & 0 \\ 0 & \widetilde{T} \end{pmatrix} \\
    & D_1 = \begin{pmatrix} \widehat{D}_1 & 0 \\ 0 & \widetilde{D}_1 \end{pmatrix}\\
    & D_2 = \begin{pmatrix} \widehat{D}_2 & 0 \\ 0 & \widetilde{D}_2 \end{pmatrix} 
    \endaligned $$
\EndIf
\EndIf
\State \Return $T, D_1, D_2$
\end{algorithmic}
\end{algorithm}

\indent Before proving exact arithmetic guarantees for \textbf{EIG}, we first provide a high-level overview of the algorithm to motivate why it works. Throughout, we rely on the analysis of each of the building blocks from the previous section as well as the requirements listed in \cref{alg:EIG}.
\begin{enumerate}
    \item Since \textbf{EIG} calls itself recursively, the first three lines check for our stopping criteria. We choose to continue divide-and-conquer until the pencil is $1 \times 1$, though as mentioned in the introduction we could choose instead to stop once the pencil is small enough to be handled by another method.
    \item The next four lines (4-7) set parameters for the algorithm. Most importantly, they determine how many steps of repeated squaring need to be taken to achieve the desired accuracy (i.e., the value of $p$ in line 7).
    \item Lines 8-15 execute a search over $g$ for a grid line that sufficiently splits the spectrum, which in our case means separating at least a fifth of the eigenvalues on each side. Since $g$ shatters $\Lambda_{\epsilon}(A,B)$, a grid line that sufficiently splits the spectrum always exists. 
    \item We check a line $\text{Re}(z) = h$ of the grid by applying the M\"{o}bius transformation $S(z) = \frac{z - (h-1)}{z - (h+1)}$ to $(A,B)$ (line 9). $S$ maps the grid line to the unit circle while sending the half plane $\left\{ \text{Re}(z) < h \right\}$ inside the unit disk. Applying this transformation to $(A,B)$ sends eigenvalues to the left/right of the dividing line inside/outside the unit circle, respectively, without changing eigenvectors. 
    \item Lines 10 and 11 apply \textbf{IRS} and \textbf{GRURV} to the transformed pencil $({\mathcal A}, {\mathcal B})$. This produces a rank-revealing factorization $UR_1^{-1}R_2V$ of the approximate projector onto the right deflating subspace corresponding to eigenvectors of $({\mathcal A}, {\mathcal B})$ with eigenvalues outside the unit disk (equivalently eigenvectors of $(A,B)$ with eigenvalues to the right of the selected grid line).
    \item In line 12, we leverage the rank-revealing guarantees of \textbf{RURV}, and by extension \textbf{GRURV}, to read off the rank of the approximate projector. Note that we do this without forming $R_1^{-1}R_2$. The grid line is selected if this rank is between $\frac{1}{5}m$ and $\frac{4}{5}m$, where $m$ is the size of the pencil (which shrinks as we recur). 
    \item In line 8 we assume that the grid line is vertical, however it is possible that only a horizontal grid line sufficiently splits the spectrum. This is covered in line 14. The remainder of the algorithm similarly assumes the split is vertical, though we note in lines 17, 20, and 22 the adjustments that need to be made if we split by $\text{Im}(z) = h$ instead of $\text{Re}(z) = h$.
    \item Once a dividing line is identified, we call \textbf{DEFLATE} twice to compute orthonormal bases for both sets of spectral projectors. To recover eigenvectors corresponding to eigenvalues to the left of the line, we apply the alternative M\"{o}bius transformation $S(z) = \frac{z - (h+1)}{z - (h-1)}$.
    \item In line 19  we compute the next pair of subproblems. We then pass these to \textbf{EIG} along with pieces of the grid $g$ and slightly adjusted parameters. In particular, note that the $\epsilon$ for which $\Lambda_{\epsilon}(A,B)$ is shattered shrinks by a factor of $\frac{4}{5}$ at each step. As we will see, this is necessary to guarantee shattering since $U_R^{(k)}$, $U_L^{(k)}$, $U_R^{(m-k)}$ and $U_L^{(m-k)}$ are only approximations of the matrices used in \cref{lem:exact projectors}.
    \item Once the recursion finishes, \textbf{EIG} reconstructs a diagonal pencil $(D_1,D_2)$ and a set of approximate right eigenvectors $T$ (line 24).
\end{enumerate}
With this outline in mind, we are now ready to state and prove our main guarantee for \textbf{EIG}. 
\begin{thm}\label{thm: EIG succeeds}
Let $(A,B)$ and $g$ be a pencil and grid satisfying the requirements of \text{\normalfont \bf EIG}. Then for any choice of $\theta \in (0,1)$ and $\beta > 0$, exact-arithmetic \text{\normalfont \bf EIG} applied to $(A,B)$ and $g$ satisfies the following with probability at least $1 - \theta$.
\begin{enumerate}
    \item The recursive procedure converges and each eigenvalue of the diagonal pencil $(D_1,D_2)$ shares a grid box with a unique eigenvalue of $(A,B)$.
    \item If $\sigma_n(B) \geq 1$, each column $t_i$ of $T$ satisfies $ ||t_i - v_i||_2 \leq \beta $ for some right unit eigenvector $v_i$ of $(A,B)$.
\end{enumerate}
\end{thm}
\begin{proof} We start by bounding the probability that the first guarantee does not hold. Since \textbf{EIG} calls itself recursively, we do this by bounding the probability of failure for one step of divide-and-conquer. In this context, success requires two events:\ first, a dividing line that sufficiently splits the spectrum must be found; second, the subsequent calls to $\textbf{EIG}$ must be valid, meaning the inputs satisfy the listed properties. \\
\indent Computing the probabilities that these occur is fairly lengthy, so to improve readability we number the steps in the proof and provide in bold a description of what each step accomplishes. Throughout, we use the assumptions on the inputs -- i.e., $A,B \in {\mathbb C}^{m \times m}$ with $m \leq n$, $||A||_2 \leq 3$, $||B||_2 \leq 3 n^{\alpha}$, and $\Lambda_{\epsilon}(A,B)$ is shattered with respect to the grid $g$, which is $s_1 \times s_2$ consisting of boxes of size $\omega$.
\begin{enumerate}
    \item[\textbf{1.}] \textbf{Any transformed pencil $({\mathcal A}, {\mathcal B})$ in EIG satisfies $d_{({\mathcal A}, {\mathcal B})} \geq \frac{2}{5} \epsilon$.} \\
    Consider first a vertical grid line $\text{Re}(z) = h$ and $({\mathcal A}, {\mathcal B}) = (A - (h-1)B, A-(h+1)B)$ as in line 9. Suppose $z \in \Lambda_{\epsilon'}({\mathcal A}, {\mathcal B})$ for some $\epsilon ' > 0$. In this case, there exist matrices $E$ and $F$ with $||E||_2, ||F||_2 \leq \epsilon'$ such that $z$ is an eigenvalue of $({\mathcal A} + E, \; {\mathcal B} + F)$. If we apply the M\"{o}bius transformation $S(z) = \frac{(h+1)z - (h-1)}{z -1}$ to this pencil, we observe that $S(z)$ is an eigenvalue of
    \begin{equation}
        ((h+1)({\mathcal A} + E) - (h-1)({\mathcal B}+F), \; ({\mathcal A} + E) - ({\mathcal B} + F)), 
        \label{eqn: exact step one part one}
    \end{equation}
    or, equivalently,
    \begin{equation} 
        (2A + (h+1)E - (h-1)F, \; 2B + E - F) .
        \label{eqn: exact step one part two}
    \end{equation}
    Dividing by two, we conclude that $S(z)$ is an eigenvalue of $(A + \frac{h+1}{2}E - \frac{h-1}{2}F, \; B + \frac{1}{2} E - \frac{1}{2} F)$, where
    \begin{equation} 
        \left| \left| \frac{h+1}{2}E - \frac{h-1}{2}F \right| \right|_2 \leq \frac{|h+1|}{2} || E||_2 +  \frac{|h-1|}{2}||F||_2 \leq \frac{\epsilon'}{2} \left( |h+1| + |h-1| \right) 
        \label{eqn: exact step one part three}
    \end{equation}
    and 
    \begin{equation}
        \left| \left| \frac{1}{2} E - \frac{1}{2} F \right| \right| _2 \leq \frac{1}{2} (|| E ||_2 + ||F ||_2) \leq \epsilon' .
        \label{eqn: exact step one part four}
    \end{equation}
    Thus, $S(z)$ belongs to $\Lambda_{\epsilon''}(A,B)$ for 
    \begin{equation} 
        \epsilon'' = \max \left\{ \epsilon', \; \frac{\epsilon'}{2}(|h+1| + |h-1|) \right\} \leq 5\epsilon',
        \label{eqn: exact step one part five}
    \end{equation}
    which means the preimage of $\Lambda_{\epsilon / 5}({\mathcal A}, {\mathcal B})$ under $S^{-1}$ is contained in $\Lambda_{\epsilon}(A,B)$. Since $\Lambda_{\epsilon}(A,B)$ is shattered with respect to $g$ and therefore does not intersect the dividing line $\text{Re}(z) = h$, we conclude that $\Lambda_{\epsilon / 5}({\mathcal A}, {\mathcal B})$ does not intersect the unit circle. By \cref{ill-posed}, we obtain $d_{({\mathcal A}, {\mathcal B})} \geq \frac{2}{5} \epsilon$. Making a similar argument for the transformed pencil in line 17 or in the case of a horizontal dividing line $\text{Im}(z) = h$ yields $d_{({\mathcal A}, {\mathcal B})} \geq \frac{2}{5} \epsilon$ for any $({\mathcal A}, {\mathcal B})$ appearing in \textbf{EIG}. In the next step, we will use this lower bound to control the error in repeated squaring. 
    \item[\textbf{2.}] \textbf{The choice of $p$ guarantees that the error in repeated squaring is at most $\delta$.} \\
    Consider the first call to \textbf{IRS}, which applies repeated squaring to the transformed pencil $({\mathcal A}, {\mathcal B})$. By \cref{thm: IRS error}, we know that as long as $p \geq \log_2 \left[ \frac{||({\mathcal A}, {\mathcal B})||_2 - d_{({\mathcal A}, {\mathcal B})}}{d_{({\mathcal A}, {\mathcal B})}} \right] $ then
    \begin{equation}
    ||(A_p + B_p)^{-1}A_p - P_{R, |z|>1}||_2 \leq ||P_{R,|z|>1}||_2 \frac{2^{p+3} \left( 1 - \frac{d_{({\mathcal A},{\mathcal B})}}{||({\mathcal A},{\mathcal B})||_2} \right)^{2^p}}{\max \left\{ 0, 1 - 2^{p+2} \left( 1 -\frac{d_{({\mathcal A},{\mathcal B})}}{||({\mathcal A},{\mathcal B})||_2} \right)^{2^p} \right\}}.
    \label{eqn: IRS guarantee}
    \end{equation}
    We just showed $d_{({\mathcal A}, {\mathcal B})} \geq \frac{2}{5} \epsilon$ and 
    \begin{equation}
        ||({\mathcal A}, {\mathcal B})||_2 \leq ||{\mathcal A}||_2 + || {\mathcal B}||_2 \leq 2||A||_2 + (|h-1|+|h+1|)||B||_2 \leq 42n^{\alpha},
        \label{eqn: exact step two norm bound}
    \end{equation}
    so to satisfy $p \geq \log_2 \left[ \frac{||({\mathcal A}, {\mathcal B})||_2 - d_{({\mathcal A}, {\mathcal B})}}{d_{({\mathcal A}, {\mathcal B})}} \right] $ it is sufficient to take $p \geq \log_2 \left( \frac{105n^{\alpha}}{\epsilon} - 1 \right) $. Similarly, to eliminate the maximum from the denominator of \eqref{eqn: IRS guarantee} it is sufficient to take $p \geq -2 \log_2 \left( \log_2 \left( \frac{105n^{\alpha}}{105n^{\alpha} - \epsilon} \right) \right) $ provided $p > 6$ (which allows us to simplify the bounds by assuming $\log_2(p+2) < \frac{1}{2}p$). \\
    
    With this, we can now turn to bounding the right hand side of \eqref{eqn: IRS guarantee}. First, we upper bound $||P_{R,|z|>1}||_2$. Recall,
    \begin{equation}
        P_{R,|z|>1} = V \begin{pmatrix} 0 & 0 \\ 0 & I_r \end{pmatrix} V^{-1} = \sum_{j=m-r+1}^m v_jw_j^H
        \label{eqn: exact step two projector def}
    \end{equation}
    for $r = \text{rank}(P_{R,|z|>1})$ and $V$ a matrix that diagonalizes $B^{-1}A$. $v_i$ and $w_i^H$ are the columns of $V$ and rows of $V^{-1}$ respectively, scaled so that $w_i^Hv_i = 1$ with $v_{m-r+1}, \ldots, v_m$ corresponding to eigenvalues of $(A,B)$ to the right of $\text{Re}(z) = h$. Since $\Lambda_{\epsilon}(A,B)$ is shattered with respect to $g$, each of these eigenvalues $\lambda_{m-r+1}, \ldots, \lambda_m$ is contained in a separate grid box of $g$. If $\Gamma_i$ is the contour of the grid box containing $\lambda_i$, this means 
    \begin{equation} 
        v_jw_j^H = \frac{1}{2\pi i} \oint_{\Gamma_j} (z - B^{-1}A)^{-1} dz 
        \label{eqn: exact step two eigenvalue condition integral}
    \end{equation}
    and therefore
    \begin{equation} 
        P_{R, |z|>1} = \frac{1}{2 \pi i} \sum_{j=m-r+1}^m \oint_{\Gamma_j} (z-B^{-1}A)^{-1}dz = \frac{1}{2 \pi i} \sum_{j=m-r+1}^m  \oint_{\Gamma_j}(zB - A)^{-1}Bdz .
        \label{eqn: exact step two projector integral}
    \end{equation}
    Thus, by the triangle inequality,
    \begin{equation}
        ||P_{R, |z|>1}||_2 \leq \frac{1}{2 \pi} \sum_{j=m-r+1}^m \left| \left| \oint_{\Gamma_j} (zB-A)^{-1}B dz \right| \right|_2. 
        \label{eqn: exact step two apply triangle}
    \end{equation}
    Moreover, applying the ML inequality to each term in this sum, we have
    \begin{equation} 
        ||P_{R,|z|>1}||_2 \leq \frac{1}{2 \pi} \sum_{j=m-r+1}^m 4 \omega \sup_{z \in \Gamma_j} ||(zB-A)^{-1}B||_2 \leq \frac{2 \omega ||B||_2}{ \pi} \sum_{j=m-r+1}^m \sup_{z \in \Gamma_j} ||(A - zB)^{-1}||_2 .
        \label{eqn: exact step two ML}
    \end{equation}
    Since shattering guarantees $\Lambda_{\epsilon}(A,B) \cap \Gamma_j = \emptyset$ and therefore $||(A - zB)^{-1}||_2 \leq \frac{1}{\epsilon (1 + |z|)} \leq \frac{1}{\epsilon}$ for all $z \in \Gamma_j$, we conclude $||P_{R,|z|>1}||_2 \leq \frac{2 \omega r ||B||_2 }{ \pi \epsilon} $. Finally using the fact that $||B||_2 \leq 3 n^{\alpha}$ and $r \leq m$, we have a final upper bound $||P_{R,|z|>1}||_2 \leq  \frac{6n^{\alpha}m \omega}{ \pi \epsilon} $. 
    
    Combining this bound with $d_{({\mathcal A}, {\mathcal B})} \geq \frac{2}{5} \epsilon$ and $||({\mathcal A}, {\mathcal B})||_2 \leq 42n^{\alpha}$, \eqref{eqn: IRS guarantee} becomes  
    \begin{equation} 
        ||(A_p + B_p)^{-1}A_p - P_{R,|z|>1}||_2 \leq \frac{6n^{\alpha} m \omega}{\pi \epsilon} \cdot \frac{ 2^{p+3} \left( 1 - \frac{\epsilon}{105n^{\alpha}} \right)^{2^p}}{ 1 - 2^{p+2} \left( 1 - \frac{\epsilon}{105n^{\alpha}} \right)^{2^p} }
        \label{eqn: exact step two refined upper bound}
    \end{equation}
    for $p$ sufficiently large (i.e., following the bounds derived above). Thus, we obtain $||(A_p + B_p)^{-1}A_p - P_{R,|z|>1}||_2 \leq \delta$ by taking 
    \begin{equation} 
        \frac{6n^{\alpha}m \omega}{\pi \epsilon} \cdot \frac{ 2^{p+3} \left( 1 - \frac{\epsilon}{105n^{\alpha}} \right)^{2^p}}{ 1 - 2^{p+2} \left( 1 - \frac{\epsilon}{105n^{\alpha}} \right)^{2^p} } \leq \delta
        \label{eqn: exact step two sufficient choice delta}
    \end{equation}
    which is equivalent to 
    \begin{equation}
        2^p \left[ \frac{p+2}{2^p} +  \log_2 \left( 1 - \frac{\epsilon}{105n^{\alpha}} \right) \right] \leq \log_2 \left( \frac{\delta \pi \epsilon}{12n^{\alpha} m \omega + \delta \pi \epsilon} \right) .
        \label{eqn: exact step two delta choice two}
    \end{equation}
    Using again the assumption that $p > 6$ and further taking $p \geq -2 \log_2 \left( - \frac{1}{2} \log_2 \left( 1 - \frac{\epsilon}{105n^{\alpha}} \right) \right)$ to ensure $\frac{p+2}{2^p} \leq - \frac{1}{2} \log_2 \left( 1 - \frac{\epsilon}{105n^{\alpha}} \right) $, we get the desired accuracy as long as
    \begin{equation} 
        2^{p-1} \log_2 \left( 1 - \frac{\epsilon}{105 n^{\alpha}} \right) \leq \log_2 \left( \frac{\delta \pi \epsilon}{12n^{\alpha} m \omega + \delta \pi \epsilon} \right)
        \label{eqn: exact step two delta choice three}
    \end{equation}
    which yields a final bound $p \geq 1 + \log_2 \left[\log_2 \left( \frac{\delta \pi \epsilon}{12n^{\alpha} m \omega + \delta \pi \epsilon} \right) / \log_2 \left( 1 - \frac{\epsilon}{105n^{\alpha}} \right) \right] $.
    
    In the preceding analysis, we derived the following four bounds on $p$:
    \begin{itemize}
        \item $p \geq \log_2 \left( \frac{105n^{\alpha}}{\epsilon} - 1\right)$ [to allow us to apply \cref{thm: IRS error}].
        \item $p \geq -2 \log_2 \left( \log_2 \left( \frac{105n^{\alpha}}{105n^{\alpha} - \epsilon} \right) \right) $ [to eliminate the maximum from the error bound].
        \item $p \geq -2 \log_2 \left( - \frac{1}{2} \log_2 \left( 1 - \frac{\epsilon}{105n^{\alpha}} \right) \right) $ [to simplify the upper bound in \eqref{eqn: IRS guarantee}].
        \item $ p \geq 1 + \log_2 \left[ \log_2 \left( \frac{\delta \pi \epsilon}{12n^{\alpha} m \omega + \delta \pi \epsilon} \right) / \log_2 \left( 1 - \frac{\epsilon}{105n^{\alpha}} \right) \right] $ [to ensure an error of at most $\delta$ given the other three bounds].
    \end{itemize}
    Since the third bound is always at least as large as the second and since we also assumed $p > 6$, we conclude that $||(A_p + B_p)^{-1}A_p - P_{R,|z|>1}||_2 \leq \delta$ for the $p$ chosen in line 7 of \textbf{EIG}. Note that since $||{\mathcal A}^H||_2 = ||{\mathcal A}||_2$, $||{\mathcal B}^H||_2 = ||{\mathcal B}||_2$, $d_{({\mathcal A}^H, {\mathcal B}^H)} = d_{({\mathcal A}, {\mathcal B})}$, and $\Lambda_{\epsilon}(A^H,B^H)$ is shattered with respect to the grid $g^H = \left\{ \overline{z} : z \in g \right\}$, the same argument guarantees that running repeated squaring on ${\mathcal A}^H$, and $ {\mathcal B}^H$ has error at most $\delta$. Similarly, the same results hold for the other transformed pencils in lines 14 and 17 since in all cases $||({\mathcal A}, {\mathcal B})||_2 \leq 42 n^{\alpha}$ and $d_{({\mathcal A}, {\mathcal B})} \geq \frac{2}{5} \epsilon$.
    \item[\textbf{3.}] \textbf{A dividing line that sufficiently splits the spectrum exists.} \\
    A dividing line sufficiently splits the spectrum if it separates at least $\frac{1}{5}m$ of the $m$ eigenvalues of $(A,B)$. Suppose that no vertical line of $g$ does this. In this case, there exists adjacent vertical lines between which more than $\frac{3}{5}m$ eigenvalues lie. Since no eigenvalues share the same grid box, this implies that a horizontal grid line must sufficiently split the spectrum. 
    \item[\textbf{4.}] \textbf{With probability at least $1 - \frac{\theta}{10n^4}$, EIG finds a dividing line that separates exactly $k$ eigenvalues to the right such that $\frac{1}{5}m \leq k \leq \frac{4}{5} m$. } \\
    To obtain a lower bound on this probability, we first compute the probability that for any grid line $\text{Re}(z) = h$ the value of $k$ at line 12 is equal to the number of eigenvalues of $(A,B)$ to the right of the line. 
    
    Suppose $(A,B)$ has $r$ eigenvalues to the right of $\text{Re}(z) = h$. In this case, we know in line 10 that $(A_p+B_p)^{-1}A_p$ approaches a rank-$r$ projector $P_{R,|z|>1}$. Now $k$ is obtained by computing a rank-revealing factorization $U_RR_1^{-1}R_2V = (A_p+B_p)^{-1}A_p$ via \textbf{GRURV} and then counting the diagonal entries of $R_1^{-1}R_2$ that have modulus above a certain threshold. With this in mind, write
    \begin{equation} 
        R_1^{-1}R_2 = \begin{bmatrix} R_{11} & R_{12} \\ 0 & R_{22} \end{bmatrix}
        \label{exact step four single matrix R}
    \end{equation}
    for $R_{11}$ an $r \times r$ matrix. Since we are working in exact arithmetic, \textbf{GRURV} satisfies all of the guarantees of exact-arithmetic \textbf{RURV}. In particular, extending \cite[Theorem 5.2]{2007} to complex matrices, 
    \begin{equation}
        \sigma_r(R_{11}) \geq \sigma_r((A_p+B_p)^{-1}A_p) \sigma_r(X_{11})
        \label{eqn: R11 lower bound}
    \end{equation}
    where $X_{11}$ is the upper left $r \times r$ block of $X = Q^HV^H$ for the SVD $(A_p + B_p)^{-1}A_p = P \Sigma Q^H $. Similarly, by \cite[Lemma 4.1]{grurv},
    \begin{equation}
        ||R_{22}||_2 \leq \frac{\sigma_{r+1}((A_p + B_p)^{-1}A_p)}{\sigma_r(X_{11})} .
        \label{eqn: R22 bound}
    \end{equation}
    Now $||(A_p + B_p)^{-1}A_p - P_{R,|z|>1}||_2 \leq \delta$ as shown in step two above, so since the rank-$r$ projector $P_{R,|z|>1}$ satisfies $\sigma_r(P_{R,|z|>1}) = 1$ and $\sigma_{r+1}(P_{R,|z|>1}) = 0$, we have by the stability of singular values $ \sigma_r((A_p + B_p)^{-1}A_p) \geq 1 - \delta $ and $\sigma_{r+1}((A_p + B_p)^{-1}A_p) \leq \delta $. Moreover, since $X$ is Haar unitary, we have by \cite[Proposition C.3]{banks2020pseudospectral}
    \begin{equation} 
        {\mathbb P} \left[ \frac{1}{\sigma_r(X_{11})} \leq \frac{\sqrt{r(m-r)}}{\nu} \right] \geq 1 - \nu^2
        \label{eqn: exact step four Harr prob bound}
    \end{equation}
    for any $\nu \in (0,1]$. Applying these to \eqref{eqn: R11 lower bound} and \eqref{eqn: R22 bound} with $\nu = \sqrt{ \frac{\theta}{10 \zeta}}\frac{1}{n^2} $ for $\zeta = 2 ( \lfloor \log_2(s) + 1 \rfloor)$ and $s = \max \left\{ s_1, s_2 \right\}$, we have
    \begin{equation} 
        \sigma_r(R_{11}) \geq (1 - \delta) \sqrt{\frac{\theta}{10 \zeta}} \frac{1}{n^2 \sqrt{r(m-r)}} \geq \sqrt{ \frac{\theta}{10 \zeta}} \frac{1 - \delta}{n^3}
        \label{eqn: exact step four smallest singular value of R11}
    \end{equation}
    and
    \begin{equation}
        ||R_{22}||_2 \leq \delta \sqrt{r(m-r)} n^2 \sqrt{ \frac{10 \zeta}{\theta}} \leq \delta n^3 \sqrt{\frac{10 \zeta}{\theta}}
        \label{eqn: exact step four R22 norm bound}
    \end{equation}
    with probability at least $ 1 - \frac{\theta}{10 \zeta n^4}$. Since the eigenvalues of any matrix are bounded in modulus above and below by its singular values and the eigenvalues of $R_{11}$ and $R_{22}$ are their diagonal entries, we conclude $ |R_{11}(i,i)| \geq \sqrt{ \frac{\theta}{10 \zeta}} \frac{1 - \delta}{n^3}$ for all $1 \leq i \leq r$ while $ |R_{22}(j,j)| \leq \delta n^3 \sqrt{ \frac{10 \zeta}{\theta}}$ for all $1 \leq j \leq m-r$ with probability at least $1 - \frac{\theta}{10 \zeta n^4}$. Since the requirement in line 6 that $\delta \leq  \frac{\theta}{2(\theta + 10n^6 \zeta)}$ guarantees $\sqrt{\frac{\theta}{10 \zeta}} \frac{1 - \delta}{n^3} > \delta n^3 \sqrt{ \frac{10 \zeta}{\theta}}$, this implies that $k = r$ with probability at least $ 1- \frac{\theta}{10 \zeta n^4}$ 

    We have shown so far that with high probability the value of $k$ at line 12 is equal to the number of eigenvalues to the right of the vertical dividing line selected four lines earlier. Since repeating this argument yields the same probability of success when checking a horizontal grid line and we know a dividing line that separates at least $\frac{1}{5}m$ eigenvalues must exist, we can therefore lower bound the probability that \textbf{EIG} finds a suitable line by requiring that the value of $k$ is accurate for all lines that are checked. Since we do at most two binary searches to find a good enough grid line, we check at most $\zeta$ lines. By a union bound, we conclude that a suitable line is found, and $k$ is computed accurately for that line, with probability at least $ 1 - \frac{\theta}{10 n^4}$. 

    \item[\textbf{5.}] \textbf{Assuming $k$ is computed correctly in line 12, there exists matrices $\widehat{U}_R^{(k)}, \widehat{U}_L^{(k)} \in {\mathbb C}^{k \times m}$ and $\widehat{U}_R^{(m-k)}, \widehat{U}_L^{(m-k)} \in {\mathbb C}^{m-k \times m}$ with orthonormal columns spanning corresponding right and left deflating subspaces of $(A,B)$ such that
    \begin{enumerate}
        \item $||U_R^{(k)} - \widehat{U}_R^{(k)}||_2 \leq \sqrt{ \sqrt{ \frac{10}{\theta}} 8n^3 \delta} $
        \item $||U_L^{(k)} - \widehat{U}_L^{(k)}||_2 \leq \sqrt{ \sqrt{ \frac{10}{\theta}} 8n^3 \delta} $
        \item $||U_R^{(m-k)} - \widehat{U}_R^{(m-k)}||_2 \leq \sqrt{ \sqrt{ \frac{10}{\theta}} 8n^3 \delta} $
        \item $||U_L^{(m-k)} - \widehat{U}_L^{(m-k)}||_2 \leq \sqrt{ \sqrt{ \frac{10}{\theta}} 8n^3 \delta} $
    \end{enumerate}
    With probability at least $1 - \frac{2\theta}{5n^4}$.}  \\
    This result comes from applying \cref{thm: deflate} twice with $\nu = \sqrt{\frac{\theta}{10n^4}}$ and taking a union bound. In both cases, we use the fact that $p$ is large enough to guarantee error in repeated squaring is a most $\delta$ and consequently 
    \begin{equation} 
        \sqrt{8 \sqrt{k(m-k)}} \sqrt{ \frac{\delta}{\nu}} \leq \sqrt{ 8n \sqrt{ \frac{10n^4}{\theta}} \delta} = \sqrt{ \sqrt{ \frac{10}{\theta}} 8n^3 \delta} .
        \label{eqn: exact step five}
    \end{equation}
    \item[\textbf{6.}]\textbf{Union bound on events that guarantee success for one step.} \\
    Let $E_a, E_b, E_c,$ and $E_d$ be events that correspond to the results (a), (b), (c), and (d) from step 5 and let $E_k$ be the event that a sufficient dividing line is found and $k$ computed accurately. We have just shown 
    \begin{equation} 
        {\mathbb P}(E_a \cap E_b \cap E_c \cap E_d \; | \; E_k) \geq 1 - \frac{2 \theta}{5n^4}.
        \label{eqn: exacts step six conditional prob}
    \end{equation}
    Since we also know ${\mathbb P}(E_k) \geq 1 - \frac{\theta}{10n^4}$, we conclude
    \begin{equation} 
        \aligned 
        {\mathbb P}(E_a \cap E_b \cap E_c \cap E_d \cap E_k) &= {\mathbb P}(E_a \cap E_b \cap E_c \cap E_d \; | \; E_k) {\mathbb P}(E_k) \\
        & \geq \left( 1 - \frac{2 \theta}{5n^4} \right) \left( 1 - \frac{\theta}{10n^4} \right) \\
        & \geq 1 - \frac{\theta}{2n^4} .
        \endaligned
        \label{eqn: exact step six Bayes}
    \end{equation}
    In the remainder of the proof, we will show that conditioning on these events guarantees success for one step of \textbf{EIG}.
    \item[\textbf{7.}] \textbf{If $E_a, E_b, E_c$, and $E_d$ hold, then $\Lambda_{4 \epsilon / 5}(A_{11}, B_{11})$ and $\Lambda_{4 \epsilon / 5} (A_{22}, B_{22})$ are both shattered with respect to $g$.} \\
    Let $\widehat{U}_R^{(k)}$, $\widehat{U}_L^{(k)}$, $\widehat{U}_R^{(m-k)}$, and $\widehat{U}_L^{(m-k)}$ be the matrices such that $||U_R^{(k)} - \widehat{U}_R^{(k)}||_2$, $||U_L^{(k)} - \widehat{U}_L^{(k)}||_2$, $||U_R^{(m-k)} - \widehat{U}_R^{(m-k)}||_2$, and $||U_L^{(m-k)} - \widehat{U}_L^{(m-k)}||_2$ are all bounded above by $\sqrt{ \sqrt{ \frac{10}{\theta}} 8n^3 \delta}$ as in step five. Since $\delta \leq \sqrt{ \frac{\theta}{10}} \frac{\epsilon^2}{7200 n^{2 \alpha + 3}}$, we can replace this upper bound with $\frac{\epsilon}{30n^{\alpha}}$. With this in mind, let
    \begin{equation} 
        \aligned 
        (\widehat{A}_{11}, \widehat{B}_{11}) &= \left( (\widehat{U}_L^{(k)})^H A \widehat{U}_R^{(k)}, \;  (\widehat{U}_L^{(k)})^H B \widehat{U}_R^{(k)} \right) \\
        (\widehat{A}_{22}, \widehat{B}_{22}) &= \left( (\widehat{U}_L^{(m-k)})^H A \widehat{U}_R^{(m-k)}, \; (\widehat{U}_L^{(m-k)})^HB \widehat{U}_R^{(m-k)} \right).
        \endaligned 
        \label{eqn: exact step seven subproblems}
    \end{equation}
    Note that by \cref{lem:exact projectors}, $\Lambda_{\epsilon}(\widehat{A}_{11}, \widehat{B}_{11})$ and $\Lambda_{\epsilon}(\widehat{A}_{22}, \widehat{B}_{22})$ are both contained in $\Lambda_{\epsilon}(A,B)$ and therefore shattered with respect to $g$. At the same time,
    \begin{equation} 
        \aligned 
        ||A_{11} - \widehat{A}_{11}||_2 &= || (U_L^{(k)})^HAU_R^{(k)} - (\widehat{U}_L^{(k)})^HA \widehat{U}_R^{(k)} ||_2 \\
        & = || (U_L^{(k)})^HAU_R^{(k)} - (\widehat{U}_L^k)^HAU_R^{(k)} + (\widehat{U}_L^{(k)})^H A U_R^{(k)} - (\widehat{U}_L^{(k)})^HA \widehat{U}_R^{(k)} ||_2 \\
        & \leq ||(U_L^{(k)} - \widehat{U}_L^{(k)})^HAU_R^{(k)}||_2 + ||(\widehat{U}_L^{(k)})^HA(U_R^{(k)} - \widehat{U}_R^{(k)}||_2 \\
        & \leq || U_L^{(k)} - \widehat{U}_L^{(k)}||_2 ||A||_2 + ||A||_2 ||U_R^{(k)} - \widehat{U}_R^{(k)} ||_2 \\
        & \leq 6n^{\alpha} \frac{\epsilon}{30n^{\alpha}} \\
        & = \frac{\epsilon}{5} 
        \endaligned
        \label{exact step seven error}
    \end{equation}
    and similarly $||B_{11} - \widehat{B}_{11}||_1 \leq \frac{\epsilon}{5} $. Thus, by \cref{lem: shattering preserved}, $\Lambda_{4 \epsilon / 5}(A_{11}, B_{11})$ is shattered with respect to $g$. Repeating this argument for $(A_{22}, B_{22})$ and $(\widehat{A}_{22}, \widehat{B}_{22})$, we conclude $\Lambda_{4 \epsilon / 5}(A_{22}, B_{22})$ is also shattered with respect to $g$. 
\end{enumerate}

\indent In the preceding analysis, we showed that for one step of recursion in \textbf{EIG}, a sufficient dividing line is found, $k$ is computed correctly, and $\Lambda_{4 \epsilon / 5}(A_{11}, B_{11})$ and $\Lambda_{4 \epsilon / 5}(A_{22}, B_{22})$ are both shattered with respect to $g$, and therefore also with respect to the half grids $g_R$ and $g_L$, with probability at least $1 - \frac{\theta}{2n^4}$. Since multiplying by matrices with orthonormal columns will preserve the norm requirements, we conclude that the subsequent calls to \textbf{EIG} in lines 21 and 23 are valid when these events occur. Hence, each recursive step succeeds with probability at least $1 - \frac{\theta}{2n^4}$. Since the recursive tree of \textbf{EIG} has depth at most $\log_{5/4}(n)$ and each step calls \textbf{EIG} twice, a union bound implies that the first guarantee of \textbf{EIG} fails with probability at most 
\begin{equation}
    2 \cdot 2^{ \log_{5 / 4}(n)} \frac{\theta}{2n^4} \leq 2n^4 \frac{\theta}{2n^4} = \theta .
    \label{eqn: exact failure prob}
\end{equation}

\indent We turn now to the second guarantee when $\sigma_n(B) \geq 1$. In this case, we will show that conditioning on the same events that ensure the first guarantee also imply the second. Since \textbf{EIG} builds the approximate eigenvectors recursively, we do this inductively. \\
\indent The base case here corresponds to $m = 1$, in which case \textbf{EIG} gets the one right unit eigenvector ($v=1$) exactly correct. Suppose now we are reconstructing the eigenvectors of $(A',B')$ from the two sub-problems $(A_{11},B_{11})$ and $(A_{22},B_{22})$ it is split into. $(A',B')$ here is any pencil obtained in the divide-and-conquer process.\\
\indent Let $\widehat{T}$ and $\widetilde{T}$ be the invertible matrices obtained from applying \textbf{EIG} to $(A_{11},B_{11})$ and $(A_{22},B_{22})$ as in lines 21 and 23 of the algorithm. Since these calls to \textbf{EIG} pass the parameter $\frac{\beta}{3} $, we can assume each column of $\widehat{T}$ or $\widetilde{T}$ is at most $\frac{\beta}{3} $ away from a true unit right eigenvector of $(A_{11},B_{11})$ or $(A_{22},B_{22})$. In addition, let 
\begin{equation} 
    T = \begin{pmatrix} U_R^{(k)} & U_R^{(m-k)} \end{pmatrix} \begin{pmatrix} \widehat{T} & 0 \\ 0 & \widetilde{T} \end{pmatrix}
    \label{eqn: eigenvector build}
\end{equation}
be the matrix of approximate eigenvectors of $(A',B')$ computed in line 24. Finally let $\widehat{U}_R^{(k)}$ and $\widehat{U}_R^{(m-k)}$ be the true matrices approximated by $U_R^{(k)}$ and $U_R^{(m-k)}$, as in step 5 above. \\
\indent Consider now a column $t_i$ of $T$. It suffices to handle the case where $t_i = U_R^{(k)} \widehat{t}_i$ for a column $\widehat{t}_i$ of $\widehat{T}$, as the same argument applies exactly if $t_i = U_R^{(m-k)} \widetilde{t}_i$ for $\widetilde{t}_i$ a column of $\widetilde{T}$. By our induction hypothesis, we know there exists a true right unit eigenvector $\widehat{v}_i$ of $(A_{11},B_{11})$ such that 
\begin{equation}
|| \widehat{t}_i - \widehat{v}_i ||_2 \leq \frac{\beta}{3} .
\label{eqn: eig 1}
\end{equation}
Now let $(\widehat{A}_{11}, \widehat{B}_{11})$ be the true problem approximated by $(A_{11}, B_{11})$. Conditioning on the same events used above, we know (following the same arguments as in steps 5 and  7)
\begin{equation}
    ||A_{11} - \widehat{A}_{11} ||_2, \; ||B_{11} - \widehat{B}_{11}||_2 \leq 6 n^{\alpha} ||U_R^{(k)} - \widehat{U}_R^{(k)}||_2  \leq 6 n^{\alpha} \sqrt{ \sqrt{ \frac{10}{\theta}} 8n^3 \delta},
    \label{eqn: nearby subproblem}
\end{equation}
which, applying the bound $\delta \leq \sqrt{ \frac{\theta}{10}} \frac{\eta^2}{288n^{2 \alpha + 3}}$, becomes
\begin{equation} 
    ||A_{11} - \widehat{A}_{11}||_2 , \; ||B_{11} - \widehat{B}_{11} ||_2 \leq 6 n^{\alpha} \sqrt{ \sqrt{ \frac{10}{\theta} } 8n^3 \sqrt{ \frac{\theta}{10}} \frac{\eta^2}{288n^{2 \alpha + 3}}} = \eta .
    \label{eqn: nearby subproblem bound}
\end{equation}
Thus, by \cref{lem: shattered eigenvectors}, there exists a right unit eigenvector $\bar{v}_i$ of $(\widehat{A}_{11}, \widehat{B}_{11})$ such that 
\begin{equation}
    ||\widehat{v}_i - \bar{v}_i||_2 \leq \frac{\sqrt{8} \omega}{\pi}\frac{\eta}{\epsilon (\epsilon - \eta)} (1 + ||\widehat{B}_{11}^{-1} \widehat{A}_{11}||_2) ||\widehat{B}_{11}||_2 .
    \label{eqn: apply eignevector perturbation bound}
\end{equation}
\indent To simplify this bound, we first observe that (again by the argument made in step 7 above) we can assume $\epsilon - \eta \geq \frac{4}{5} \epsilon$. In addition, $|| \widehat{A}_{11}||_2 \leq 3$ and $|| \widehat{B}_{11}||_2 \leq 3n^{\alpha}$. Finally, since $\sigma_n(B) \geq 1$ and any split of divide-and-conquer can decrease the smallest singular value by at most $\eta$ (by the stability of singular values), and since the decision tree of \textbf{EIG} has depth at most $\log_{5/4}(n)$, the bound $\eta \leq \frac{1}{2 \log_{5/4}(n)} $ ensures
\begin{equation}
    \sigma_n(\widehat{B}_{11}) \geq 1 - \log_{5/4}(n) \eta \geq \frac{1}{2} . 
    \label{stability of smallest singular value}
\end{equation}
Putting everything together, we have
\begin{equation}
||\widehat{v}_i - \bar{v}_i||_2 \leq \frac{\sqrt{8} \omega}{\pi}  \frac{\eta}{\frac{4}{5} \epsilon^2}  \left(1 + \frac{||\widehat{A}_{11}||_2}{\sigma_n(\widehat{B}_{11})} \right) ||\widehat{B}_{11}||_2  \leq \frac{105 \sqrt{8}}{4 \pi}  \frac{\omega n^{\alpha}}{\epsilon^2} \eta \leq \frac{\beta}{3}.
\label{eqn: eig 2}
\end{equation}
since $ \eta \leq \frac{4 \pi}{315 \sqrt{8}} \frac{\beta \epsilon^2}{\omega n^{\alpha}}$. Now let $v_i = \widehat{U}_R^{(k)}\bar{v}_i$, which is  a true right unit eigenvector of $(A',B')$. By construction, we have
\begin{equation} 
    \aligned 
    ||t_i - v_i||_2 &= ||U_R^{(k)}\widehat{t}_i - \widehat{U}_R^{(k)} \bar{v}_i ||_2 \\
    & =  ||U_R^{(k)} \widehat{t}_i - U_R^{(k)} \widehat{v}_i + U_R^{(k)} \widehat{v}_i - \widehat{U}_R^{(k)} \widehat{v}_i + \widehat{U}_R^{(k)} \widehat{v}_i - \widehat{U}_R^{(k)} \bar{v}_i||_2 \\
    & \leq || U_R^{(k)} (\widehat{t}_i - \widehat{v}_i)||_2  + || (U_R^{(k)} - \widehat{U}_R^{(k)})\widehat{v}_i||_2 + || \widehat{U}_R^{(k)} (\widehat{v}_i - \bar{v}_i)||_2. \\
    & \leq ||\widehat{t}_i - \widehat{v}_i||_2 + ||U_R^{(k)} - \widehat{U}_R^{(k)}||_2 + ||\widehat{v}_i - \bar{v}_i||_2 .
    \endaligned 
    \label{eqn: nearby true eigenvector} 
\end{equation}
Applying  \eqref{eqn: eig 1} and \eqref{eqn: eig 2} to this and using the fact that $|| U_R^{(k)} - \widehat{U}_R^{(k)}||_2 \leq \frac{\beta}{3}$, we conclude
\begin{equation} 
    ||t_i - v_i||_2 \leq \frac{\beta}{3} + \frac{\beta}{3} + \frac{\beta}{3} = \beta .
    \label{eqn: nearby eigenvector triangle}
\end{equation}
By induction, we obtain the same bound for the approximate and true eigenvectors of $(A,B)$.
\end{proof}

\indent The condition on the eigenvector guarantee in \cref{thm: EIG succeeds} (i.e., that $\sigma_n(B) \geq 1$) may seem restrictive, but it reflects the use case in the following section. That is, while it is possible to adjust the parameters of \textbf{EIG} to allow for less strict lower bounds on $\sigma_n(B)$, we plan to apply $\textbf{EIG}$ to the perturbed and scaled $(\widetilde{A}, n^{\alpha} \widetilde{B})$, where by construction $\sigma_n(n^{\alpha} \widetilde{B}) \geq 1$ with high probability.

\section{Randomized Pencil Diagonalization}
\label{section diag}
We can now prove our main result, presented in the introduction as \cref{thm: main result}. First, we state our algorithm for producing an approximate diagonalization and show it succeeds in exact arithmetic with high probability. We follow that with corresponding forward error guarantees and an evaluation of the complexity of the routine. 
\subsection{Diagonalization Routine}
In the previous section, we constructed an eigensolver that produces approximations of the eigenvalues and eigenvectors of $(\widetilde{A}, n^{\alpha} \widetilde{B})$, contained in the outputs $(D_1, D_2)$ and $T$ respectively. Since the $n^{\alpha}$ scaling only changes eigenvalues and can  be easily undone, this effectively yields eigenvalues/eigenvectors for $(\widetilde{A}, \widetilde{B})$. That is, $T$ contains approximate eigenvectors of $(\widetilde{A}, \widetilde{B})$ while the corresponding eigenvalues belong to $(n^{\alpha} D_1, D_2)$ or equivalently $(D_1, n^{- \alpha} D_2)$. This prompts the question:\ how accurate are these approximations? \\
\indent Unfortunately, as noted earlier, little can be said. While perturbation results do exist for the simple eigenvalues of a regular pencil or for the eigenvalues of a regular pencil obtained by perturbing a diagonalizable one \cite[Theorems VI.2.2 and VI.2.6]{stewart1990matrix}, these are typically stated in terms of a chordal metric on the Riemann sphere. Our version of Bauer-Fike (\cref{thm: pencil BF}) provides a route to some Euclidean norm bounds, which we state below, but its dependence on the the smallest singular value of $B$ means it is only useful if $B$ is far from singular. This reflects a reality of our algorithm:\ \textbf{EIG} always produces a set of finite eigenvalues, and we cannot identify which ones are close to true finite eigenvalues of $(A,B)$ and which ones are meant to be infinite. \\
\indent Instead, we leverage \textbf{EIG} to produce an approximate diagonalization of $(A,B)$. Recall from \cref{section: background} that if $B$ is invertible and $T$ contains right eigenvectors of $(A,B)$ then $S = BT$ produces a diagonalization of the pencil such that $S^{-1}AT = D_1$ and $S^{-1}BT = I$.  Using again the fact that $\widetilde{B}$ is invertible almost surely, we conclude that if $T$ is the output of \textbf{EIG} then it and $\widetilde{B}T$ should approximately diagonalize $(\widetilde{A}, \widetilde{B})$. Assuming $||A - \widetilde{A}||_2$ and $||B - \widetilde{B}||_2$ are small, this produces an approximate diagonalization of $(A,B)$. We state \textbf{RPD}, a routine that wraps around \textbf{EIG} to produce this diagonalization, as \cref{alg: RPD} below.
\label{section RPD}
\begin{algorithm}
\caption{Randomized Pencil Diagonalization (\textbf{RPD})\\
\textbf{Input:}  $A,B \in {\mathbb C}^{n \times n}$ and $\varepsilon < 1$ a desired accuracy. \\
\textbf{Requires:} $||A||_2, ||B||_2 \leq 1$ \\
\textbf{Output:} Nonsingular $S, T$ and diagonal $D$ such that $||A - SDT^{-1}||_2, ||B - ST^{-1}||_2 \leq \varepsilon$ with high probability. } \label{alg: RPD}
\begin{algorithmic}[1]
\State $\gamma = \frac{\varepsilon}{16}$
\State $\alpha = \frac{ \left\lceil 2\log_n(1/\gamma) + 3 \right\rceil}{2}$
\State $\epsilon = \gamma^5 / (64n^{\frac{11 \alpha + 25}{3}} + \gamma^5)$
\State $\beta =  \frac{\varepsilon \gamma^2}{24 (1 + 4 \gamma)} n^{- 3 \alpha - 5} $ 
\State $\omega = \frac{\gamma^4}{4}n^{-\frac{8 \alpha + 13}{3}}$ 
\State Draw two independent Ginibre matrices $G_1, G_2 \in {\mathbb C}^{n \times n}$
\State $(\widetilde{A}, \widetilde{B}) = (A + \gamma G_1, B + \gamma G_2)$
\State Draw $z$ uniformly from the box of side length $\omega$ cornered at $-4 - 4i$
\State $g = \text{grid}(z, \omega, \lceil 8 / \omega \rceil, \lceil 8 / \omega \rceil)$
\State $[T, D_1, D_2] = \textbf{EIG}(n, \widetilde{A}, n^{\alpha} \widetilde{B}, \epsilon, \alpha, g, \beta, 1/n)$ 
\For{$i = 1:n$}
    \State $D(i,i) = n^{\alpha}\frac{D_1(i,i)}{D_2(i,i)}$
\EndFor
\State $S = \widetilde{B}T$ 
\State \Return $S, T, D$
\end{algorithmic}
\end{algorithm}

\begin{figure}
    \centering
    \resizebox{.95\linewidth}{!}{
   \begin{tikzpicture}[node distance=1.5cm, auto]
     \node (start) [startstop] {\textbf{Start:} $A,B \in {\mathbb C}^{n \times n}$ with $||A||_2, ||B||_2 \leq 1$.};
     \node (pro1) [process, below of=start] { Build the random grid $g$ and perturb:
      $(\widetilde{A}, \widetilde{B}) = (A + \gamma G_1, B + \gamma G_2)$.};
      \node (input1) [io, below of=pro1, yshift=.2cm]{\textbf{Inputs:} $m = n$; $(\widetilde{A}, n^{\alpha} \widetilde{B})$; $g$.};
      \node (eig) [point, below of=input1, yshift=.2cm] {\textbf{EIG}};
      \node (dec2) [decision, below of=eig, yshift=-.4cm] {$m = 1$?};
      \node (pro2) [process, left of=dec2, yshift= -2cm, xshift = -3.5cm] {Pick a grid line $\text{Re}(z) = h$ and apply a M\"{o}bius transformation that sends the line to the unit circle. };
      \node (stop) [startstop2, right of=pro2, xshift=6.5cm] {\textbf{End EIG:} Recursively build $T$, $D_1$, and $D_2$.};
      \node (irs) [process, below of=pro2] {Apply \textbf{IRS} and \textbf{GRURV} to the transformed pencil and read off the rank $k$.};
      \node (dec1) [decision, below of=irs, yshift=-.7cm] {Good split?};
      \node (def1) [process, below of=dec1, yshift=-.7cm] {Apply \textbf{DEFLATE} to find right and left spaces corresponding to eigenvalues with $\text{Re}(z) > h$.};
      \node (pro5) [process, below of=def1, yshift=-.2cm] {Apply a second M\"{o}bius transformation and call \textbf{DEFLATE} again to find right and left spaces corresponding to eigenvalues with $\text{Re}(z) < h$.};
      \node (pro6) [process, below of=pro5, yshift=-.2cm] {Split the problem into two smaller ones: $(A_{11}, B_{11})$ and $(A_{22}, B_{22})$};
      \node (inputs2) [io, below of=pro6] {\textbf{Inputs:}
      $k$; $(A_{11}, B_{11})$; $g_R$. \; \;  $m-k$; $(A_{22}, B_{22})$; $g_L$.};
      \node (stop2) [startstop, below of=stop]{\textbf{Stop:} Output $T$, $S = \widetilde{B}T$, and $D = n^{\alpha} D_2^{-1}D_2$.};
     \draw [arrow] (start) -- (pro1);
     \draw [arrow] (pro1) -- (input1);
     \draw [arrow] (input1) -- (eig);
     \draw [arrow] (eig) -- (dec2);
     \draw [arrow] (pro2) -- (irs);
     \draw [arrow] (irs) -- (dec1);
     \draw [arrow] (dec1) -- node[anchor=east] {yes} (def1);
     \draw [arrow] (def1) -- (pro5);
     \draw [arrow] (pro5) -- (pro6);
     \path[line, rounded corners] (dec2.east) -- node[anchor=south]{yes} ++ (1.65cm,0) -|  (stop);
     \path[line, rounded corners] (dec2.west) -- node[anchor=south]{no} ++ (-1.65cm,0) -| (pro2);
     \path [line, rounded corners] (dec1.east) -- node[anchor=south] {no} ++ (3.8cm,0) |- (pro2.east); 
     \draw [arrow] (pro6) -- (inputs2);
     \path [line, rounded corners] (inputs2.west) -- ++ (-2.5cm,0) |- (eig.west);
     \draw [arrow] (stop) -- (stop2);
   \end{tikzpicture}
   }
   \vspace{.5cm}
    \caption{A diagram of \textbf{RPD} (\cref{alg: RPD}) for producing an approximate diagonalization of $(A,B)$. As in the pseudocode, we assume for simplicity that each split in \textbf{EIG} is made by a vertical grid line.}
    \label{fig: RPD flow}
\end{figure}
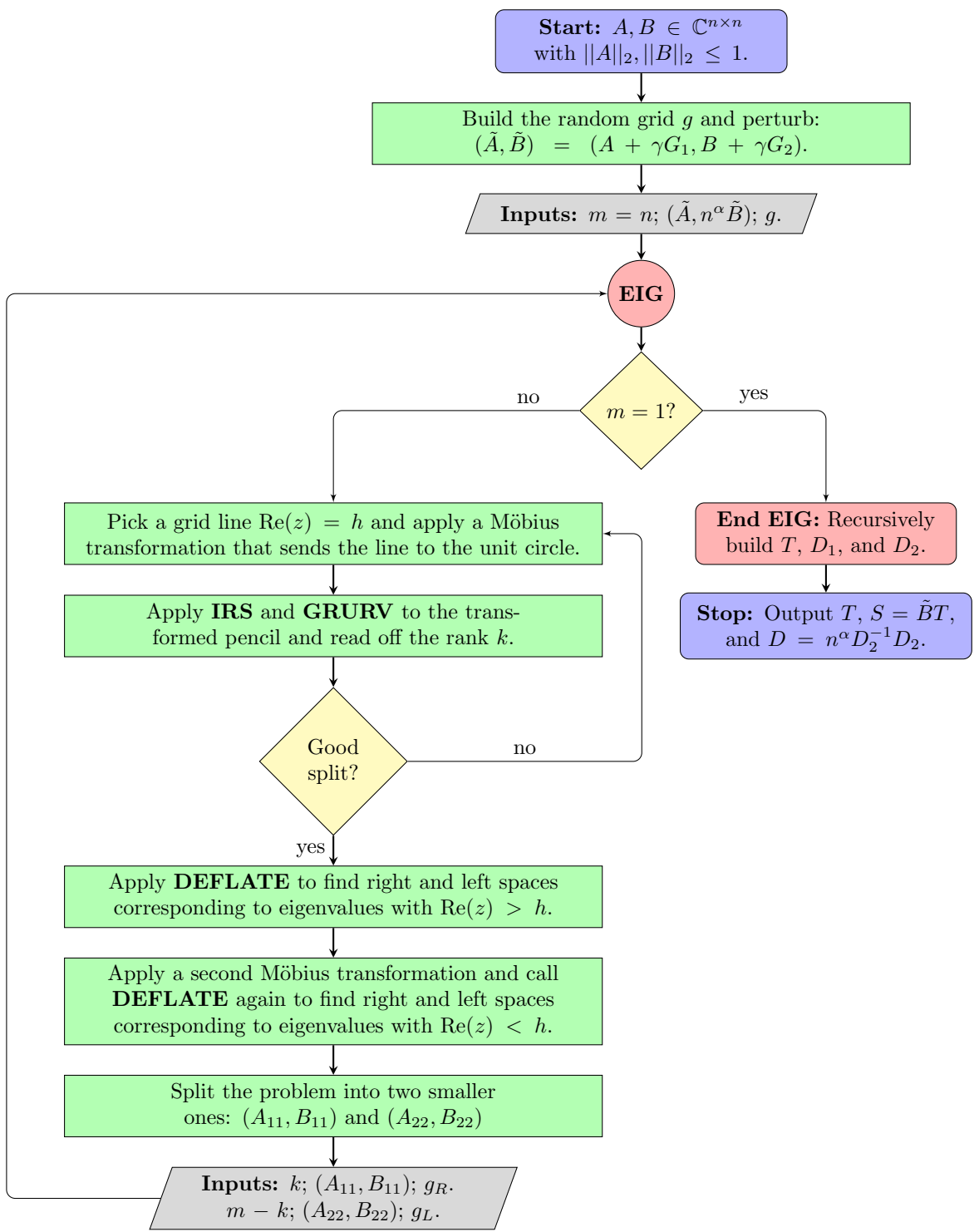

Note that the assumption $||A||_2, ||B||_2 \leq 1$ is essentially made for convenience; we can obtain a diagonalization of any pencil $(A,B)$ via \textbf{RPD} by first normalizing the matrices accordingly. In \cref{thm: RPD guarantees}, we show that \textbf{RPD} produces an approximate diagonalization with the given accuracy in exact arithmetic, thereby proving the bulk of our main result, \cref{thm: main result}. For reference, \cref{fig: RPD flow} provides a high level overview of \textbf{RPD}, including the details of its call to \textbf{EIG}.

\begin{thm}\label{thm: RPD guarantees}
  For any $A,B \in {\mathbb C}^{n \times n}$ with $||A||_2, ||B||_2 \leq 1$, the outputs of exact-arithmetic \text{\normalfont \bf RPD} satisfy
  $$ ||A - SDT^{-1}||_2, \; ||B - ST^{-1}||_2 \leq \varepsilon $$
  with probability at least
  $$ \left[1 - \frac{82}{n} - \frac{531441}{16n^2} \right] \left[ 1 - \frac{1}{n} - 4e^{-n} \right] \left[ 1 - \frac{1}{n} \right] .$$
\end{thm}
\begin{proof}
Consider the perturbed and scaled pencil $(\widetilde{A}, n^{\alpha} \widetilde{B})$. By \cref{thm: shattering} and its proof, with probability at least $\left( 1 - \frac{82}{n} - \frac{531441}{16n^2} \right) \left(1 - \frac{n^{2-2\alpha}}{\gamma^2} - 4e^{-n} \right) $ we have the following:
\begin{enumerate}
    \item $||G_1||_2, ||G_2||_2 \leq 4$,
    \item $||\widetilde{A}||_2 \leq 3$,
    \item $||n^{\alpha} \widetilde{B}||_2 \leq 3n^{\alpha}$,
    \item $ \sigma_n(n^{\alpha} \widetilde{B}) \geq 1 $,
    \item $\Lambda(\widetilde{A}, n^{\alpha} \widetilde{B}) \subseteq B_3(0)$,
    \item $\kappa_V(n^{-\alpha} \widetilde{B}^{-1} \widetilde{A}) \leq \frac{n^{\alpha+2}}{\gamma}$,
    \item $\Lambda_{\epsilon}(\widetilde{A}, n^{\alpha} \widetilde{B})$ is shattered with respect to the grid $g$ (for $\epsilon$ as in line 3). 
\end{enumerate}
Conditioning on these events, we observe that 2, 3, and 7 ensure that the call to \textbf{EIG} in line 10 is valid, meaning with probability at least $1 - \frac{1}{n}$ we can add the guarantees of \cref{thm: EIG succeeds} to our list:
\begin{itemize}
    \item[8.] Each eigenvalue of $(D_1,D_2)$ shares a grid box of $g$ with a true eigenvalue of $(\widetilde{A}, n^{\alpha} \widetilde{B})$.
    \item[9.] Since $\sigma_n(n^{\alpha} \widetilde{B}) \geq 1$ (item 4) each column $t_i$ of $T$ satisfies $||t_i - v_i||_2 \leq \beta$ for $v_i$ a true unit right eigenvector of $(\widetilde{A}, n^{\alpha} \widetilde{B})$.
\end{itemize}
By the definition of conditional probability, all nine of these events occur simultaneously with probability at least 
\begin{equation}
    \left[ 1 - \frac{82}{n} - \frac{531441}{16n^2} \right] \left[ 1 - \frac{n^{2-2\alpha}}{\gamma^2} - 4e^{-n} \right] \left[ 1 - \frac{1}{n} \right] .
    \label{eqn: conditional prob for RPD}
\end{equation}
Since the choice of $\alpha$ in line 2 guarantees
\begin{equation}
    \frac{n^{2 - 2\alpha}}{\gamma^2} \leq \frac{1}{n} 
    \label{eqn: alpha guarantee}
\end{equation}
\eqref{eqn: conditional prob for RPD} can be bounded from below by 
\begin{equation}
    \left[ 1 - \frac{82}{n} - \frac{531441}{16n^2} \right] \left[ 1 - \frac{1}{n} - 4e^{-n} \right] \left[ 1 - \frac{1}{n} \right] .
    \label{eqn: n dependent conditional prob for RPD}
\end{equation}
To complete the proof we show that the nine items listed above guarantee both $||A - SDT^{-1}||_2 \leq \varepsilon$ and $||B - ST^{-1}||_2 \leq \varepsilon$. The second of these is trivial:\ since $\gamma = \frac{\varepsilon}{16}$ and $ST^{-1} = \widetilde{B}$ we have
\begin{equation}
    ||B - ST^{-1}||_2 = ||B - \widetilde{B}||_2 = ||\gamma G_2||_2 \leq 4 \gamma = \frac{\varepsilon}{4} < \varepsilon .
    \label{eqn: RPD guarantee for B}    
\end{equation}
To show the same for $A$ and $SDT^{-1}$ we use the following key fact. Let $V$ be the matrix whose columns contain, in order, the true right unit eigenvectors of $(\widetilde{A}, n^{\alpha} \widetilde{B})$ guaranteed by item 9. Since $\widetilde{B}$ is invertible with probability one and $(\widetilde{A}, \widetilde{B})$ and $(\widetilde{A}, n^{\alpha} \widetilde{B})$ have the same set of eigenvectors, $V$ diagonalizes $\widetilde{B}^{-1} \widetilde{A}$. Thus, there exists a diagonal matrix $\Lambda$ such that $\widetilde{B}^{-1} \widetilde{A} = V \Lambda V^{-1}$. Moreover, since the eigenvalue of $(\widetilde{A}, n^{\alpha} \widetilde{B})$ corresponding to $v_i$ shares a grid box of $g$ with the eigenvalue of $(D_1, D_2)$ corresponding to $t_i$ (this was how $v_i$ was found in the proof of \cref{thm: EIG succeeds}), each diagonal entry of $\Lambda$ is at most $\sqrt{2} n^{\alpha} \omega$ away from the corresponding diagonal entry of $D = n^{\alpha} D_2^{-1}D_1$. Note that the eigenvalues of $(D_1,D_2)$ are all contained in $g$, which guarantees that $D_2$ is invertible. With all of this in mind, expand $||A - SDT^{-1}||_2$ as follows:
\begin{equation}
    \aligned 
    ||A - SDT^{-1}||_2 &= ||A - \widetilde{A} + \widetilde{A} - SDT^{-1}||_2 \\
    & \leq ||A - \widetilde{A}||_2 + ||\widetilde{A} - SDT^{-1}||_2 \\
    &\leq ||A - \widetilde{A}||_2 + ||\widetilde{B}||_2 ||\widetilde{B}^{-1}\widetilde{A} - TDT^{-1}||_2 \\
    & \leq 4 \gamma + (1 + 4 \gamma) ||V \Lambda V^{-1} - TDT^{-1}||_2 \\
    & \leq 4 \gamma + (1 + 4 \gamma) \left[ ||(V - T) \Lambda V^{-1}||_2 + ||T(\Lambda - D)V^{-1}||_2 + ||TD(V^{-1} - T^{-1})||_2 \right].
    \endaligned 
    \label{eqn: RPD guarantee for A 1}
\end{equation}
To simplify this bound in terms of our parameters, we observe the following.
\begin{itemize}
    \item Since the columns of $V$ and $T$ are unit vectors, $1 \leq ||V||_2, ||T||_2 \leq \sqrt{n}$.
    \item Since $V$ diagonalizes $\widetilde{B}^{-1}\widetilde{A}$ and therefore also $n^{-\alpha} \widetilde{B}^{-1} \widetilde{A}$ and $\kappa_V(n^{-\alpha} \widetilde{B}^{-1} \widetilde{A}) \leq \frac{n^{\alpha+2}}{\gamma}$, \cref{rem: condition number} allows us to assume  
    \begin{equation}
        ||V||_2 ||V^{-1}||_2 \leq \frac{n^{\alpha + 2}}{\gamma}
        \label{eqn: V inverse norm guarantee}
    \end{equation}
    without changing probabilities. Combined with the previous point, this implies $||V^{-1}||_2 \leq \frac{n^{\alpha+2}}{\gamma}$.
    \item The columns of $T$ and $V$ satisfy $||t_i - v_i||_2 \leq \beta$ so $||T - V||_2 \leq \sqrt{n} \beta$
    \item By the stability of singular values, $||T - V||_2 \leq \sqrt{n} \beta $ implies $\sigma_n(T) \geq \sigma_n(V) - \sqrt{n} \beta $. Combining this with our upper bound on $||V^{-1}||_2$ yields
    \begin{equation}
        ||T^{-1}||_2 \leq \frac{n^{\alpha+2}}{\gamma - n^{\frac{2 \alpha+5}{2}}\beta} .
        \label{eqn: T inverse norm guarantee}
    \end{equation}
    \item Because each diagonal entry of $\Lambda$ is at most $\sqrt{2}n^{\alpha} \omega$ from the corresponding diagonal entry of $D$, $||\Lambda - D||_2 \leq \sqrt{2}n^{\alpha} \omega$. 
    \item Since both $(\widetilde{A}, n^{\alpha} \widetilde{B})$ and $(D_1, D_2)$ have eigenvalues in $B_3(0)$, $||\Lambda||_2, ||D||_2 \leq 3n^{\alpha} $.
    \item Finally, $||V^{-1} - T^{-1}||_2 = ||T^{-1}(T-V)V^{-1}||_2 \leq ||T^{-1}||_2 ||T-V||_2 ||V^{-1}||_2$.
\end{itemize}
Together, these imply the following bounds:
\begin{subequations}
    \begin{align}
    & ||(V-T) \Lambda V^{-1}||_2 \leq ||V-T||_2 ||\Lambda||_2 ||V^{-1}||_2 \leq \sqrt{n} \beta \cdot 3n^{\alpha} \cdot \frac{n^{\alpha+2}}{\gamma} = \frac{3 \beta}{\gamma} n^{\frac{4 \alpha + 5}{2}} \label{subeq: piece 1} \\
    & ||T(\Lambda - D)V^{-1}||_2 \leq ||T||_2 ||\Lambda - D||_2 ||V^{-1}||_2 \leq \sqrt{n} \cdot \sqrt{2}n^{\alpha} \omega \cdot \frac{n^{\alpha+2}}{\gamma} = \frac{\sqrt{2} \omega}{\gamma} n^{\frac{4 \alpha + 5}{2}} \label{subeq: piece 2} \\
    & ||TD(V^{-1} - T^{-1})||_2 \leq \sqrt{n} \cdot 3n^{\alpha} \cdot \frac{n^{\alpha+2}}{\gamma - n^{\frac{2 \alpha + 5}{2}} \beta} \cdot \sqrt{n} \beta \cdot \frac{n^{\alpha+2}}{\gamma} = \frac{ 3 \beta n^{3 \alpha +5}}{\gamma(\gamma - n^{\frac{2 \alpha+5}{2}} \beta)}. \label{subeq: piece 3} 
    \end{align}
\end{subequations}
Now the choice of $\beta$ in line 4 ensures $\beta \leq \frac{\varepsilon \gamma}{12(1 + 4 \gamma)}n^{-\frac{4 \alpha + 5}{2}}$, so \eqref{subeq: piece 1} simplifies to 
\begin{equation}
    ||(V-T) \Lambda V^{-1}||_2 \leq 3 \cdot \frac{\varepsilon \gamma}{12(1+4 \gamma) n^{\frac{4 \alpha + 5}{2}}} \cdot \frac{n^{\frac{4 \alpha + 5}{2}}}{\gamma} = \frac{\varepsilon}{4(1 + 4 \gamma)}.
    \label{eqn: final piece 1}
\end{equation}
Similarly, using this time the fact that $\beta < \frac{\gamma}{2} n^{-\frac{2 \alpha + 5}{2}}$ and therefore $\gamma - n^{\frac{2 \alpha + 5}{2}}\beta  > \frac{\gamma}{2} $, \eqref{subeq: piece 3} becomes
\begin{equation}
    ||TD(V^{-1} - T^{-1})||_2 \leq \frac{6}{\gamma^2} n^{3 \alpha + 5} \cdot \frac{\varepsilon \gamma^2}{24 (1 + 4 \gamma)} n^{-3 \alpha - 5} = \frac{\varepsilon}{4(1 + 4 \gamma)} .
    \label{eqn: final piece 3}
\end{equation}
Finally, $\omega =  \frac{\gamma^4}{4}n^{-\frac{8 \alpha + 13}{3}}$ implies $\omega \leq \frac{\varepsilon \gamma}{4 \sqrt{2} ( 1 + 4 \gamma)} n^{-\frac{4 \alpha + 5}{2}}$, which allows us to upper bound \eqref{subeq: piece 2} as
\begin{equation}
    ||T (\Lambda - D)V^{-1} ||_2 \leq \sqrt{2} \cdot \frac{\varepsilon \gamma}{4 \sqrt{2}(1+4 \gamma) n^{\frac{4 \alpha + 5}{2}}} \cdot \frac{n^{\frac{4 \alpha + 5}{2}}}{\gamma} = \frac{\varepsilon}{4(1+4 \gamma)}.
    \label{eqn: final piece 2}
\end{equation}
Applying these to \eqref{eqn: RPD guarantee for A 1}, we obtain
\begin{equation}
    ||A - SDT^{-1}||_2 \leq 4 \gamma + (1 + 4 \gamma) \left[ \frac{\varepsilon}{4(1 + 4 \gamma)} + \frac{\varepsilon}{4(1 + 4 \gamma)} + \frac{\varepsilon}{4(1 + 4 \gamma)} \right] = 4 \gamma + \frac{3}{4} \varepsilon .
    \label{eqn: RPD guarantee for A 2}
\end{equation}
Since $\gamma = \frac{\varepsilon}{16}$ we conclude $||A - SDT^{-1}||_2 \leq \varepsilon $.
\end{proof}

\subsection{Forward Error Guarantee and Asymptotic Complexity}

We might hope that \cref{thm: RPD guarantees} can provide a forward error guarantee in terms of the eigenvector condition number $\kappa_V(A,B)$. As we demonstrate below, this is possible when the input pencil is well-behaved (in particular regular with a full set of distinct and finite eigenvalues). The resulting \cref{thm: forward error bound} parallels a similar result of Banks et al.\ \cite[Proposition 1.1]{banks2020pseudospectral}.

\begin{thm} \label{thm: forward error bound}
Let $A,B, \widetilde{A}, \widetilde{B} \in {\mathbb C}^{n \times n}$ with $||A||_2 \leq 1$. Assume that $B$ is invertible and that the pencil $(A,B)$ has a full set of distinct, finite eigenpairs $ \left\{ (\lambda_i, v_i) \right\}_{i = 1}^n$. If $||A - \widetilde{A}||_2 \leq \varepsilon \sigma_n(B) \text{\normalfont gap}(A,B)$ and $||B - \widetilde{B}||_2 \leq \varepsilon \sigma_n^2(B) \text{\normalfont gap}(A,B)$ for
$$ \varepsilon \leq \frac{1}{32 \kappa_V(A,B)} $$
then $(\widetilde{A}, \widetilde{B})$ has no repeated eigenvalues and any set of corresponding eigenpairs $\{ (\widetilde{\lambda}_i, \widetilde{v}_i) \}_{i=1}^n$ satisfies
$$ |\lambda_i - \widetilde{\lambda}_i| \leq 4 \text{\normalfont gap}(A,B) \kappa_V(A,B) \varepsilon \; \; \text{and}  \; \;  ||v_i - \widetilde{v}_i||_2 \leq 24 n \kappa_V(A,B) \varepsilon  \; \; \text{for all} \; \; i = 1, \ldots, n$$
after multiplying the $v_i$ by phases and reordering, if necessary. 
\end{thm}
\begin{proof}
    We first note that $\varepsilon  < \frac{1}{32 \kappa_V(A,B)}$ implies
    \begin{equation}
        ||B - \widetilde{B}||_2 \leq \frac{\sigma_n^2(B) \text{gap}(A,B)}{32 \kappa_V(A,B)} \leq \frac{\sigma_n(B)}{16 \kappa_V(A,B)} \leq \frac{\sigma_n(B)}{2},
        \label{eqn: B minus B tilde}
    \end{equation}
    where the second inequality follows from $\text{gap}(A,B) \leq 2 ||B^{-1}A||_2 \leq \frac{2}{\sigma_n(B)}$ (since $||A||_2 \leq 1$). Thus, \cref{lem: stability of singular values} guarantees $\sigma_n(\widetilde{B}) \geq \sigma_n(B) - ||B - \widetilde{B}||_2 \geq \frac{\sigma_n(B)}{2}$, so $\widetilde{B}$ is invertible. With this in mind, consider the matrices $B^{-1}A$ and $\widetilde{B}^{-1}\widetilde{A}$, which are diagonalizable with the same eigenpairs as $(A,B)$ and $(\widetilde{A}, \widetilde{B})$ respectively. Using again the fact that $||A||_2 \leq 1 $, we have
    \begin{equation}
        \aligned 
        ||B^{-1}A - \widetilde{B}^{-1}\widetilde{A}||_2 &= || B^{-1}A - \widetilde{B}^{-1}A + \widetilde{B}^{-1}A - \widetilde{B}^{-1} \widetilde{A}||_2 \\
        & \leq ||B^{-1} - \widetilde{B}^{-1}||_2 ||A||_2 + ||\widetilde{B}^{-1}||_2 ||A - \widetilde{A}||_2 \\
        & \leq ||\widetilde{B}^{-1}||_2||\widetilde{B} - B||_2 ||B^{-1}||_2 + || \widetilde{B}^{-1}||_2 ||A - \widetilde{A}||_2 \\
        & \leq \frac{2}{\sigma_n(B)} \varepsilon \sigma_n^2(B) \text{gap}(A,B) \frac{1}{\sigma_n(B)} + \frac{2}{\sigma_n(B)} \varepsilon \sigma_n(B) \text{gap}(A,B)\\
        & \leq 4 \varepsilon \text{gap}(A,B).
        \endaligned 
        \label{eqn: formed product difference}
    \end{equation}
    This implies that the eigenvalues of $\widetilde{B}^{-1}\widetilde{A}$ belong to the $4 \varepsilon \text{gap}(A,B)$-pseudospectrum of $B^{-1}A$, and moreover we can continuously deform the eigenvalues of $B^{-1}A$ to those of $\widetilde{B}^{-1}\widetilde{A}$ without leaving this pseudospectrum. Since $\varepsilon < \frac{1}{4}$ and $B^{-1}A$ is diagonalizable, single matrix Bauer-Fike (\cref{thm: bauer fike}) therefore implies that $\widetilde{B}^{-1}\widetilde{A}$, and by extension $(\widetilde{A}, \widetilde{B})$, has a full set of distinct eigenvalues $\lambda_1, \ldots, \lambda_n$. Bauer-Fike also guarantees that these eigenvalues can be ordered so that each $\widetilde{\lambda}_i$ satisfies
    \begin{equation}
        | \lambda_i - \widetilde{\lambda}_i| \leq 4 \text{gap}(A,B) \kappa_V(B^{-1}A) \varepsilon
        \label{eqn: apply BF}
    \end{equation}
    for a corresponding eigenvalue $\lambda_i$ of $B^{-1}A$. Recalling that $\kappa_V(B^{-1}A) = \kappa_V(A,B)$, we obtain the first guarantee. \\
    \indent The second guarantee follows from the proof of \cite[Proposition 1.1]{banks2020pseudospectral} since the requirement $\varepsilon < \frac{1}{32 \kappa_V(A,B)}$ ensures by \eqref{eqn: formed product difference} that
    \begin{equation}
        ||B^{-1}A - \widetilde{B}^{-1}\widetilde{A}||_2 \leq \frac{\text{gap}(B^{-1}A)}{8 \kappa_V(B^{-1}A)} .
        \label{eqn: produce difference with bound}
    \end{equation}
    Note that while \cite[Proposition 1.1]{banks2020pseudospectral} comes with norm assumptions on the matrices, these are not used in the proof of the eigenvector guarantee. 
\end{proof}

\indent Banks et al.\ state \cite[Proposition 1.1]{banks2020pseudospectral} -- the single matrix counterpart to \cref{thm: forward error bound} -- in terms of an eigenvalue problem condition number:
\begin{equation}
    \kappa_{\text{eig}}(A) = \frac{\kappa_V(A)}{\text{gap}(A)}.
    \label{eqn: Banks et al condition number}
\end{equation}
With this in mind, we could rewrite \cref{thm: forward error bound} in terms of the similar quantity
\begin{equation}
    \frac{\kappa_V(A,B)}{\sigma_n^3(B) \text{gap}(A,B)} .
    \label{eqn: possible condition number}
\end{equation}
While \eqref{eqn: possible condition number} does capture some of the challenges of the generalized eigenvalue problem (i.e., it is infinite if $B$ is singular or $(A,B)$ has a repeated eigenvalue) it is not invariant under scaling of $A$ and $B$ and is therefore not suitable as a condition number for the problem. For this reason, we choose instead to state the forward error guarantees in terms of $\kappa_V(A,B)$. \\
\indent The proof of \cite[Proposition 1.1]{banks2020pseudospectral} relies on some basic perturbation results for individual matrices summarized in \cite{Greenbaum}. Developing corresponding results for the generalized problem may allow for a similar forward error guarantee with looser requirements on $||A - \widetilde{A}||_2$ and $||B - \widetilde{B}||_2$. In particular, we might hope to reduce or eliminate the factors of $\sigma_n(B)$. Absent these improvements, \cref{thm: forward error bound} is fairly restrictive, requiring not only that $(A,B)$ is regular and diagonalizable but that $B$ is nonsingular. In part, this is a consequence of our version of Bauer-Fike, which makes many of the same assumptions and produces the same guarantees if used in the proof above. We consider a generalization to arbitrary regular pencils in \cref{appendix: BF}. The singular pencil case is even more complex, requiring a different perturbation analysis altogether (see for example work of Demmel and Kågström \cite{DEMMEL1987139} or Lotz and Noferini \cite{wilkinson_bus}). \\
\indent To wrap up this section, we compute the asymptotic complexity of \textbf{RPD} in terms of $n$, the size of the pencil $(A,B)$, and $\varepsilon$, the accuracy of the approximate diagonalization. Throughout, we assume that we have access to black-box algorithms for multiplying two $n \times n$ matrices and computing the QR factorization of an $n \times n$ matrix, which require $T_{\text{MM}}(n)$ and $T_{\text{QR}}(n)$ arithmetic operations, respectively. For simplicity, we also assume access to  QL/RQ algorithms (used in \textbf{GRURV}) that require $T_{\text{QR}}(n)$ operations. \\
\indent \cref{prop: op count} shows that \textbf{RPD} runs in nearly matrix multiplication time. Its proof is somewhat more subtle than we have implied so far; rather than simply arguing that each step of divide-and-conquer runs in nearly matrix multiplication time and that only logarithmically many steps are needed, we apply a sharper geometric sum that leverages the shrinking problem size guaranteed by significant eigenvalue splits.
\begin{prop} \label{prop: op count}
Exact-arithmetic \text{\normalfont \bf{RPD}} requires at most $ O\left(\log^2 \left( \frac{n}{\varepsilon} \right) T_{\text{\normalfont MM}}(n) \right) $
arithmetic operations. 
\end{prop}
\begin{proof}
    We track only matrix multiplication and QR, as all other building blocks of \textbf{RPD} have smaller complexity. We begin by noting that line 2 of \cref{alg: RPD} implies $\alpha = {\mathcal O} \left( \log_n \left( \frac{1}{\gamma} \right) \right)$ and $\gamma = \Theta(\varepsilon)$, so 
    \begin{equation}
        n^{\alpha} = O\left( n^{ \log_n(1/\gamma)} \right) =  O \left( \frac{1}{\gamma} \right) = O\left( \frac{1}{\varepsilon} \right) .
        \label{eqn: n alpha big O bound}
    \end{equation}
    Meanwhile in line 3 we set
    \begin{equation}
        \epsilon > \frac{\gamma^5}{65n^{\frac{11 \alpha + 25}{3}}} = \Omega \left( \frac{\varepsilon^{26/3}}{n^{25/3}} \right) .
        \label{eqn: epsilon big Omega}
    \end{equation}
    Together, these imply that  the number of steps of repeated squaring required at any point in the recursion can be bounded asymptotically as 
    \begin{equation}
        p =  O \left( \log \left( \frac{n^{\alpha}}{\epsilon} \right) \right) = O \left( \log \left( \frac{n^{25/3}}{\varepsilon^{29/3}} \right) \right) = O \left( \log \left( \frac{n}{\varepsilon} \right) \right) .
        \label{eqn: p big O bound}
    \end{equation}
    \indent Consider now working through one step of divide-and-conquer. Lines 8-15 of \textbf{EIG} make up the bulk of the work, executing a search over the grid lines for one that sufficiently splits the spectrum. For each line that is checked, we make one call to \textbf{IRS} and one call to \textbf{GRURV}; each step of repeated squaring consists of one $2m \times 2m$ QR factorization and two $m \times m$ matrix multiplications, while applying \textbf{GRURV} to a product of two $m \times m$ matrices requires $3T_{\text{QR}}(m) + 2T_{\text{MM}}(m)$ operations. Combining these with \eqref{eqn: p big O bound}, we conclude that each grid line checked results in 
    \begin{equation}
     O\left( \log \left( \frac{n}{\varepsilon} \right) \left[ T_{\text{QR}}(2m) + 2T_{\text{MM}}(m) \right] + 2T_{\text{QR}}(m) + 2T_{\text{MM}}(m) \right) = O \left( \log \left( \frac{n}{\varepsilon} \right) \left[ T_{\text{QR}}(m) + T_{\text{MM}}(m) \right] \right) 
        \label{eqn: line check big O}
    \end{equation}
    operations. Since we check at most $O\left( \log \left( \frac{1}{\omega} \right) \right)$ grid lines each time and 
    \begin{equation}
        \omega = \frac{\gamma^4}{4n^{ \frac{8 \alpha + 13}{3}}} = \Omega \left( \frac{\varepsilon^{20/3}}{n^{13/3}} \right),  
        \label{eqn: omega bit Omega bound}
    \end{equation}
    lines 8-15 of \textbf{EIG} take at most
    \begin{equation}
        O \left( \log^2 \left( \frac{n}{\varepsilon} \right) \left[ T_{\text{QR}}(m) + T_{\text{MM}}(m) \right] \right) = O \left( \log^2 \left( \frac{n}{\varepsilon} \right) T_{\text{MM}}(m) \right)
        \label{eqn: line check big O total}
    \end{equation}
    operations, where we simplify by noting $T_{\text{QR}}(m) = O(T_{\text{MM}}(n))$ \cite[\S 4.1]{2007}. Since the remainder of one step of \textbf{EIG} -- i.e., the subsequent calls to \textbf{DEFLATE} -- has complexity equal to that of checking one grid line, we conclude that \eqref{eqn: line check big O total} is the asymptotic complexity of one step of divide-and-conquer.  \\
    \indent To complete the proof, we sum this expression recursively. Since the $\log^2 \left( \frac{n}{\epsilon} \right)$ term of \eqref{eqn: line check big O total} is independent of $m$, this reduces to summing $T_{\text{MM}}(m)$ over all subproblems produced by \textbf{EIG}. With this in mind, set $T_{\text{MM}}(m) = O(m^{\xi})$ for $\xi \in [2,3]$ and suppose we divide an $m \times m$ pencil into subproblems of size $m_1$ and $m_2$. Since we enforce a significant split, we are guaranteed $\frac{1}{5}m \leq m_1,m_2 \leq \frac{4}{5}m$ and therefore
    \begin{equation}
       m_1^{\xi} + m_2^{\xi} = m_1^{\xi} + (m-m_1)^{\xi} \leq \left(\frac{4}{5} m \right)^{\xi} + \left( \frac{1}{5} m \right)^{\xi} \leq \frac{17}{25} m^{\xi},
        \label{eqn: subproblem matrix multiplication}
    \end{equation}
    where the last inequality is obtained by applying $\xi \geq 2$. Consequently, a sum of $m^{\xi}$ over all subproblems can be bounded by 
    \begin{equation}
        \sum_{k=0}^{\infty} n^{\xi} \left( \frac{17}{25} \right)^k  = n^{\xi} \sum_{k=0}^{\infty} \left(\frac{17}{25}\right)^k = \frac{25}{8} n^{\xi}
    \end{equation}
    and therefore $\sum_m T_{\text{MM}}(m) = O(T_{\text{MM}}(n))$. Applying this to \eqref{eqn: line check big O total} yields the final complexity.
\end{proof}

\section{Numerical Examples}
\label{section experiments}
In this section, we consider several examples to investigate how pseudospectral divide-and-conquer performs in practice. Our first task is to adjust the parameters of \textbf{RPD} and \textbf{EIG}, as the values listed in the pseudocode of \cref{alg:EIG,alg: RPD} -- though necessary in the proof of \cref{thm: RPD guarantees} -- are prohibitively restrictive for implementation. Here, we make the following relaxations.
\begin{itemize}
    \item First, we eliminate the $n^{\alpha}$ scaling, testing examples with eigenvalues exclusively (or predominantly) in $B_3(0)$.
    \item Extracting the main dependence on $\gamma$ and $n$, we set $\epsilon = \beta = \omega = \gamma / n$ and we limit the number of steps of repeated squaring to $p = \lceil \log_2(n/ \epsilon) \rceil$.
    \item Finally, we drop the factor of $1/n^3$ from the criteria used to compute $k$ in line 12 of \textbf{EIG}. 
\end{itemize}

\indent In light of these simplifications, we present the following experiments as simply a proof of concept. Accordingly, we do not consider run times nor do we use an explicitly parallel implementation of the algorithm. Throughout, all results were obtained in Matlab version R2023a.\footnote{Our implementation is available here: \href{https://github.com/ry-schneider/Randomized_Pencil_Diagonalization}{https://github.com/ry-schneider/Randomized\_Pencil\_Diagonalization}}

\subsection{Model Problems}
\indent We start by using \textbf{RPD} as stated (i.e., running to subproblems of size $1 \times 1$) on the following $50 \times 50$ model problems.
\begin{enumerate}
    \item[\textbf{1.}] \textbf{Planted Spectrum:} First, we consider a pencil with equally spaced, real eigenvalues in the interval $[-2,2]$. To obtain $(A,B)$, we fill a diagonal matrix $\Lambda$ as $\Lambda(j,j) = -2 + \frac{4}{49}(j-1)$ and set $A = X\Lambda Y^{-1}$ and $B = XY^{-1}$ for $X$ and $Y$ two independent, complex Gaussian matrices. In accordance with \textbf{RPD}, we then normalize the pencil so that $||A||_2, ||B||_2 \leq 1$. We can think of this example as the best-case scenario, where $\text{gap}(A,B)$ is large and $B$ is far from singular.  
    \item[\textbf{2.}] \textbf{Jordan Block:}  Next, we consider a pencil $(A,B)$ with $B = I$ and $A = J_{50}(0)$ for $J_{50}(0)$ a Jordan block with eigenvalue zero. In contrast to the planted spectrum example, this tests a generalized eigenvalue problem with $\text{gap}(A,B) = 0$.
\end{enumerate}

\begin{figure}[p]
    \centering
    \begin{subfigure}[t]{\linewidth}
        \includegraphics[width=\textwidth]{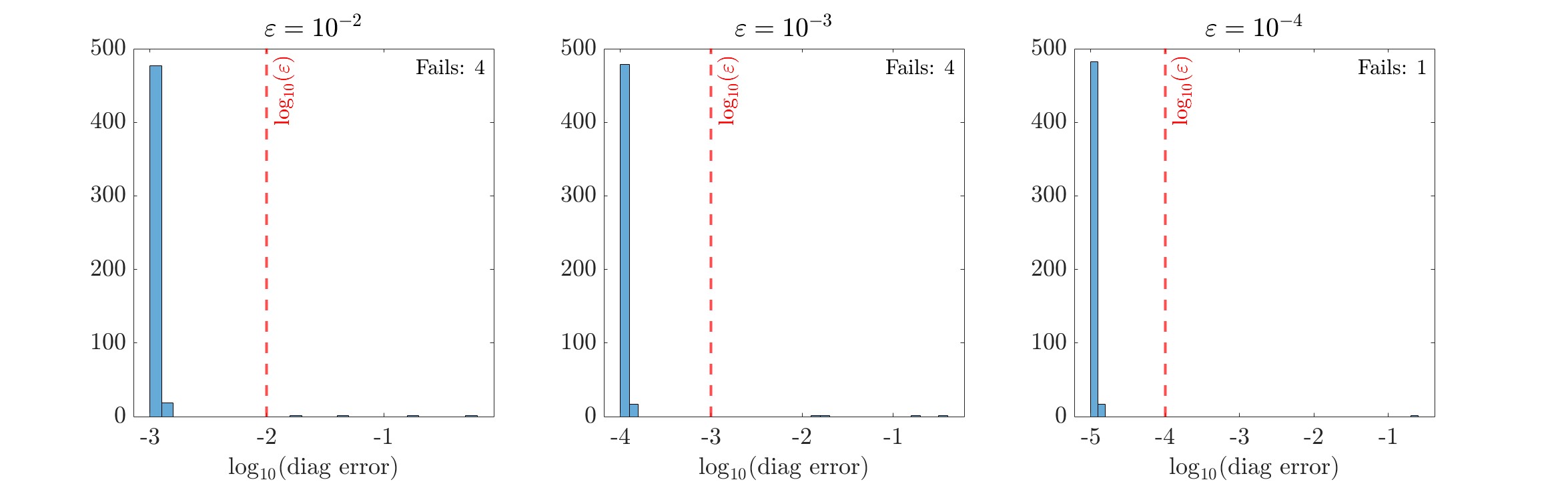}
        \caption{Frequency of diagonalization error \eqref{eqn: diag error}. We mark the $\log_{10}(\varepsilon)$ threshold for success in red and count the number of failed runs. For this example $||A||_2 = 1$ and $||B||_2 = 0.9104$. }
    \end{subfigure}
    \begin{subfigure}[t]{\linewidth}
        \includegraphics[width=\textwidth]{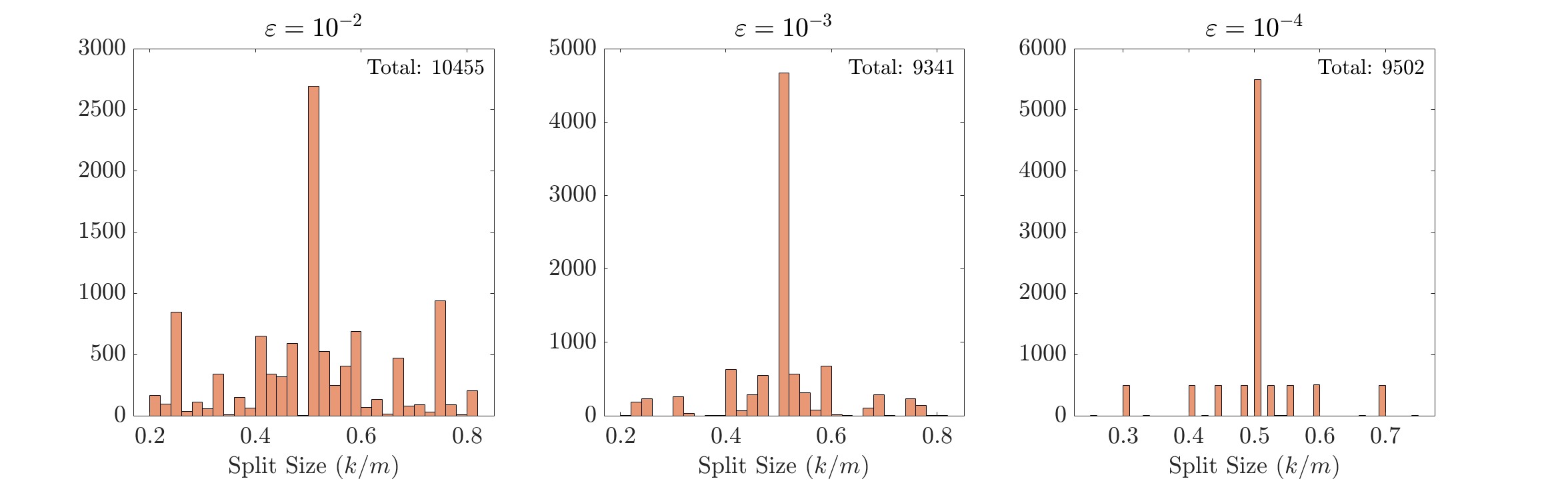}
        \caption{Frequency of relative eigenvalue split sizes (i.e., $k/m$ in \textbf{EIG}). We do not include subproblems with $m \leq 3$, as these can only be split in one way. Since the total number of splits is variable (and dependent on the splits themselves) we record it at the top of each plot. Dividing this total by 500 gives a rough average number of splits per run.}
    \end{subfigure}
    \begin{subfigure}[t]{\linewidth}
        \includegraphics[width=\textwidth]{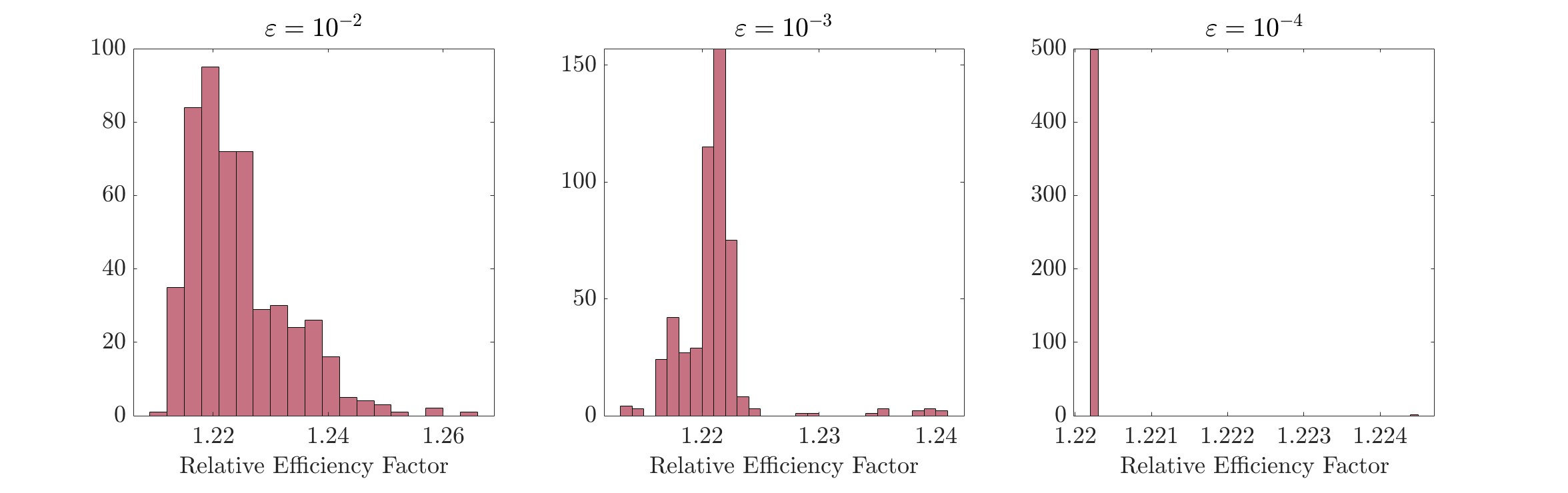}
        \caption{Frequency of relative efficiency factor (i.e., the pseudo-flop count divided by its optimal value).}
    \end{subfigure}
    \caption{Performance data for \textbf{RPD} on the $50 \times 50$ planted-spectrum example with decreasing values of $\varepsilon$. Each plot corresponds to 500 runs of \textbf{RPD}.}
    \label{fig: increasing_accuracy}
\end{figure}

\begin{figure}[p]
    \centering
    \begin{subfigure}[t]{\linewidth}
        \includegraphics[width=\textwidth]{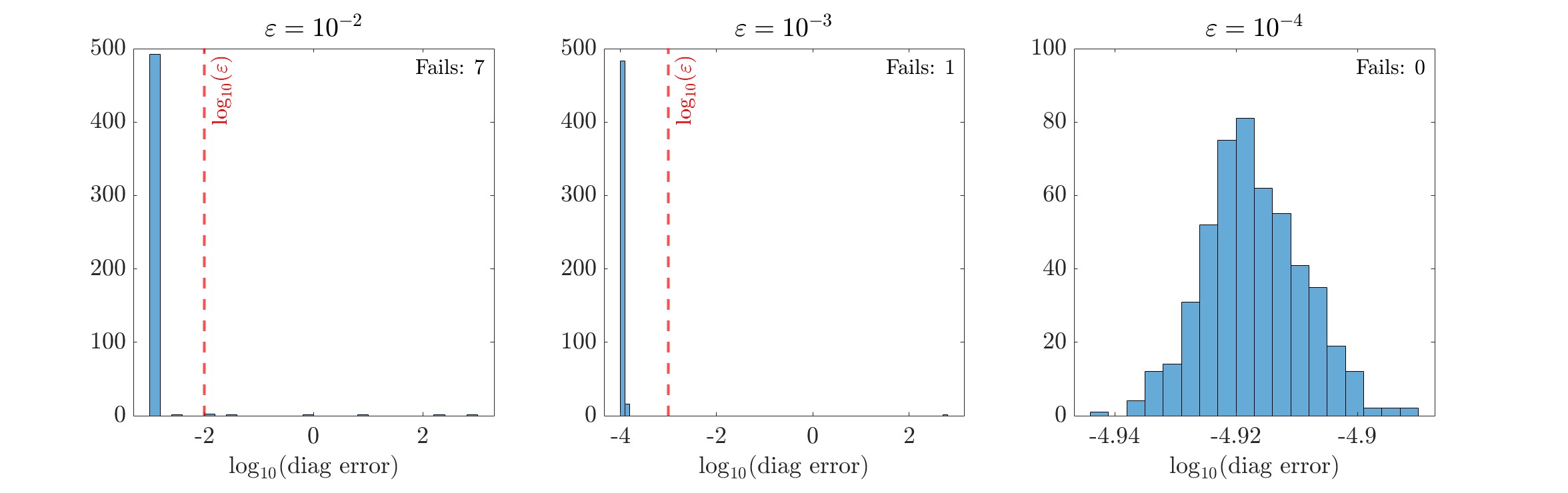}
        \caption{Same as (a) of \cref{fig: increasing_accuracy}. For this example $||A||_2 = ||B||_2 = 1$.}
    \end{subfigure}
    \begin{subfigure}[t]{\linewidth}
        \includegraphics[width=\textwidth]{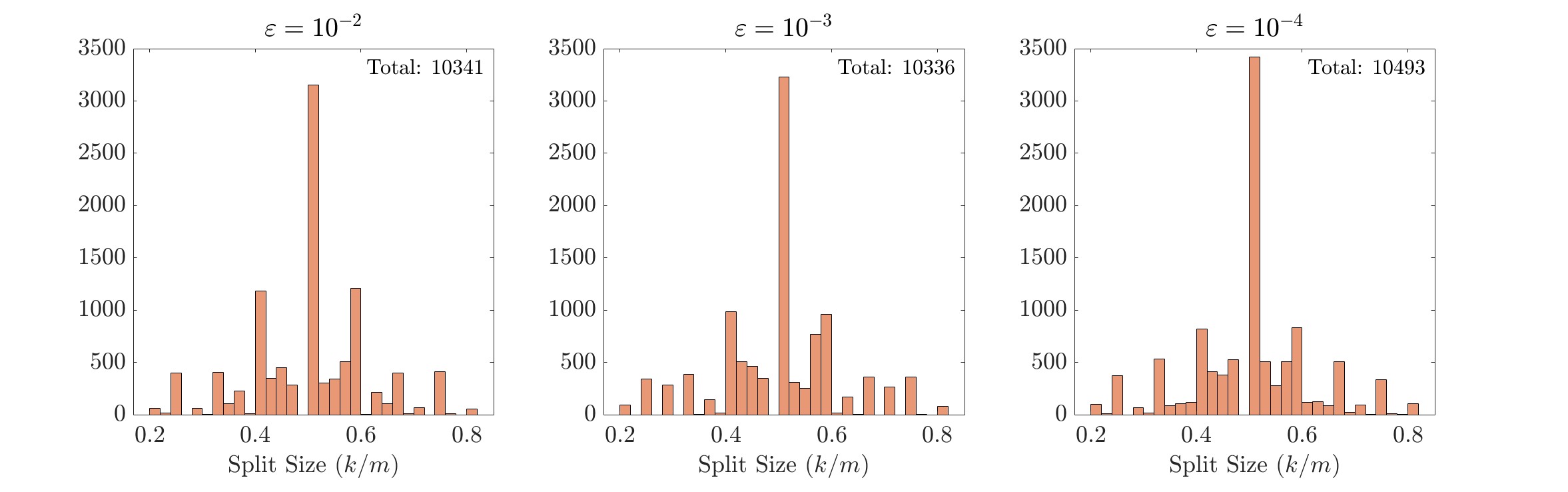}
        \caption{Same as (b) of \cref{fig: increasing_accuracy}.}
    \end{subfigure}
    \begin{subfigure}[t]{\linewidth}
        \includegraphics[width=\textwidth]{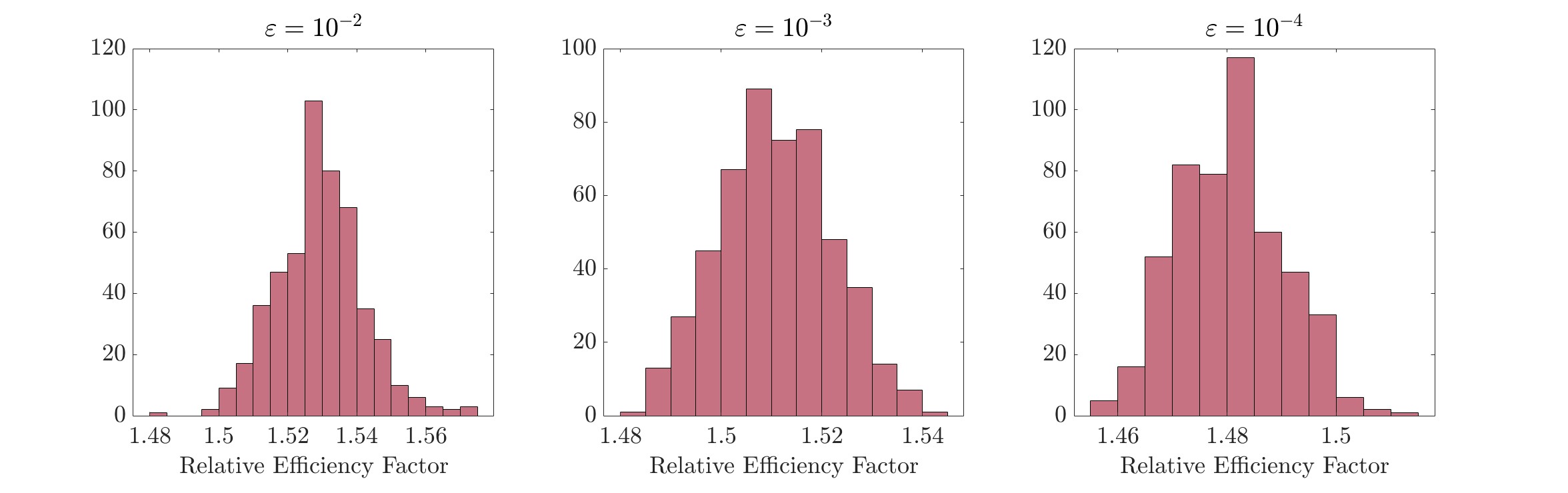}
        \caption{Same as (c) of \cref{fig: increasing_accuracy}}
    \end{subfigure}
    \caption{A repeat of \cref{fig: increasing_accuracy} for the $50 \times 50$ Jordan block example.}
    \label{fig: jordan_block}
\end{figure}

\indent We track the performance of divide-and-conquer on these examples in several ways. First and foremost, we want to verify a finite-precision counterpart to \cref{thm: RPD guarantees}:\ does \textbf{RPD} reliably produce accurate diagonalizations of each pencil? With this in mind, we compute diagonalization error as
\begin{equation}
    \log_{10} \left( \max \left\{ ||A - SDT^{-1}||_2, ||B - ST^{-1}||_2 \right\} \right),
    \label{eqn: diag error}
\end{equation}
for $S,D$, and $T$ the outputs of \textbf{RPD}, and we consider a run to be successful if this error is at most $\log_{10}(\varepsilon)$. Tracking the number of failed runs yields an empirical failure probability for \textbf{RPD} (with the simplifications made above). Note that \eqref{eqn: diag error} is only meaningful if $||A||_2$ and $||B||_2$ are roughly equal and close to one, as is the case in our examples. \\
\indent Next, we want to measure the efficiency of the divide-and-conquer process. One way to do this is to catalog the relative split size at each step (i.e., $k/m$ in \textbf{EIG}). While we know that \textbf{EIG} guarantees that the relative split is at least $0.2$ and at most $0.8$, divide-and-conquer is most efficient if relative splits are close to 0.5 at each step. \\
\indent Of course, the split size tells only part of the story. Recalling the proof of \cref{prop: op count}, \textbf{EIG} spends most of its time finding a dividing line; thus, even if splits are reliably near 50/50, the algorithm may be slow if too many lines are checked. Assuming access to $O(n^3)$ algorithms for matrix multiplication and QR, one step of our implementation requires $O(\log( \frac{n}{\varepsilon}) m^3 l)$ operations, where $l$ is the number of grid lines checked and $m$ is the size of the current subproblem.  Ignoring the $\log ( \frac{n}{\varepsilon} )$ factor (as it will cancel in our eventual measure of efficiency) we can do a pseudo-flop count by summing $m^3l$ over all steps with $m > 1$, and we can easily compute an optimal value for this count by requiring at each step that $l = 1$ and that the split is as close to 50/50 as possible. Dividing the actual count by the optimal one produces what we call a \textit{relative efficiency factor} for each run, which tells us roughly how many times more work \textbf{RPD} is doing than the best-case scenario. \\
\indent Histograms of each of these measures of performance are presented for both examples in \cref{fig: increasing_accuracy,fig: jordan_block}. In each test, we run \textbf{RPD} 500 times and present results for decreasing values of $\varepsilon$. With only a handful of failed runs on each problem, the results are compelling:\ \textbf{RPD} reliably diagonalizes both pencils, and divide-and-conquer appears to favor near-optimal eigenvalue splits. This carries through to the relative efficiency. Our rough flop count shows that \textbf{RPD} executes only slightly more than the optimal amount of work to produce these diagonalizations. Note that the number of failed runs appears to decrease with $\varepsilon$. While this may seem counterintuitive, it is a byproduct of our relaxed parameters, which become more restrictive (or equivalently more sensitive) as $\varepsilon$ shrinks.

\begin{figure}[t]
    \centering
    \begin{subfigure}[t]{\linewidth}
        \includegraphics[width=\textwidth]{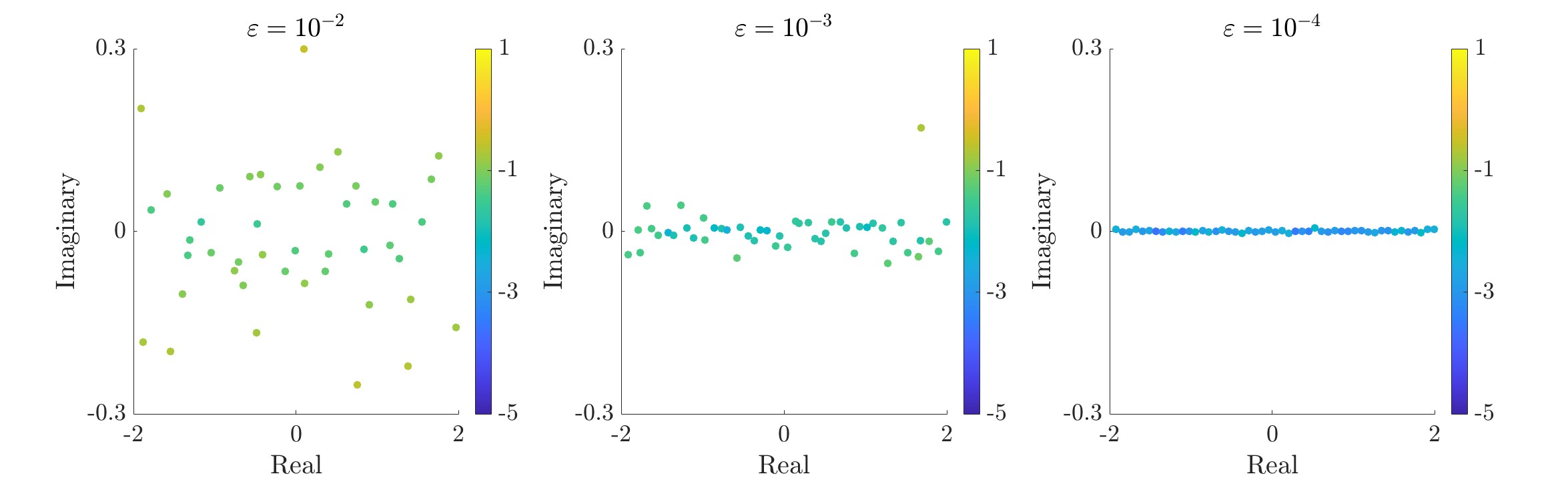}
        \caption{Eigenvalue approximations for the planted spectrum example obtained from \textbf{RPD} with different values of $\varepsilon$. Each approximate eigenvalue is colored according to its accuracy, which is computed as $\log_{10}| \widetilde{\lambda}_i - \lambda_i|$ for $\lambda_i$ and $\widetilde{\lambda}_i$ the true and approximate eigenvalues ordered by their real parts.}
    \end{subfigure}
    \begin{subfigure}[t]{\linewidth}
        \includegraphics[width=\textwidth]{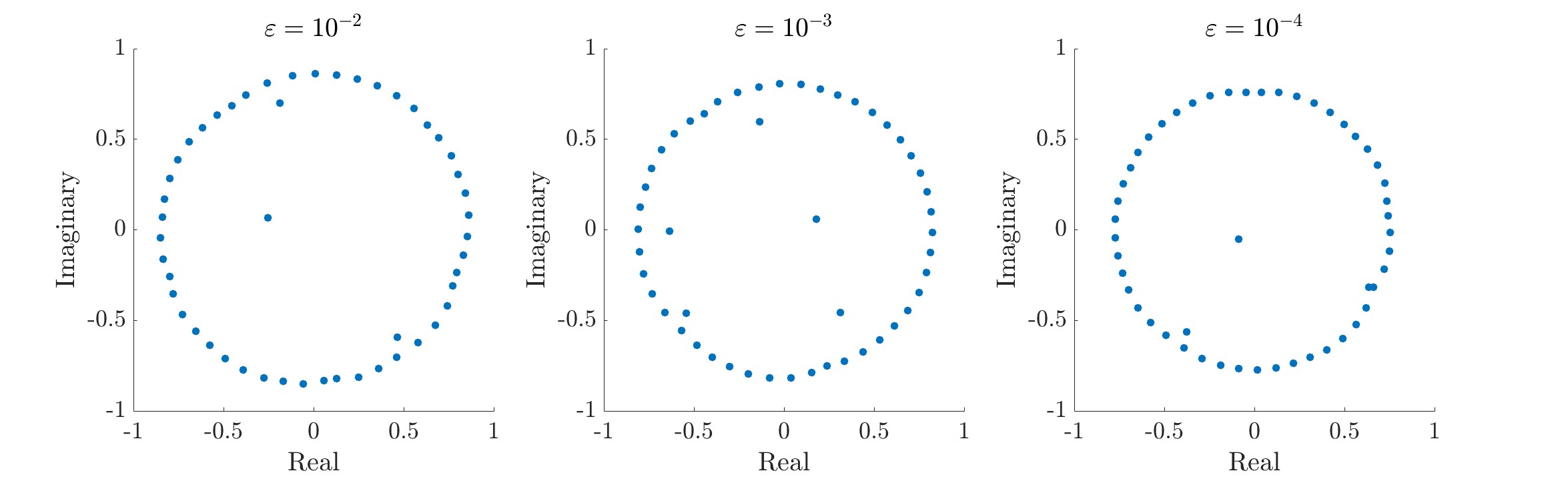}
        \caption{Eigenvalue approximations for the Jordan block example obtained from \textbf{RPD} with different values of $\varepsilon$. The only true eigenvalue for this problem is zero.}
    \end{subfigure}
    \caption{Eigenvalue approximation data for \textbf{RPD} applied to the planted spectrum and Jordan block examples. We present here only results for successful runs -- i.e., each set of approximate eigenvalues corresponds to a diagonalization with $||A - SDT^{-1}||_2 \leq \varepsilon$ and $||B - ST^{-1}||_2 \leq \varepsilon$.}
    \label{fig:eig_approx_data}
\end{figure}

\indent While these results are promising, we might be more interested in probing the boundaries of \cref{thm: forward error bound}. That is, when \textbf{RPD} succeeds, how accurate are the corresponding sets of approximate eigenvalues? With this in mind \cref{fig:eig_approx_data} provides eigenvalue approximation data for both model problems, where in each case we consider only approximations produced by successful runs. The results in these plots trace out nicely the challenge of extracting accurate eigenvalues from an accurate diagonalization. In the best case -- the planted spectrum example -- increasingly accurate diagonalizations provide correspondingly better eigenvalue approximations, as promised by \cref{thm: forward error bound}. The Jordan block example, on the other hand, demonstrates that when $\text{gap}(A,B) = 0$ we cannot hope to recover repeated eigenvalues with any confidence, though this is also the case for classical backwards-stable algorithms like QZ.

\subsection{Large $n$ and Infinite Eigenvalues}
Running \textbf{RPD} down to subproblems of size $1 \times 1$ is useful from a theoretical perspective but unlikely to be done in practice. We turn next to a more realistic use case, where $n$ is large and only a few splits are made before passing off to QZ. In this setting, we test the algorithm on pencils with an eigenvalue at infinity (i.e., where $B$ is singular). To do this, we construct a $1000 \times 1000$ pencil by drawing $A$ and $B$ randomly, computing a singular value decomposition $B = U\Sigma V^H$, and setting 
\begin{equation}
    B = B - \sigma_{\min}(B) uv^H
    \label{eqn: make B singular}
\end{equation}
for $u$ and $v$ the last columns of $U$ and $V$, respectively. By construction, this forces $B$ to be singular without changing its remaining singular values (and critically ensuring its norm remains comparable to $A$). As in the previous examples, we normalize $A$ and $B$ before calling \textbf{RPD}.
\begin{remark}
    Because we omit the $n^{\alpha}$ scaling on this example, the perturbed pencil $(\widetilde{A}, \widetilde{B})$ is very likely to have an eigenvalue outside of the shattering grid. While this appears problematic it is ultimately harmless. When only a handful of splits are needed, $(\widetilde{A}, \widetilde{B})$ can have many -- even a large fraction -- of eigenvalues outside of the grid, and divide-and-conquer will only falter if these fall predominantly in a specific region (for example between two vertical grid lines but above all of the horizontal ones). In part, this justifies our choice to omit the $n^{\alpha}$ scaling; while it is necessary to state an algorithm that provably runs to scalar subproblems, it is overly restrictive in practice -- driving eigenvalues that are initially small closer together and necessitating a finer grid. 
\end{remark}

\begin{figure}[t]
    \centering
    \includegraphics[width=\linewidth]{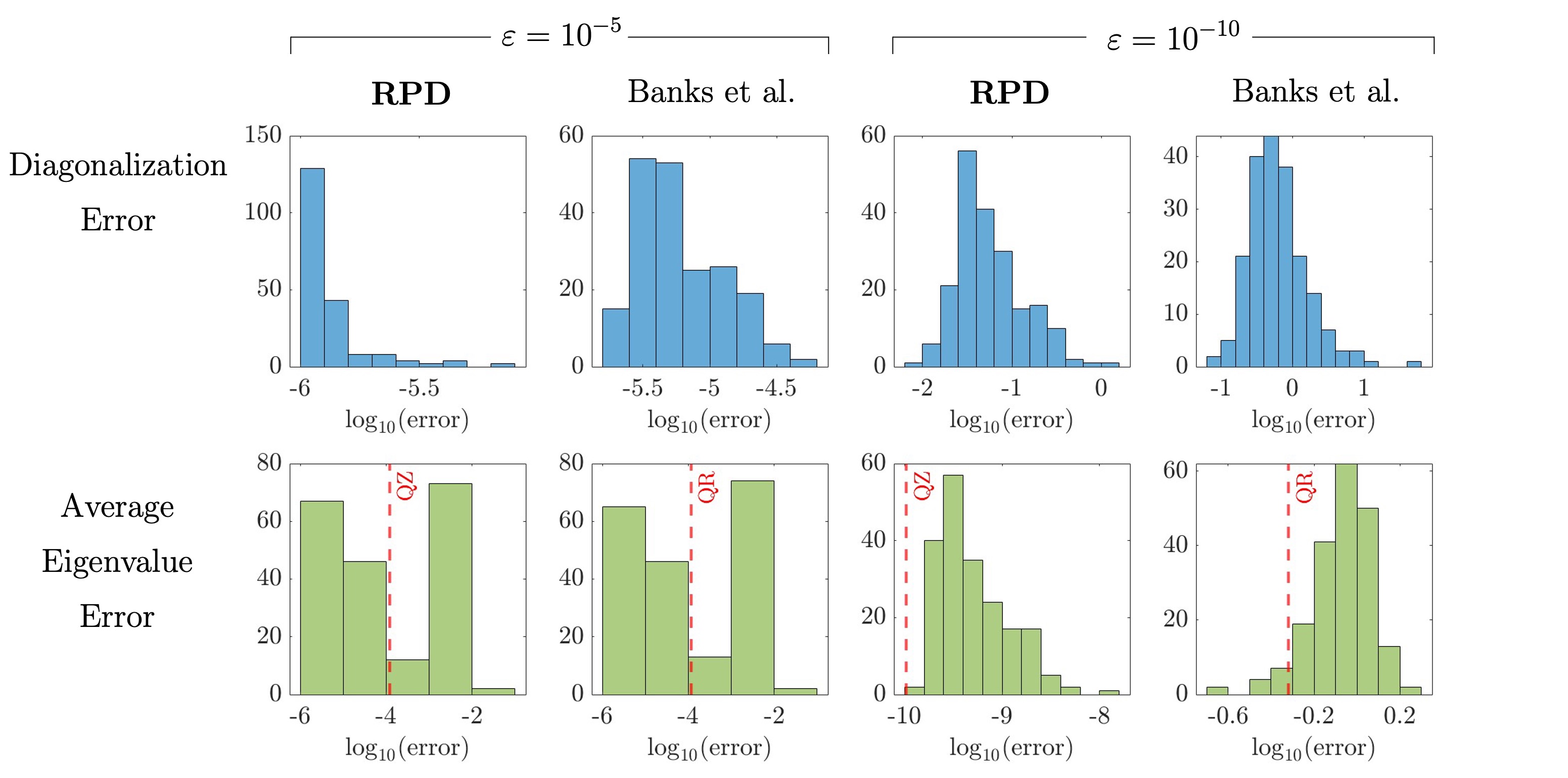}
    \caption{Performance data for \textbf{RPD} versus the single matrix, divide-and-conquer algorithm of Banks et al. \cite{banks2020pseudospectral}. Here, $1000 \times 1000$ pencils are constructed by drawing $A$ and $B$ randomly and subtracting a rank one matrix from $B$ to force it to be singular (without changing its remaining singular values). For twenty random draws of $A$ and $B$, we present ten runs through both algorithms, with first $\varepsilon = 10^{-5}$ and subsequently $\varepsilon = 10^{-10}$. Diagonalization error is computed according to \eqref{eqn: diag error} while eigenvalue accuracy is measured by ordering the true and approximate eigenvalues by magnitude and computing the average absolute error, excluding the eigenvalue at infinity. For the latter, we mark the  error achieved by QZ/QR (computed in the same way) when applied to  $(\widetilde{A},\widetilde{B})$ and $\widetilde{B}^{-1} \widetilde{A}$, respectively.}
    \label{fig:alg_comp}
\end{figure}

In addition to testing larger values of $n$, we also use this example to further justify our decision to avoid matrix inversion. To that end -- and in light of \cref{prop: shattering_for_prod} -- we compare \textbf{RPD} to an alternative algorithm that proceeds as follows:
\begin{enumerate}
    \item Perturb to obtain $(\widetilde{A}, \widetilde{B})$.
    \item Form the product $X = \widetilde{B}^{-1} \widetilde{A}$.
    \item Apply the single-matrix, divide-and-conquer routine of Banks et al.\ \cite[Algorithm EIG]{banks2020pseudospectral}.
\end{enumerate}
To ensure a fair comparison, we run both algorithms with the same perturbations (meaning the same $\widetilde{A}, \widetilde{B}$) and the same grid. We also restrict the number of steps of the Newton iteration for the sign function -- the counterpart to \textbf{IRS} used by Banks et al.\ -- to $\lceil \log_2(n / \epsilon) \rceil$. In both cases, we run divide-and-conquer to subproblems with $m \leq 250$, at which point we default to QZ and QR, respectively, by calling Matlab's \texttt{eig}. In the tests presented below, both algorithms averaged five splits before defaulting to \texttt{eig}. \\
\indent  Once again, we track diagonalization error via \eqref{eqn: diag error}. Importantly, we are interested in producing a diagonalization 
 of $(A,B)$ even if the single-matrix algorithm of Banks et al.\ is used. In that case, a matrix $T$ containing right eigenvectors (the output of single-matrix divide-and-conquer) gives rise to a corresponding matrix of left eigenvectors $S = \widetilde{B}T$, which together approximately diagonalize $(A,B)$. This mirrors exactly how the diagonalization is obtained in \textbf{RPD}. \\
\indent In addition, we record the eigenvalue error associated with each diagonalization. This is done by ordering the approximate and true eigenvalues by magnitude and computing the average, absolute error over the spectrum, excluding the eigenvalue at infinity. \cref{fig:alg_comp} records eigenvalue and diagonalization errors for both algorithms and for two choices of $\varepsilon$. Each plot corresponds to 200 trials, derived from twenty random draws of $A$ and $B$ run through each algorithm ten times.

When $\varepsilon = 10^{-5}$ we see little difference between the algorithms; though \textbf{RPD} is slightly more accurate, both approaches produce successful diagonalizations and comparably good eigenvalue approximations. Recalling that the size of the perturbation is determined by $\varepsilon$, this appears to be a consequence of the relatively large perturbation applied to $B$, which ensures that $\widetilde{B}$ is well-conditioned and that the error incurred by forming $\widetilde{B}^{-1}\widetilde{A}$ is small. In this case -- or more generally in situations where $B$ is known to be well-conditioned -- both algorithms are viable and essentially equivalent. \\
\indent In contrast, when $\varepsilon$ is much smaller \textbf{RPD} shows clear advantages. While neither can produce an accurate enough diagonalization -- an indication that our relaxed parameters are too loose for this regime -- \textbf{RPD} is consistently an order of magnitude better than its single-matrix alternative, and its corresponding eigenvalue approximations are remarkably accurate. Again this seems attributable to the conditioning of $\widetilde{B}$; here, forming $\widetilde{B}^{-1}\widetilde{A}$ not only forces divide-and-conquer to work with  a poorly-conditioned matrix but also introduces error that meaningfully shifts the eigenvalues away from those of $(A,B)$. \\
\indent Indeed, the difference in eigenvalue error present when $\varepsilon = 10^{-10}$ traces a similar gap between QZ and QR, as marked on the histograms. This indicates that the poor eigenvalue recovery of Banks et al.\ is due primarily to the gap between the eigenvalues of $\widetilde{B}^{-1}\widetilde{A}$ and $(A,B)$, which is essentially what the error in QR represents. Consequently, we cannot hope that by improving the diagonalization produced by Banks et al.\ -- which should be possible by adjusting the parameters -- we will see a similar improvement in the eigenvalues. On this example, then, we expect that this approach will break down as $\varepsilon$ becomes small; while initially decreasing $\varepsilon$ may improve diagonalization and eigenvalue accuracy, $\widetilde{B}$ will eventually become poorly-conditioned enough to guarantee that the eigenvalues of $\widetilde{B}^{-1}\widetilde{A}$ are far from those of $(A,B)$. Contrast this with \textbf{RPD}, which does not suffer from the same drawback and can reproduce the finite eigenvalues of $(A,B)$ to high accuracy.  \\
\indent Note that the example considered here is a larger version of the one covered by \cref{fig: zoomed_pseudospectra_comp2}. Together, they capture the danger of operating with inversion:\ when $\widetilde{B}$ is poorly conditioned, not only are the pseudospectra of $\widetilde{B}^{-1}\widetilde{A}$ unwieldy, but the eigenvalues they collapse to may significantly stray from those of the input pencil. 

\subsection{Singular Pencils}
\indent To this point we have exclusively tested divide-and-conquer on regular pencils. Since singular pencils are a reality in many applications, we consider as a final test the following example taken from Lotz and Noferini \cite{wilkinson_bus}:
\begin{equation}
     A = \begin{pmatrix} 2 & - 1 & -5 & -1 \\
    6 & -2 & -11 & -2 \\
    5 & 0 & -2 & 0  \\
    3 & 1 & 3 & 1 \end{pmatrix}, \; \; B = \begin{pmatrix} 1 & -1 & -4 & -2 \\
    2 & -3 & -12 & -6 \\
    -1 & -3 & -11 & -6 \\
    -2 & -2 & -7 & -4 \end{pmatrix} .
    \label{eqn: singular_example}
\end{equation}
The singular structure of this pencil is revealed by its Kronecker canonical form \cite{kronecker}:
\begin{equation}
    A - \lambda B = \begin{pmatrix} -3 & 1 & 1 & 1 \\ -8 & 3 & 2 & 0 \\ -5 & 3 & 0 & 1 \\ -2 & 2 & -1 & 1 \end{pmatrix} \begin{pmatrix} 1 - \lambda & 0 & 0 & 0 \\ 0 & -\lambda & 1 & 0 \\ 0 & 0 & - \lambda & 1 \\ 0 & 0 & 0 & 0 \end{pmatrix} \begin{pmatrix} - 1 & 0 & 1 & 0 \\ 0 & 1 & -4 & -2 \\ 0 & 0 & 1 & 0 \\ 1 & -1 & 0 & 1 \end{pmatrix}^{-1}
    \label{eqn: KCF},
\end{equation}
which indicates that $(A,B)$ has only one simple eigenvalue at $\lambda = 1$. \\
\indent Practically speaking, recovering this eigenvalue via divide-and-conquer should be difficult due to the initial perturbation made by \textbf{RPD}. In fact, Lotz and Noferini show that an arbitrarily small (nonrandom) perturbation to $(A,B)$ can send its eigenvalues to any four points in the complex plane. In spite of this, divide-and-conquer finds the true eigenvalue to roughly five digits of precision when run with $\varepsilon = 10^{-6}$, as shown in \cref{table:singular_pencil}. 

\renewcommand{\arraystretch}{1.2}
\begin{table}[t]
\begin{center}
\begin{tabular}{S[table-format=-1.6] llll}
\hline
\hline
 &
\multicolumn{3}{c}{\hspace{.28cm}Divide-and-Conquer (\cref{alg: RPD})} \\          
\multicolumn{1}{c}{\hspace{.2cm}QZ\cite{QZ}} & \multicolumn{1}{c}{\hspace{.28cm}Run 1} & \multicolumn{1}{c}{\hspace{.28cm}Run 2} & \multicolumn{1}{c}{\hspace{.28cm}Run 3}\\
\hline
-2.089013 & \hspace{.28cm}$0.999997 - 0.000007i$ & \hspace{.28cm}$1.000001 - 0.000008i$ & \hspace{.28cm}$0.999994 - 0.000032i$ \\
1.000000 & \hspace{.28cm}$0.463026 - 0.216072i$ & $-0.064526 + 0.383530i$ & \hspace{.28cm}$0.816203 - 0.082329i$ \\
0.445724 & $-0.036567 - 0.391359i$ & $-0.236744 - 0.298458i$ & \hspace{.28cm}$0.178283 - 0.299444i$ \\
-0.014976 & $-1.987468 - 0.361121i$ & $-1.152406 + 0.981613i$ & $-2.368016 + 0.144692i$ \\
\hline 
\hline
\end{tabular}
\end{center}
\caption{Eigenvalues of the singular pencil \eqref{eqn: singular_example} as computed by QZ and pseudospectral divide-and-conquer. We present three successful runs of \textbf{RPD}, each with $\varepsilon = 10^{-6}$.} \label{table:singular_pencil}
\end{table}

\begin{figure}[t]
    \centering
    \includegraphics[width=\linewidth]{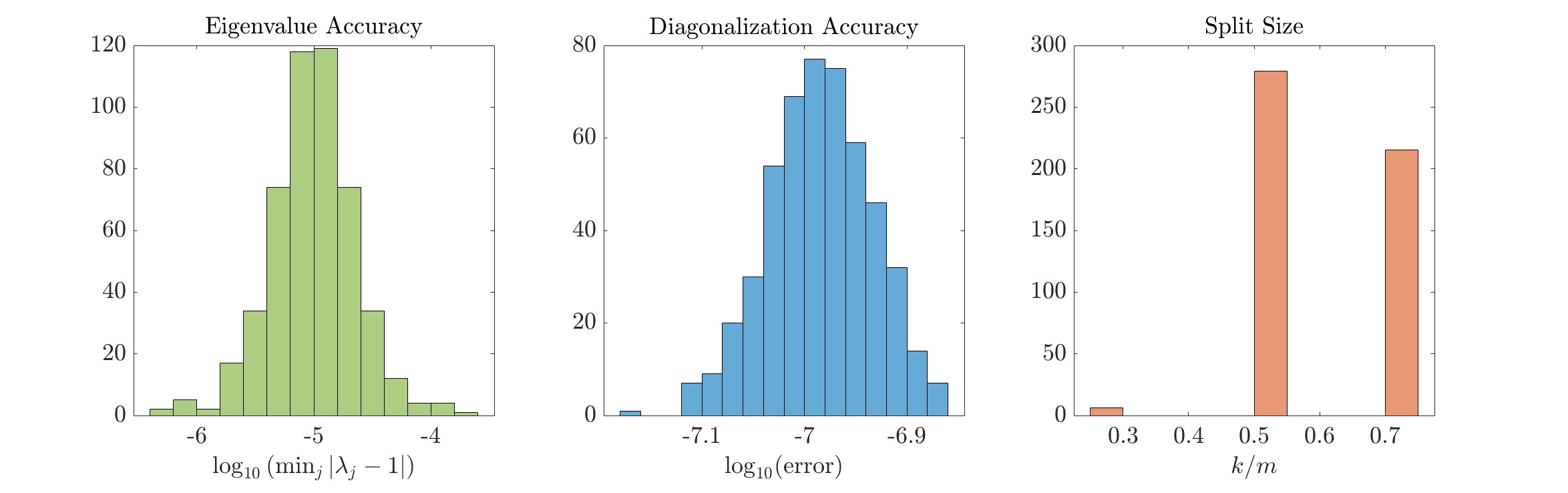}
    \caption{Statistics for 500 runs of \textbf{RPD} applied to a normalized version of the singular pencil \eqref{eqn: singular_example} with $||A||_2 = 0.6940$, $||B||_2 = 1$, and $\varepsilon = 10^{-6}$. For the first plot, $\left\{ \lambda_i\right\}_{1 \leq i \leq 4}$ are the approximate eigenvalues obtained via divide-and-conquer. Once again, no runs defaulted to calling \texttt{eig}. As in the previous examples, we consider only splits on subproblems with $m > 3$ for the third histogram; as a result, it records only the first split of each run, and the total number of splits is 500 (the number of runs).}
    \label{fig: singular_pencil_histograms}
\end{figure}

\indent Of course, QZ also finds $\lambda = 1$ to a remarkable 14 digits of precision. Motivated by this observation, Lotz and Noferini develop a theory of weak condition numbers to help explain the apparent stability of this eigenvalue. The results in \cref{table:singular_pencil} reflect the spirit of these condition numbers; though perturbations exist that produce arbitrary eigenvalues, a typical random one is likely to yield an eigenvalue near one.\footnote{Earlier work of Demmel and Kågström \cite{DEMMEL1987139} took a different approach to establishing the stability of QZ on singular pencils, demonstrating that QZ will preserve Kronecker canonical form, and therefore successfully recover true eigenvalues, provided round-off errors are small enough. This nevertheless cannot explain the strong performance of \textbf{RPD}, as Gaussian perturbations \textit{will} change the Kronecker structure of the pencil with high probability.}  Randomization has another advantage here:\ while QZ always produces the same three spurious ``eigenvalues," divide-and-conquer does not -- meaning multiple runs can help distinguish true eigenvalues from fake. Finally, we provide in \cref{fig: singular_pencil_histograms} a few empirical statistics for 500 runs of \textbf{RPD} on a normalized version of $(A,B)$.

\section{Conclusion}
In this paper, we have constructed an inverse-free routine for producing an approximate diagonalization of any matrix pencil $(A,B)$. The method, which rests on a divide-and-conquer eigensolver for the generalized eigenvalue problem, involves several threads of randomness. The most important of these is a random perturbation made to the input matrices, which regularizes the problem and allows the divide-and-conquer process to succeed with high probability. Leveraging this result, our approach produces a provably accurate approximate diagonalization (in a backward-error sense) for any input pencil in exact arithmetic, again with high probability. Moreover, the algorithm runs in nearly matrix multiplication time, establishing $O(\log^2 ( \frac{n}{\varepsilon}) T_{\text{MM}}(n))$ as an upper bound on the complexity of producing -- without relying on inversion -- an exact-arithmetic diagonalization of an $n \times n$ pencil with backward error $\varepsilon$. Replacing $B$ with the identity, this provides another nearly matrix multiplication time algorithm for diagonalizing a single matrix, representing an inverse-free counterpart to the algorithm of Banks et al.\ \cite{banks2020pseudospectral}. \\
\indent For interested readers, we note that there are a number of open problems and next steps related to this work, which we discuss briefly below.
\begin{itemize}
    \item \textbf{Finite-Precision Analysis:} The most immediate open question concerns a floating-point stability bound for \textbf{RPD}. While we expect that the method will succeed with similar guarantees in finite precision if stated properly, we have not yet completed an end-to-end analysis. \cref{appendix: finite arithmetic} explores some of the work done in this direction.
    \item \textbf{Hermitian/Deterministic Variants:} Throughout, we made no assumptions on the matrices of our pencil. Nevertheless, much attention has been paid to the Hermitian eigenvalue problem, where $A$ is Hermitian and $B$ is Hermitian positive definite, meaning a version of the method adapted to that setting may be worth developing. The main difference there is the use of \textbf{DEFLATE}, as the left and right eigenspaces of a Hermitian pencil are the same. More generally, we might hope to devise a deterministic way to produce the pseudospectral shattering result central to the algorithm.
    \item \textbf{Perturbation Results for the Generalized Eigenvalue Problem:} In completing this work, we found that many of the useful perturbation results for the single-matrix eigenvalue problem do not have counterparts for the generalized problem. In particular, there is a lack of results stated in ways that are practical for numerical purposes (i.e., using the Euclidean norm as opposed to a chordal metric). We hope that this work inspires more research in this direction, and we expect that developing better tools may improve the analysis presented here. 
\end{itemize}
\section{Acknowledgements}
This work was supported by NSF grant DMS 2154099 and by Graduate Fellowships for STEM Diversity (GFSD). We thank Jess Banks, Jorge Garza-Vargas, and Nikhil Srivastava for their advice and feedback as well as Anne Greenbaum for helpful correspondence. We also thank two anonymous referees for their detailed comments, corrections, and suggestions, which have improved the quality of the paper. 

\appendix
\section{Alternative Versions of Bauer-Fike for Matrix Pencils}
\label{appendix: BF}

In this appendix, we discuss alternative formulations of \cref{thm: pencil BF}. As mentioned above, our version of Bauer-Fike is stated for easy use in the proof of \cref{thm: shattering} and therefore appears somewhat unusual compared to existing perturbation results for the generalized eigenvalue problem. Nevertheless, it is essentially the same as standard bounds, as we demonstrate below. To do this, we introduce common notation from perturbation theory as presented in Stewart and Sun \cite{stewart1990matrix}.  \\
\indent We begin by re-casting each eigenvalue $\lambda$ as an ordered pair $\langle \alpha, \beta \rangle$ such that $\lambda = \frac{\alpha}{\beta}$. In this notation, solutions to the generalized eigenvalue problem satisfy $\beta Av = \alpha B v$. The main benefit of representing eigenvalues in this way is the natural inclusion of those at infinity, which now correspond to $\beta = 0$. Moreover, this definition links nicely with matrix pencil diagonalization; if $(A,B) = (SD_1T^{-1}, SD_2T^{-1})$, each diagonal entry of $D_1$ and $D_2$ contains a corresponding $\alpha$ and $\beta$. Note that if $(A,B)$ is regular, we can never have $\alpha = \beta = 0$.  \\
\indent Of course, there are infinitely many choices of $\alpha$ and $\beta$ that yield the same eigenvalue $\lambda$. For this reason, we can think of each eigenvalue $\langle \alpha, \beta \rangle$ as a projective line -- i.e., the subspace spanned by $(\alpha \; \beta)^T$. Accordingly, we measure distance between two eigenvalues $\langle \alpha_1, \beta_1 \rangle$ and $\langle \alpha_2, \beta_2 \rangle$ with the chordal metric:
\begin{equation}
    \chi \left( \langle \alpha_1, \beta_1 \rangle, \langle \alpha_2, \beta_2 \rangle \right) = \frac{|\alpha_1 \beta_2 - \beta_1 \alpha_2|}{\sqrt{ |\alpha_1|^2 + |\beta_1|^2} \sqrt{|\alpha_2|^2 + |\beta_2|^2}}.
    \label{eqn: chordal metric}
\end{equation}
Dividing both the numerator and denominator of \eqref{eqn: chordal metric} by $|\beta_1 \beta_2|$, we observe that if $\lambda_1 = \frac{\alpha_1}{\beta_1}$ and $\lambda_2 = \frac{\alpha_2}{\beta_2}$ then 
\begin{equation}
    \chi \left( \langle \alpha_1, \beta_1 \rangle, \langle \alpha_2, \beta_2 \rangle \right) =  \frac{|\lambda_1 - \lambda_2|}{\sqrt{|\lambda_1|^2 + 1}\sqrt{|\lambda_2|^2 + 1}} .
    \label{eqn: chordal metric in terms of points}
\end{equation}
In other words, the chordal distance between $\langle \alpha_1, \beta_1 \rangle$ and $\langle \alpha_2, \beta_2 \rangle$ is half the Euclidean distance between the images of $\lambda_1$ and $\lambda_2$ under the stereographic projection. This ensures that the distance between any two eigenvalues is at most one, including eigenvalues at infinity. \\
\indent Though perturbation results in terms of $\chi$ originate with Stewart \cite{Stewart_gershgorin}, the first version of Bauer-Fike stated in terms of the chordal distance for general matrices is due to Elsner and Sun \cite{ELSNER1982341}, which we summarize below.

\begin{thm}[Elsner and Sun 1982] \label{thm: E&S BF}
Let $(A,B)$ be a diagonalizable pencil with eigenvalues $\langle \alpha_i, \beta_i \rangle$. If $T$ is any matrix of right eigenvectors of $(A,B)$ and $\langle \widetilde{\alpha}_i, \widetilde{\beta}_i \rangle$ are the eigenvalues of the regular pencil $(\widetilde{A}, \widetilde{B})$, then
$$ \max_i \min_j \chi \left( \langle \alpha_i, \beta_i \rangle, \langle \widetilde{\alpha}_j, \widetilde{\beta}_j \rangle  \right) \leq \kappa(T) ||(AA^H + BB^H)^{-1/2}||_2 ||(A-\widetilde{A}, B - \widetilde{B})||_2 $$
\end{thm}

While \cref{thm: E&S BF} is more general than \cref{thm: pencil BF}, the two results are essentially the same under the assumptions made in \cref{section:pseudospectra}. To demonstrate this, suppose that $(A,B)$ is diagonalizable and that $B$ is nonsingular with $\epsilon < \sigma_n(B)$. If $z \in \Lambda_{\epsilon}(A,B)$ then $z$ is an eigenvalue of a pencil $(\widetilde{A}, \widetilde{B})$ with $||A - \widetilde{A}||_2, ||B - \widetilde{B}||_2 \leq \epsilon$. Thus, if $\lambda$ is the eigenvalue of $(A,B)$ closest to $z$, \cref{thm: E&S BF} and \eqref{eqn: chordal metric in terms of points} imply
\begin{equation}
    \frac{|z - \lambda|}{\sqrt{|z|^2 + 1}\sqrt{
|\lambda|^2 + 1}} \leq \sqrt{2} \epsilon \kappa(T) ||(AA^H + BB^H)^{-1/2}||_2
\label{eqn: same BF step one}
\end{equation}
or, equivalently,
\begin{equation}
    |z - \lambda| \leq \sqrt{2} \epsilon \kappa(T) ||(AA^H + BB^H)^{-1/2}||_2 \sqrt{|z|^2 + 1} \sqrt{|\lambda|^2 + 1} .
\label{eqn: same BF step two}
\end{equation}
Applying \cref{lem: pseudo upper bound} to bound $|z|$ and $|\lambda|$ and noting $||(AA^H + BB^H)^{-1/2}||_2 = \frac{1}{\sigma_n(A,B)}$, we conclude that $\Lambda_{\epsilon}(A,B)$ is contained in balls around the eigenvalues of $(A,B)$ of radius
\begin{equation}
    \sqrt{2} \epsilon \frac{ \kappa(T) }{\sigma_n(A,B)} \left[ 1 + \left( \frac{\epsilon ||B^{-1}||_2 + ||B^{-1}A||_2}{1 - \epsilon ||B^{-1}||_2} \right)^2 \right] .
    \label{eqn: same BF step three}
\end{equation}
\indent On one hand, the $\sigma_n(A,B)$ in \eqref{eqn: same BF step three} is an improvement over the factor of $||B^{-1}||_2$ in \cref{thm: pencil BF}, as $\sigma_n(A,B) \geq \sigma_n(B)$. At the same time, we pay the price in an additional factor of $\sqrt{2}$ and the square on the piece coming from \cref{lem: pseudo upper bound}. Since $\sigma_n(A,B) \geq \sigma_n(A)$, \eqref{eqn: same BF step three} is likely a better bound if we know that $A$ is far from singular, though in the proof of \cref{thm: shattering} the only bounds on $\sigma_n(A)$ available are the same as those for $\sigma_n(B)$. This means that, up to constants, a proof of shattering using \cref{thm: E&S BF} is the same as the one given above. In particular, Elsner and Sun's result cannot be used to improve the dependence of $\epsilon$ on $n$. Similarly, it would not change \cref{thm: forward error bound}. \\
\indent Nevertheless, there is an important benefit of the generality of \cref{thm: E&S BF}:\ the ability to control changes in the finite eigenvalues of a pencil with some at infinity. Up to this point, we have pursued a standard Bauer-Fike result that bounds $\Lambda_{\epsilon}(A,B)$ by a union of balls around the eigenvalues of $(A,B)$. Because we know $\Lambda_{\epsilon}(A,B)$ is unbounded if $B$ is singular or $\epsilon > \sigma_n(B)$, \cref{thm: pencil BF} and the alternate version \eqref{eqn: same BF step three} do not consider the case where $B$ is singular. This ultimately carries through to our forward error bound \cref{thm: forward error bound}, which similarly requires that $B$ is invertible. While we know that we can't obtain any kind of forward error guarantees for the infinite eigenvalues of $(A,B)$ with singular $B$, we would like to extend \cref{thm: forward error bound} to the finite eigenvalues if possible. \cref{thm: E&S BF} offers one way to do this, producing the following proposition.
\begin{prop}
    Let $(A,B)$ be a regular, diagonalizable pencil and let $(\widetilde{A}, \widetilde{B})$ be a nearby pencil with $||A - \widetilde{A}||_2, ||B - \widetilde{B}||_2 \leq \varepsilon$. Suppose that $\langle \alpha, \beta \rangle$ and $\langle \widetilde{\alpha}, \widetilde{\beta} \rangle$ are eigenvalues of $(A,B)$ and $(\widetilde{A}, \widetilde{B})$, where $\langle \alpha, \beta \rangle$ is the closest eigenvalue of $(A,B)$ to $\langle \widetilde{\alpha}, \widetilde{\beta} \rangle$ in the chordal metric $\chi$. If $\lambda = \frac{\alpha}{\beta}$ and $\widetilde{\lambda} = \frac{\widetilde{\alpha}}{\widetilde{\beta}}$ are both finite, then
     $$ |\lambda - \widetilde{\lambda}| \leq \sqrt{2} \varepsilon \frac{\kappa_V(A,B)}{\sigma_n(A,B)} \sqrt{|\lambda|^2+1}\sqrt{|\widetilde{\lambda}|^2 +1} .$$
\end{prop}
\begin{proof}
    This follows directly from \cref{thm: E&S BF}, noting that we can replace $\kappa(T)$ with $\kappa_V(A,B)$ by taking an infimum. 
\end{proof}
Of course, this result is only useful if we can bound $|\lambda|$ and $|\widetilde{\lambda}|$, which is not necessarily straightforward. We do note, however, that $\sigma_n(A,B)$ is guaranteed to be nonzero since $(A,B)$ is regular. \\
\indent To finish this section, we summarize other known variants of Bauer-Fike for the generalized eigenvalue problem (though this is by no means an exhaustive list). First, we note that \cref{thm: E&S BF} is the version stated by Stewart and Sun \cite[Theorem VI.2.6]{stewart1990matrix} rewritten in terms of the canonical angle between the row spaces of $(A,B)$ and $(\widetilde{A}, \widetilde{B})$. Stewart and Sun also include a version that swaps $\kappa(T)$ for $||T||_2||S^{-1}||_2$, where $(S^{-1}AT, S^{-1}BT)$ is a diagonalization of $(A,B)$ with the columns of $S^{-1}$ and $T$ normalized appropriately (see \cite[Theorem VI.2.7]{stewart1990matrix}). This yields the simpler bound
\begin{equation}
    \max_i \min_j \chi \left( \langle \alpha_i, \beta_i \rangle, \langle \widetilde{\alpha}_j, \widetilde{\beta}_j \rangle \right) \leq ||S^{-1}||_2||T||_2||(A-\widetilde{A}, B-\widetilde{B})||_2 .
\end{equation}
Note however that $||S^{-1}||_2$ hides a dependence on the conditioning of $A$ and $B$ (recall for example that we construct our diagonalization above by taking $S = BT$). \\
\indent Minor improvements on \cref{thm: E&S BF} can be found in work of Elsner and Lancaster \cite{Elsner_and_Lancaster}. Meanwhile, Chu \cite{Chu1} stated their own version of Bauer-Fike consisting of four separate bounds, depending on whether the initial and perturbed eigenvalues are finite/infinite, which was subsequently generalized to matrix polynomials \cite{Chu2}. Finally, we note a recent sharp version of Shi and Wei \cite{SHI20123218} stated in terms of the sign-complex spectral radius.  

\section{Finite-Precision Analysis}
\label{appendix: finite arithmetic}
In this appendix, we analyze key building blocks of \textbf{EIG} and \textbf{RPD} in finite-precision arithmetic. We first outline our model of floating-point computation before proving results for shattering and \textbf{GRURV}. A rigorous stability bound for \textbf{IRS} under similar assumptions is covered in \cite{My_thesis}. \\
\indent While a full end-to-end analysis of \textbf{RPD} remains open, these bounds suggest that it likely requires provably lower precision than methods built on inversion, in particular its single-matrix counterpart of Banks et al.\ (see the discussion in \cite[Chapter 6]{My_thesis}).
\subsection{Black-Box Assumptions}
 Throughout, we assume a floating-point arithmetic where 
\begin{equation}
    fl(x \circ y) = (x \circ y)(1 + \Delta), \; \; |\Delta| \leq {\bf u} 
    \label{eqn: floating point}
\end{equation}
for basic operations $\circ \in \left\{ +, -, \times, \div \right\}$ and a machine precision ${\bf u}(\varepsilon, n)$, which is a function of the desired accuracy $\varepsilon$ and the size of the problem $n$. We assume that $\sqrt{x}$ can be evaluated in this way as well. \eqref{eqn: floating point} is a standard formulation for finite-precision computations (see for example \cite{Higham}). \\
\indent Throughout, we can think of ${\bf u}$ as inverse polynomial in $n$. In practice, we choose ${\bf u}$ based on our error analysis to ensure either true backward stability or, if that's not possible, the slightly weaker logarithmic stability of Demmel, Dumitriu, and Holtz \cite{2007}. 
With this in mind, we further assume access to black-box algorithms for Gaussian sampling, matrix multiplication, and QR/QL, which we summarize below. These follow Banks et al.\ \cite{banks2020pseudospectral} to allow for easy connection with their floating-point bounds.

\begin{assump}[Gaussian Sampling]
    There exists a $c_{\text{N}}$-stable Gaussian sampler $\text{N}(\sigma)$ that takes $\sigma \in {\mathbb R}_{\geq 0}$ and outputs  $\text{N}(\sigma)$ satisfying $ |\text{N}(\sigma) - {\mathcal N}| \leq c_{\text{N}} \sigma {\bf u} $ for some ${\mathcal N} \sim {\mathcal N}(0, \sigma^2)$.
    \label{sampling assumption}
\end{assump}

\begin{assump}[Matrix Multiplication]
    There exists a $\mu_{\text{MM}}(n)$-stable $n \times n$ multiplication algorithm ${\bf MM}( \cdot, \cdot)$ satisfying 
    $$||{\bf MM}(A,B) - AB ||_2 \leq \mu_{\text{MM}}(n) {\bf u}||A||_2 ||B||_2 $$
    in $T_{\text{MM}}(n)$ arithmetic operations.
    \label{MM assumption}
\end{assump}

\begin{assump}[QR Factorization]
    There exists a $\mu_{\text{QR}}(n)$-stable $n \times n$ QR algorithm ${\bf QR}(\cdot)$ satisfying
    \begin{enumerate}
        \item $[Q,R] = {\bf QR}(A)$
        \item $R$ is exactly upper triangular
        \item There exist $A'$ and unitary $Q'$ such that $A' = Q'R$ with 
            $$ ||Q'-Q||_2 \leq \mu_{\text{QR}}(n) {\bf u} \; \; \text{and} \; \; ||A' - A||_2 \leq \mu_{\text{QR}}(n) {\bf u} ||A||_2, $$
    \end{enumerate}
    in $T_{\text{QR}}(n)$ arithmetic operations.
    \label{QR assumption}
\end{assump}
For simplicity, we further assume access to a QL algorithm $\textbf{QL}(\cdot)$ that satisfies similar guarantees as in \cref{QR assumption} with the same parameter $\mu_{\text{QR}}(n)$. Note that for these routines we can always guarantee a truly triangular result in finite-precision by forcing entries above/below the diagonal to be zero.  Finally, we include a black-box inversion algorithm to allow us to quantify the benefit of working inverse-free.
\begin{assump}[Matrix Inversion]
    There exists a $(\mu_{\text{INV}}(n), c_{\text{INV}})$-stable $n \times n$ inversion algorithm ${\bf INV}(\cdot)$ satisfying
    $$||{\bf INV}(A) - A^{-1}||_2 \leq \mu_{\text{INV}}(n) {\bf u} \kappa_2(A)^{c_{\text{INV}} \log(n)} ||A^{-1}||_2$$
    in $T_{\text{INV}}(n)$ arithmetic operations. 
    \label{INV assumption}   
\end{assump}
Before proceeding with our analysis, we pause to note that the assumptions listed above are all compatible with fast matrix multiplication \cite{2007}, meaning the bounds cover any implementation  that uses a fast matrix multiplication routine, including \cite{STRASSEN1969,COPPERSMITH1990251,fastermm,williams2023new}.

\subsection{Finite-Precision Shattering}
In this section, we consider how floating-point computations impact pseudospectral shattering. We start with $\Lambda_{\epsilon}(\widetilde{A},n^{\alpha}\widetilde{B})$. A straightforward extension of \cite[Theorem 3.13]{banks2020pseudospectral} implies the following finite-precision counterpart to \cref{thm: shattering}.
\begin{thm}\label{thm: finite_shattering}
    Let $A,B \in {\mathbb C}^{n \times n}$ with $||A||_2,||B||_2 \leq 1$ and let $0 < \gamma < \frac{1}{2}$. Further, let $\omega = \frac{\gamma^4}{4} n^{\frac{-8 \alpha + 13}{3}}$ and construct the grid $g = \text{\normalfont grid}(z, \omega, \lceil 8/ \omega \rceil, \lceil 8/ \omega \rceil)$ for $\alpha > 0$ and $z$ chosen uniformly at random from the square with bottom left corner $-4-4i$ and side length $\omega$. On a floating-point machine with precision ${\bf u}$, suppose $G_1,G_2 \in {\mathbb C}^{n \times n}$ satisfy $G_k(i,j) = \text{\normalfont N}(\frac{1}{\sqrt{n}})$ for $1 \leq i,j \leq n$ and $k = 1,2$. If $\widetilde{A} = A + \gamma G_1$ and $\widetilde{B} = B + \gamma G_2$ (again in finite precision) then $\Lambda_{\epsilon}(\widetilde{A}, n^{\alpha} \widetilde{B})$ is shattered with respect to $g$ for
    $$ \epsilon = \frac{1}{2} \cdot \frac{\gamma^5}{64 n^{\frac{11 \alpha + 25}{3}} + \gamma^5} $$
    with probability at least $ \left[ 1 - \frac{82}{n} - \frac{531441}{16n^2} \right] \left[ 1 - \frac{n^{2 - 2\alpha}}{\gamma^2} - 4e^{-n} \right] $ provided
    $$ {\bf u} \leq \frac{1}{2(3 + c_{\text{\normalfont N}})n^{\alpha + \frac{1}{2}}} \cdot \frac{\gamma^5}{64 n^{\frac{11 \alpha + 25}{3}} + \gamma^5} .$$
\end{thm}
\begin{proof}
    Following the proof of \cite[Theorem 3.13]{banks2020pseudospectral}, $\widetilde{A}$ and $\widetilde{B}$ are at most $(3 + c_{\text{N}})\sqrt{n}{\bf u}$ away (in the spectral norm) from their exact-arithmetic counterparts. Accommodating the $n^{\alpha}$ scaling on $\widetilde{B}$ and recalling that shattering is achieved in exact arithmetic for $\epsilon = \frac{\gamma^5}{64 n^{\frac{11 \alpha + 25}{3}} + \gamma^5}$ with probability at least $ \left[ 1 - \frac{82}{n} - \frac{531441}{16n^2} \right] \left[ 1 - \frac{n^{2 - 2\alpha}}{\gamma^2} - 4e^{-n} \right]$, \cref{lem: shattering preserved} implies that it is sufficient to take
    \begin{equation}
        (3 + c_{\text{N}})n^{\alpha + \frac{1}{2}} {\bf u} \leq \frac{1}{2} \cdot \frac{\gamma^5}{64 n^{\frac{11 \alpha + 25}{3}} + \gamma^5}, 
        \label{eqn: pencil_shattering_requirement}
    \end{equation}
    which is equivalent to the listed requirement on ${\bf u}$.
\end{proof}

We consider next $\widetilde{B}^{-1} \widetilde{A}$. In this case, we must account for error incurred by inverting $\widetilde{B}$ and multiplying by $\widetilde{A}$ on top of the error already baked into $\widetilde{A}$ and $\widetilde{B}$. With this in mind, we first prove an intermediate result.
\begin{lem}\label{lem: product_error}
    Suppose $A_1,A_2,B_1,B_2 \in {\mathbb C}^{n \times n}$ with $||A_1 - A_2||_2, ||B_1-B_2||_2 \leq \delta$ for $\delta < \sigma_n(B_1)$. Then
    $$ ||B_1^{-1}A_1 - B_2^{-1}A_2||_2 \leq \delta ||B_1^{-1}||_2 \left( 1 + \frac{||A_1||_2 + \delta}{\sigma_n(B_1)-\delta} \right). $$
\end{lem}
\begin{proof}
    We have
    \begin{equation}
        \aligned
        ||B_1^{-1}A_1 - B_2^{-1}A_2||_2 &= ||B_1^{-1}A_1 - B_1^{-1}A_2 + B_1^{-1}A_2 - B_2^{-1}A_2||_2 \\
        & \leq ||B_1^{-1}||_2||A_1 - A_2||_2 + ||B_1^{-1} - B_2^{-1}||_2||A_2||_2 \\
        & \leq \delta ||B_1^{-1}||_2 + ||B_1^{-1}||_2||B_2 - B_1||_2 ||B_2^{-1}||_2 ||A_2||_2 \\
        & \leq \delta ||B_1^{-1}||_2 \left( 1 + \frac{||A_2||_2}{\sigma_n(B_2)} \right) \\
        & \leq \delta ||B_1^{-1}||_2 \left( 1 + \frac{||A_1||_2 + \delta}{\sigma_n(B_1) - \delta} \right),
        \endaligned
        \label{eqn: expand_product}
    \end{equation}
    where the last inequality follows from stability of singular values. 
\end{proof}

Combining \cref{lem: product_error} with \cref{MM assumption} and \cref{INV assumption} yields the following finite-precision counterpart to \cref{prop: shattering_for_prod}.
\begin{thm}\label{thm: finite_product_shattering}
    Let $A,B \in {\mathbb C}^{n \times n}$ with $||A||_2, ||B||_2 \leq 1$ and let $\gamma \in (0,\frac{1}{2})$. Further let $\omega = \frac{\gamma^4}{4}n^{-\frac{8 \alpha+13}{3}}$ and construct the grid $g = \text{\normalfont grid}(z,\omega,\lceil 8/\omega \rceil, \lceil 8/\omega \rceil)$ for $\alpha > 0$ and $z$ chosen uniformly at random from the square with bottom left corner $-4 -4i$ and side length $\omega$. On a floating-point machine with precision ${\bf u}$, suppose $G_1,G_2 \in {\mathbb C}^{n \times n}$ satisfy $G_k(i,j) = \text{\normalfont N}(\frac{1}{\sqrt{n}})$ for $1 \leq i,j \leq n$ and $k = 1,2$. If $\widetilde{A} = A + \gamma G_1$, $\widetilde{B} = B + \gamma G_2$ (again in finite precision) and further
    $$ M = \text{\normalfont {\bf MM}}(\text{\normalfont {\bf INV}}(\widetilde{B}), \widetilde{A}),$$
     then $\Lambda_{\epsilon}(n^{-\alpha}M)$ is shattered with respect to $g$ for 
    $$ \epsilon = \frac{\gamma^5}{32} n^{- \frac{11 \alpha + 25}{3}} $$
    with probability at least $\left[ 1 - \frac{82}{n} - \frac{531441}{16n^2} \right] \left[ 1 - \frac{n^{2 - 2\alpha}}{\gamma^2} - 4e^{-n} \right]$ provided
    $$ {\bf u} \leq \frac{1}{21 \mu_{\text{\normalfont INV}}(n)} \cdot \frac{1}{(6n^{\alpha} + 1)^{c_{\text{\normalfont INV}} \log(n)}} \cdot \frac{\gamma^5}{64} n^{- \frac{11 \alpha + 25}{3}} .$$
\end{thm}
\begin{proof}
    Set $(A_2,B_2) = (\widetilde{A},\widetilde{B})$ and let $(A_1,B_1)$ be the corresponding exact-arithmetic pencil (perturbed with true Ginibre matrices). Further let $X = n^{-\alpha}B_1^{-1}A_1$ and $X' = n^{-\alpha} B_2^{-1}A_2$. Here, $X$ is the exact-arithmetic product covered by \cref{prop: shattering_for_prod} while $X'$ is the exact product corresponding to $(A_2,B_2)$, equivalently a floating-point version of $X$ that assumes exact inversion and matrix multiplication. Throughout, we assume access to the events that guarantee shattering in \cref{prop: shattering_for_prod}, in particular $\sigma_n(B_1) \geq n^{-\alpha}$, and $||A_1||_2,||B_1||_2 \leq 3$, which occur with probability at least $\left[ 1 - \frac{82}{n} - \frac{531441}{16n^2} \right] \left[ 1 - \frac{n^{2 - 2\alpha}}{\gamma^2} - 4e^{-n} \right]$. \\
    \indent As in the proof of \cref{thm: finite_shattering}, we start by noting $||A_1 - A_2||_2, ||B_1 - B_2||_2 \leq (3 + c_{\text{N}})\sqrt{n}{\bf u}$. Since ${\bf u} \leq \frac{1}{2(3 + c_{\text{N}})n^{\alpha + \frac{1}{2}}}$ and therefore $(3 + c_{\text{N}})\sqrt{n} {\bf u} \leq \frac{1}{2n^{\alpha}} < \sigma_n(B_1)$, \cref{lem: product_error} implies
    \begin{equation}
        ||X - X'||_2 \leq 2(3 + c_{\text{N}})\sqrt{n} {\bf u} (1 + 3n^{\alpha}).
        \label{eqn: finite_product_piece_one}
    \end{equation}
    With this in mind, we next seek a bound on $||n^{-\alpha}M - X'||_2$. To do this, let $C = \textbf{INV}(B_2)$. Applying \cref{INV assumption}, we have
    \begin{equation}
        ||C - B_2^{-1}||_2 \leq \mu_{\text{INV}}(n) {\bf u} \kappa_2(B_2)^{c_{\text{INV}}\log(n)} ||B_2^{-1}||_2.
        \label{eqn: C_error}
    \end{equation}
    By the stability of singular values, $||B_2||_2 \leq 3 + (3 + c_{\text{N}})\sqrt{n}{\bf u}$ and $\sigma_n(B_2) \geq n^{-\alpha} - (3 + c_{\text{N}})\sqrt{n} {\bf u}$, so this can be simplified to
    \begin{equation}
        \aligned 
        ||C - B_2^{-1}||_2 &\leq \mu_{\text{INV}}(n) {\bf u} \left( \frac{(3 + (3 + c_{\text{N}})\sqrt{n}{\bf u})n^{\alpha}}{1 - [(3 + c_{\text{N}})\sqrt{n} {\bf u}] n^{\alpha}} \right)^{c_{\text{INV}} \log(n)} \cdot \frac{n^{\alpha}}{1 - [(3 + c_{\text{N}})\sqrt{n} {\bf u} ]n^{\alpha}} \\
        &\leq \mu_{\text{INV}}(n) {\bf u}(6 n^{\alpha}+1)^{c_{\text{INV}} \log(n)} (2n^{\alpha}),
        \endaligned 
        \label{eqn: improved_C_bound}
    \end{equation}
    where we again use the fact that $(3 + c_{\text{N}})\sqrt{n} {\bf u} \leq \frac{1}{2n^{\alpha}}$. Now $M = \textbf{MM}(C,A_2)$, so by \cref{MM assumption}
    \begin{equation}
        ||M - CA_2||_2 \leq \mu_{\text{MM}}(n) {\bf u} ||C||_2||A_2||_2.
        \label{eqn: first_M_bound}
    \end{equation}
    Bounding $||C||_2$ via \eqref{eqn: improved_C_bound} as 
    \begin{equation}
        ||C||_2 \leq ||B_2^{-1}||_2 + \mu_{\text{INV}}(n)  {\bf u}(6n^{\alpha} + 1)^{c_{\text{INV}} \log(n)} (2n^{\alpha}) \leq 2n^{\alpha}\left[ 1 + \mu_{\text{INV}}(n) {\bf u}(6n^{\alpha} + 1)^{c_{\text{INV}} \log(n)}  \right]
        \label{eqn: inverse_norm_bound}
    \end{equation}
    and further noting $||A_2||_2 \leq ||A||_1 + (3 + c_{\text{N}})\sqrt{n}{\bf u} \leq 3 + (3 + c_{\text{N}}) \sqrt{n} {\bf u}$, we conclude
    \begin{equation}
        \aligned 
        ||M - CA_2||_2 &\leq \mu_{\text{MM}}(n){\bf u} (6n^{\alpha} + 1) \left[ 1 + \mu_{\text{INV}}(n) {\bf u}(6n^{\alpha} + 1)^{c_{\text{INV}} \log(n)}  \right] \\
         & = \mu_{\text{MM}}(n) {\bf u} (6n^{\alpha}+1) + \mu_{\text{MM}}(n) \mu_{\text{INV}}(n) {\bf u}^2 (6n^{\alpha} + 1)^{c_{\text{INV}} \log(n) + 1}.
         \endaligned 
        \label{eqn: improved_M_bound}
    \end{equation}
    Putting everything together, we have
    \begin{equation}
        \aligned 
        ||n^{-\alpha}M - &X'||_2 = ||n^{-\alpha}M - n^{-\alpha}CA_2 + n^{-\alpha}CA_2 - n^{-\alpha} B_2^{-1}A_2||_2 \\
        & \leq n^{-\alpha} ||M - CA_2||_2 + n^{-\alpha} ||C - B_2^{-1}||_2 ||A_2||_2 \\
        & \leq 7 \left[ \mu_{\text{MM}}(n) {\bf u} + \mu_{\text{MM}}(n) \mu_{\text{INV}}(n) {\bf u}^2 (6n^{\alpha} + 1)^{c_{\text{INV}} \log(n)} + \mu_{\text{INV}}(n) {\bf u} (6n^{\alpha}+1)^{c_{\text{INV}} \log(n)} \right]
        \endaligned 
        \label{eqn: M_vs_Xprime}
    \end{equation}
    after applying $n^{-\alpha}(6n^{\alpha} + 1) \leq 7$ and $||A_2||_2 \leq 3 + (3 + c_{\text{N}})\sqrt{n} {\bf u} \leq \frac{7}{2}$ to simplify constants. The last term in this expression clearly dominates; assuming each piece of \eqref{eqn: M_vs_Xprime} is bounded by the last one, we obtain
    \begin{equation}
        ||n^{-\alpha}M - X'||_2 \leq  21 \mu_{\text{INV}}(n) {\bf u} (6n^{\alpha} + 1)^{c_{\text{INV}} \log(n)}.
        \label{eqn: simplified_M_bound}
    \end{equation}
    Combining this with \eqref{eqn: finite_product_piece_one} yields our final error bound:
    \begin{equation}
        ||n^{-\alpha} M - X||_2 \leq 21 \mu_{\text{INV}}(n) {\bf u}(6n^{\alpha} + 1)^{c_{\text{INV}} \log(n)} + 2(3 + c_{\text{N}})\sqrt{n} {\bf u} (1 + 3n^{\alpha}).
        \label{eqn: final_M_bound}
    \end{equation}
    \indent Since \cref{prop: shattering_for_prod} implies that $\Lambda_{\epsilon}(X)$ is shattered for $\epsilon = \frac{\gamma^5}{16} n^{- \frac{11 \alpha + 25}{3}}$, we ensure shattering for $\Lambda_{\epsilon}(n^{-\alpha} M)$ with $\epsilon = \frac{\gamma^5}{32}n^{- \frac{11 \alpha + 25}{3}}$ by requiring that each piece of \eqref{eqn: final_M_bound} is bounded by $\frac{\gamma^5}{64}n^{- \frac{11 \alpha + 25}{3}}$, which is guaranteed as long as
    \begin{equation}
        {\bf u} \leq \frac{1}{21 \mu_{\text{INV}}(n)} \cdot \frac{1}{(6n^{\alpha} + 1)^{c_{\text{INV}} \log(n)}} \cdot \frac{\gamma^5}{64} n^{- \frac{11 \alpha + 25}{3}}.
        \label{eqn: final_precision_req}
    \end{equation}
    Note that this requirement on ${\bf u}$ satisfies the assumptions used throughout the proof. 
\end{proof}

Comparing \cref{thm: finite_product_shattering} with \cref{thm: finite_shattering} makes clear the practical cost of forming $\widetilde{B}^{-1}\widetilde{A}$. In particular, finite-precision shattering for the pencil requires $O(\log(n/\epsilon))$ bits of precision versus $O(\log(n/\epsilon) + \log^2(n))$ for $n^{-\alpha} \widetilde{B}^{-1} \widetilde{A}$, where the polylog increase in the latter is rooted in the logarithmic stability of \textbf{INV}. This captures rigorously the phenomenon displayed in \cref{fig: zoomed_pseudospectra_comp2}.

\subsection{Two-Matrix GRURV}
In our algorithm, \textbf{EIG} only ever applies {\bf GRURV} to a pair of matrices in one of the following ways:\ ${\bf GRURV}(2, A_1, A_2, -1, 1)$ or ${\bf GRURV}(2, A_1, A_2, 1, -1)$. Under the floating-point setting established above, the first of these proceeds as follows. \\

\noindent ${\bf GRURV}(2, A_1, A_2, -1, 1)$:
\begin{enumerate}
    \item $[U_2, R_2, V] = {\bf RURV}(A_2)$
    \item $X = {\bf MM}(A_1^H, U_2)$
    \item $[U, R_1^H] = {\bf QL}(X)$
\end{enumerate}
\indent There are three sources of error here: {\bf RURV} in step one, matrix multiplication in step two, and QL in step three. The latter two are controlled by \cref{MM assumption} and \cref{QR assumption} respectively. To control error in finite precision {\bf RURV}, we use work of Banks et al.\ \cite[Lemma C.17]{banks2020pseudospectral} (which is itself built on \cref{sampling assumption}). Under the additional assumptions
\begin{equation}
    4 ||A_2||_2 \max \left\{ c_{\text{N}} \mu_{\text{MM}}(n) {\bf u}, c_{\text{N}} \mu_{\text{QR}}(n) {\bf u} \right\} \leq \frac{1}{4} \leq ||A_2||_2 \; \; \text{and} \; \; 1 \leq \min \left\{ \mu_{\text{MM}}(n), \mu_{\text{QR}}(n), c_{\text{N}} \right\} 
    \label{eqn: FA GRURV assumptions}
\end{equation}
this yields the following.
\begin{thm}
Given $A_1,A_2 \in {\mathbb C}^{n \times n}$, let $[U, R_1, R_2, V] = \text{\normalfont \bf GRURV}(2, A_1, A_2,-1,1)$ on a floating-point machine with precision ${\mathbf u}$ satisfying \eqref{eqn: FA GRURV assumptions}. Then there exist $\widetilde{A}_1$, $\widetilde{A}_2 \in {\mathbb C}^{n \times n}$ and unitary $\widetilde{U}$, $\widetilde{V} \in {\mathbb C}^{n \times n}$ such that $\widetilde{A}_1^{-1} \widetilde{A}_2 = \widetilde{U} R_1^{-1}R_2 \widetilde{V}$ and 
\begin{enumerate}
    \item $\widetilde{V}$ is Haar distributed.
    \item $||U - \widetilde{U}||_2 \leq \mu_{\text{\normalfont QR}}(n) {\bf u}$.
    \item $||A_1 - \widetilde{A}_1||_2 \leq \left[ \mu_{\text{\normalfont QR}}(n) {\bf u} + \mu_{\text{\normalfont MM}}(n) {\bf u}(1 + \mu_{\text{\normalfont QR}}(n) {\bf u}) + \mu_{\text{\normalfont QR}}(n) {\bf u} (1 + \mu_{\text{\normalfont MM}}(n){\bf u})(1 + \mu_{\text{\normalfont QR}}(n) {\bf u}) \right] ||A_1||_2 .$
    \item For every $1 > \alpha > 0$ and $t > 2 \sqrt{2} + 1$, the event that both 
    \begin{itemize}
        \item $||V - \widetilde{V}||_2 \leq \frac{8 tn^{3/2}}{\alpha} c_{\text{\normalfont N}} \mu_{\text{\normalfont QR}}(n){\bf u} + \frac{10n^2}{\alpha} {\bf u} $
        \item $||A_2 - \widetilde{A}_2||_2 \leq \left( \frac{9 tn^{3/2}}{\alpha} c_{\text{\normalfont N}} \mu_{\text{\normalfont QR}}(n) {\bf u} + 2 \mu_{\text{\normalfont MM}}(n) {\bf u} + \frac{10n^2}{\alpha} c_{\text{\normalfont N}} {\bf u} \right) ||A_2||_2 $
    \end{itemize}
    occurs with probability at least $1 - 2e \alpha^2 - 2e^{-t^2n}$.
\end{enumerate}
\label{thm: two matrix GRURV}
\end{thm}
\begin{proof} Beginning with finite-precision guarantees for \textbf{RURV}, \cite[Lemma C.17]{banks2020pseudospectral} implies that there exists unitary matrices $\widetilde{U}_2$ and $\widetilde{V}$ and a matrix $\widetilde{A}_2$ such that $\widetilde{A}_2 = \widetilde{U}_2 R_2 \widetilde{V}$, where $\widetilde{V}$ is Haar distributed and $||U_2 - \widetilde{U}_2||_2$, $||V - \widetilde{V}||_2$, and $||A_2 - \widetilde{A}_2||_2$ can all be controlled. In particular, 
\begin{equation}
    ||U_2 - \widetilde{U}_2||_2 \leq \mu_{\text{QR}}(n) {\bf u} .
    \label{eqn: two grurv 1}
\end{equation}
Consequently, \cref{MM assumption} guarantees that in step two above
\begin{equation}
    ||X - A_1^HU_2||_2 \leq \mu_{\text{MM}}(n) {\bf u}||U_2||_2 ||A_1||_2  \leq \mu_{\text{MM}}(n) {\bf u} (1 + \mu_{\text{QR}}(n) {\bf u})||A_1||_2 .
    \label{eqn: two grurv 2}
\end{equation}
Lastly, \cref{QR assumption} implies that there exists a unitary matrix $\widetilde{U}$ and a matrix $\widetilde{X}$ such that $\widetilde{X} = \widetilde{U} R_1^H$ with
\begin{equation}
    ||X - \widetilde{X}||_2 \leq \mu_{\text{QR}}(n) {\bf u} ||X||_2 .
    \label{eqn: two grurv 3}
\end{equation}
With this in mind, let $\widetilde{A}_1 = \widetilde{U}_2 \widetilde{X}^H$. By construction, we have
\begin{equation}
    \widetilde{A}_1^{-1} \widetilde{A}_2 = (\widetilde{U}_2 \widetilde{X}^H)^{-1}\widetilde{U}_2 R_2 \widetilde{V} = (R_1 \widetilde{U}^H)^{-1}R_2 \widetilde{V} = \widetilde{U} R_1^{-1}R_2 \widetilde{V}, 
    \label{eqn: two grurv 4}
\end{equation}
which means $\widetilde{U} R_1^{-1}R_2 \widetilde{V}$ is a generalized rank-revealing factorization of $\widetilde{A}_1^{-1} \widetilde{A}_2$. Now the distance between $A_2$ and $\widetilde{A_2}$ follows directly from \cite[Lemma C.17]{banks2020pseudospectral}, as does the analogous result for $\widetilde{V}$. Additionally, we observe
\begin{equation}
    \aligned
    ||A_1 - \widetilde{A}_1||_2 &= ||A_1 - \widetilde{U}_2 \widetilde{X}^H||_2 \\
    & \leq ||A_1^H \widetilde{U}_2 - \widetilde{X}||_2 \\
    & = ||A_1^H \widetilde{U}_2 - A_1^HU_2 + A_1^HU_2 - X + X - \widetilde{X} ||_2 \\
    & \leq ||A_1^H(\widetilde{U}_2 - U_2)||_2 + ||A_1^HU_2 - X||_2 + ||X - \widetilde{X}||_2 \\
    & \leq ||A_1||_2 \mu_{\text{QR}}(n) {\bf u} + \mu_{\text{MM}}(n){\bf u}(1 + \mu_{\text{QR}}(n) {\bf u})||A_1||_2 + \mu_{\text{QR}}(n) {\bf u}||X||_2 .
    \endaligned
    \label{eqn: two grurv 5}
\end{equation}
Since (\ref{eqn: two grurv 2}) implies
\begin{equation}
    \aligned
    ||X||_2 &\leq \mu_{\text{MM}}(n) {\bf u}(1 + \mu_{\text{QR}}(n) {\bf u})||A_1||_2 + ||A_1^H U_2||_2 \\
    &  \leq  \mu_{\text{MM}}(n) {\bf u}(1 + \mu_{\text{QR}}(n) {\bf u})||A_1||_2 + ||A_1||_2(1 + \mu_{\text{QR}}(n) {\bf u}) \\
    & = (1 + \mu_{\text{MM}}(n){\bf u})(1 + \mu_{\text{QR}}(n){\bf u}) ||A_1||_2,
    \endaligned 
    \label{eqn: two grurv 6}
\end{equation}
this becomes
\begin{equation}
    \resizebox{.9\textwidth}{!}{$||A_1 - \widetilde{A}_1||_2 \leq \left[ \mu_{\text{QR}}(n) {\bf u} + \mu_{\text{MM}}(n) {\bf u}(1 + \mu_{\text{QR}}(n) {\bf u}) + \mu_{\text{QR}}(n) {\bf u} (1 + \mu_{\text{MM}}(n){\bf u})(1 + \mu_{\text{QR}}(n) {\bf u}) \right] ||A_1||_2$}.
    \label{eqn: two grurv 7}
\end{equation}
Finally recalling \eqref{eqn: two grurv 1}, we have all of the listed bounds.
\end{proof}
While \cref{thm: two matrix GRURV} applies explicitly to $\textbf{GRURV}(2, A_1, A_2, -1, 1)$ it also covers $\textbf{GRURV}(2, A_1, A_2, 1, -1)$; the latter only swaps \textbf{RURV} and \textbf{QL} for \textbf{RULV} and \textbf{QR}, thereby providing the same guarantees since \textbf{QR} and \textbf{QL} satisfy essentially the same assumptions.

\bibliographystyle{abbrv}{}
\bibliography{bib}
\end{document}